\documentclass[12pt,english,reqno]{amsart}
\usepackage{amsmath, amsthm, amssymb}

\advance \topmargin by -\headheight
\advance \topmargin by -\headsep
          
\evensidemargin \oddsidemargin
\theoremstyle{definition}

\theoremstyle{remark}

\numberwithin{equation}{section}

\numberwithin{equation}{section} 
\numberwithin{figure}{section} 
\theoremstyle{plain}
\newtheorem{thm}{Theorem}[section]
  \theoremstyle{plain}
  
  \theoremstyle{plain}
  \newtheorem{lem}[thm]{Lemma}
  \theoremstyle{plain}
  \newtheorem{cor}[thm]{Corollary}

\numberwithin{equation}{section}

\usepackage{babel}
\newcommand{\wt}[1]{\widetilde{#1}}
\newcommand{\sz}[1]{\langle\negthinspace\langle#1\rangle\negthinspace\rangle}
\newcommand{\wh}{\widehat}
\newcommand{\md}[1]{(\text{mod}\;#1)}
\newcommand{\dagg}{\sideset{}{^{\dagger}}}

\newcommand{\mi}[1]{\min\left\{#1\right\}}

\newcommand{\mc}[1]{\mathcal{#1}}

\newcommand{\floor}[1]{\lfloor#1\rfloor}

\newcommand{\Sing}{\text{Sing}\;}

\newcommand{\size}[1]{\langle#1\rangle}
\newcommand{\mb}[1]{\mathbb{#1}}
\newcommand{\ma}[1]{\max\left\{#1\right\}}
\newcommand{\mf}[1]{\mathfrak{#1}}
\newcommand{\mr}[1]{\mathrm{#1}}
\newcommand{\bs}[1]{\boldsymbol{#1}}
\newcommand{\ord}{\text{ord}\;}
\newcommand{\rk}[1]{\text{rank}\;#1}

\newcommand{\res}{\text{res}\;}

\newcommand{\vol}[1]{\text{Vol}\;#1}

\begin{document}
\title{Birch's theorem in function fields}
\author{Siu-lun Alan Lee}
\address{School of Mathematics, University of Bristol, University Walk,
 Bristol BS8 1TW, United Kingdom}
\email{sl5701@bris.ac.uk}
\begin{abstract}
We establish an aysmptotic formula for the number of points with coordinates in $\mb{F}_q[t]$ on a complete intersection of degree $d$ defined over $\mb{F}_q[t]$, with explicit error term, provided that the characteristic of $\mb{F}_q$ is greater than $d$, the codimension of the singular locus of the complete intersection is large enough, and this intersection has a non-singular point at each place of $\mb{F}_q[t]$. In particular, when this complete intersection is non-singular, we show that it satisfies weak approximation.
\end{abstract}
\thanks{The author is supported by the University of Bristol Overseas Centenary
Postgraduate Research Scholarship.}
\maketitle
\section{Introduction}
In his seminal paper \cite{birch}, Birch considered the problem of counting points with integer coordinates on the complete intersection of a number of hypersurfaces defined over $\mb{Z}$. He established an asymptotic formula for the number of such points in an expanding box, under some conditions on the geometry of the complete intersection. On the basis of the close resemblance between the rings $\mb{Z}$ and $\mb{F}_q[t]$, we consider the analogue of Birch's problem over $\mb{F}_q[t]$. To state our problem more precisely, let $s,d,R$ be positive integers with $s>R$ and $d\geq2$, and let $q$ be a power of a prime $p$, so that the field $\mb{F}_q$ has  characteristic $p$. Write $\bs{x}=(x_1,...,x_s)$, and let $X$ denote the complete intersection defined by the simultaneous equations
\begin{equation}\label{eq:X}
F_1(\bs{x})=...=F_R(\bs{x})=0,\end{equation}
where $F_1,...,F_R$ are forms of degree $d$ in $s$ variables with coefficients in $\mb{F}_q[t]$. The aim of this paper is to generalise Birch's result to the function field setting, as well as to establish weak approximation on $X$. 
\par Before we can state our results precisely, we have to introduce some further notation. For simplicity, write $\mb{A}=\mb{F}_q[t]$ and $\mb{K}=\mb{F}_q(t)$. For any non-zero $x\in\mb{A}$, write $\ord x$ for its degree. By convention, we define $\ord0=-\infty$. We extend the function $\ord$to $\mb{K}$ by letting $\ord(x/y)=\ord x-\ord y$ whenever $x,y\in\mb{A}$ with $y\neq0$. We then define an absolute value $\size{\cdot}:\mb{K}\rightarrow\mb{R}$ by taking
\begin{equation}\label{eq:size}
\size{\alpha}=q^{\ord\alpha}\qquad\text{when }\alpha\in\mb{K}.\end{equation}
Let $\mb{K}_{\infty}$ denote the completion of $\mb{K}$ under this absolute value. We extend this absolute value to $\mb{K}_{\infty}^m$, for any positive integer $m$, by
\begin{equation}\label{eq:size of tuple}
\size{\bs{\alpha}}=\ma{\size{\alpha_j}:j=1,...,m}.\end{equation}
As a source of motivation, the sets $\mb{A}$, $\mb{K}$,  $\mb{K}_{\infty}$ are the function field analogues of $\mb{Z}$, $\mb{Q}$ and $\mb{R}$ respectively. With this in mind, we call any tuple with coordinates in $\mb{A}$, $\mb{K}$ and $\mb{K}_{\infty}$ an \emph{integral point}, a \emph{rational point} and a \emph{real point} respectively. It is also convenient to write $\wh{Y}=q^Y$ for any real number $Y$. On the other hand, we need an analogue of the $p$-adic numbers over function fields\footnote{This is the only instance in the paper where $p$ refers to an arbitrary rational prime, and has nothing to do with the characteristic of $\mb{F}_q$.}. In the sequel, we reserve the symbol $\varpi$ for monic irreducible polynomials in $\mb{A}$. For any $\varpi$, let $\mb{K}_{\varpi}$ denote the completion of $\mb{K}$ with respect to the place at $\varpi$. We refer to this completion as the \emph{$\varpi$-adic completion} of $\mb{K}$, and the tuples with coordinates in $\mb{K}_{\varpi}$ as \emph{$\varpi$-adic points}. Finally let $\bs{F}=(F_1,...,F_R)$. For any $s$-tuple $\bs{x}$, write $\nabla\bs{F}(\bs{x})$ for the matrix $(\partial F_i/\partial x_j)$ of derivatives of the forms $F_i$, where $i=1,...,R$ and $j=1,...,s$. Let $\Sing X$ denote the singular locus of $X$. More precisely, $\Sing X$ consists of all points $\bs{x}\in X$, such that
\[\rk\nabla\bs{F}(\bs{x})<R.\]
Any point in $X\backslash\Sing X$ is called a \emph{non-singular point} of $X$. Note that our singular locus differs from Birch's singular locus in \cite{birch} in notation. But Aleksandrov and Moroz (see Corollary of \cite{moroz}) have shown that the two singular loci in fact coincide.
\par By the Lang-Tsen theorem (see Theorem 3.6 of \cite{greenberg}), whenever $s>Rd^2$, there exists a non-trivial rational point on $X$. It is natural now to seek an asymptotic formula for the number of rational points on $X$. In the sequel, let $P$ be a large positive real number that is allowed to tend to infinity. Write $\rho(P)$ for the number of integral points on $X$ with $\size{\bs{x}}<\wh{P}$.
The asymptotic formula for this quantity is unveiled below.
\begin{thm}\label{thm:asymptotic formula}
Write
\begin{equation}\label{eq:assumption on k}
k=s-\dim\mr{Sing}\;X-R(R+1)(d-1)2^{d-1}.\end{equation}
Assume that $k\geq1$, $p>d$, and that $X$ contains a non-singular real point as well as a non-singular $\varpi$-adic point for every $\varpi$. Then, as $P\rightarrow\infty$, we have
\[\rho(P)=\wh{P}^{s-Rd}(\mf{S}\mc{J}+O(\wh{P}^{-k/(2^{d+1}R(d-1))})),\]
where $\mf{S}, \mc{J}$ and the implied constant in the $O(\cdot)$-notation are positive constants depending only on $s,d,R,q$ and $X$.
\end{thm}
The quantities $\mf{S}$ and $\mc{J}$ are defined in \S7 and \S8 respectively. The former, called the \emph{singular series},  is equal to the product of densities of $\varpi$-adic solutions of the equations \eqref{eq:X}, over all $\varpi$. Meanwhile, the \emph{singular integral} $\mc{J}$ is equal to the density of real solutions of the same equations. The number of variables required as asserted in Theorem \ref{thm:asymptotic formula} is identical to the analogous situation, where $\mb{A}$ is replaced by $\mb{Z}$, treated in Theorem 1 of \cite{birch}. Aside from establishing an asymptotic formula with $\mb{A}$ replacing $\mb{Z}$, the major improvement here is the appearance of an explicit error term in our asymptotic formula. Now that we work over a field with positive characteristic $p$, our results inevitably put a condition on $p$. In the theorem above, the condition that $p$ be larger than the degree $d$ of the forms involved is due to the Weyl differencing argument employed in \S3. 
\par For any field $L\supseteq\mb{K}$, write $X(L)=X\cap L^s$. Recall that $X$ satisfies \emph{weak approximation} if the image of the canonical embedding
\[X(\mb{K})\hookrightarrow X(\mb{K}_{\infty})\times\prod_{\varpi}X(\mb{K}_{\varpi})\]
is dense. By adapting the proof of Corollary 1 in \S5 of \cite{skinner} to the function field setting, we prove that the complete intersection $X$ given in \eqref{eq:X} indeed satisfies this property.
\begin{thm}\label{thm:weak approx}
Let $X$ be non-singular. Suppose $p>d$ and 
\[s\geq R(R+1)(d-1)2^{d-1}.\]
Assume further that $X(\mb{K}_{\infty})\neq\phi$ and $X(\mb{K}_{\varpi})\neq\phi$ for every $\varpi$. Then weak approximation holds on $X$.
\end{thm}
As far as the author is concerned, this result does not seem to be a trivial consequence of the Lang-Tsen theorem. In Theorem 1 of \cite{hassett}, Hassett and Tschinkel proved that the function $\phi:\mb{N}\rightarrow\mb{N}$, given by $\phi(1)=1$ and
\[\phi(d)=\binom{\phi(d-1)+d-1}{\phi(d-1)}\qquad(d>1),\]
satisfies the property that, when $\mb{K}$ is replaced by the function field of a smooth curve over any algebraically closed field, $R=1$ and $s>\phi(d)$, the hypersurface $X$ satisfies weak approximation. So when $X$ is a quadratic (resp. cubic) hypersurface, their result requires $3$ (resp. $7$) variables, whereas Theorem \ref{thm:weak approx} requires $4$ (resp. $16$) variables. However, as soon as $d\geq4$, the quantity $\phi(d)$ is much larger than $R(R+1)(d-1)2^{d-1}=2^d(d-1)$.
\par The plan of this paper is as follows. In \S2 the foundations for the use of the Hardy-Littlewood method in the function field setting are put in place, with the view towards relating the counting functions involved to an integral of an exponential sum. In \S3 the technique of Weyl differencing is employed in  establishing an upper bound for our exponential sum that is sufficiently generic with respect to changes in parameters. In \S4 we introduce the key ingredient of the Hardy-Littlewood method, namely the dissection of a suitable compact group into major and minor arcs. In \S5 an estimate for the contribution from the minor arcs to the integral in question is obtained. In the following sections, we turn our attention towards the major arc contribution. In \S6 the exponential sum of interest is expressed on the major arcs in terms of functions that are suitably well behaved, so as to streamline our subsequent analysis. In \S7 and \S8 respectively we investigate the singular series and singular integral. Via an application of the upper bound obtained in \S3, the major arc contribution is expressed in terms of these two quantities, allowing for a permissible error in the process. These two quantities are shown to be positive, so as to render our major arc analysis meaningful. Finally, in \S9, the minor arc estimate in \S5 is combined with the asymptotic formula for the major arc contribution. With different choices of parameters, Theorems \ref{thm:asymptotic formula} and \ref{thm:weak approx} follow at once. Two appendices are included at the end of the paper. They provide essential components for our analysis in the main text. The first deals with the change-of-variable property for integrals over function fields. The second concerns the properties of lattices defined over function fields. 
\par The author would like to thank his PhD supervisor Professor Trevor Wooley for his patient guidance and invaluable comments on his work. Without his help, this paper would never have been possible.
\section{The Hardy-Littlewood method over function fields}
The application of the Hardy-Littlewood method in the setting where integers are replaced by polynomials in $\mb{A}$ requires some explanation. Our goal here is to introduce such notation and basic notions as are subsequently needed to initiate discussion of the key components of this version of the circle method. The material here is taken from section 1 of \cite{kubota} and section 2 of \cite{liuwooley}.
\par By definition, every element of $\mb{K}_{\infty}$ has the form
\begin{equation}\label{eq:form of alpha}
\alpha=\sum_{i=-\infty}^{n}a_it^i,\end{equation}
where $n\in\mb{Z}$ and $a_i\in\mb{F}_q$ for all integers $i\leq n$. Put
\begin{equation}\label{eq:torus}
\mb{T}=\left\{\alpha\in\mb{K}_{\infty}:\size{\alpha}<1\right\}.\end{equation}
If $\alpha\in\mb{K}_{\infty}$ is of the form \eqref{eq:form of alpha}, and $M$ is an integer, we put
\begin{equation}\label{eq:integral part}
\lfloor\alpha\rfloor_M=\sum_{i=M}^{n}a_it^i\end{equation}
and $\left\{\alpha\right\}_M=\alpha-\floor{\alpha}_M$. Evidently $\floor{\alpha}_0\in\mb{A}$ and $\left\{\alpha\right\}_0\in\mb{T}$. It is also convenient to write $\sz{\alpha}=\size{\left\{\alpha\right\}_0}$ for any $\alpha\in\mb{K}_{\infty}$. For any integers $m$ and $M$ with $m$ positive, we extend the operators $\lfloor\cdot\rfloor_M$, $\left\{\cdot\right\}_M$ and $\sz{\cdot}$ to $\mb{K}_{\infty}^m$, by
\[\lfloor\bs{\alpha}\rfloor_M=(\lfloor\alpha_1\rfloor_M,...,\lfloor\alpha_m\rfloor_M),\]
\[\left\{\bs{\alpha}\right\}_M=(\left\{\alpha_1\right\}_M,...,\left\{\alpha_m\right\}_M)\]
and
\begin{equation}\label{eq:sz of tuple}
\sz{\bs{\alpha}}=\ma{\sz{\alpha_j}:j=1,...,m}.\end{equation}
As usual, we write $e(z)$ for $e^{2\pi iz}$ when $z$ is real. We need a similar exponential function on the function field analogue $\mathbb{K}_{\infty}$ of $\mathbb{R}$. For any element $\alpha$ of $\mathbb{K}_{\infty}$ given by \eqref{eq:form of alpha}, define its residue $\res\alpha$ by
\begin{equation}
\res\alpha=\begin{cases}
a_{-1},\qquad & \text{ when }n\geq-1,\\
0, & \text{ when }n<-1.\end{cases}\label{defn:res}\end{equation}
Let $\text{tr}:\mathbb{F}_q\rightarrow\mathbb{F}_p$ be the trace map of the field extension $\mathbb{F}_q$ over $\mathbb{F}_p$. Define the map $e_q:\mathbb{F}_q\rightarrow\mb{C}$ by $e_q(a)=e(\text{tr}(a)/p)$. We then define the exponential function $E:\mathbb{K}_{\infty}\rightarrow\mathbb{C}$ by
\begin{equation}
E(\alpha)=e_q(\res\alpha).\label{defn:e}\end{equation}
Then for any $\gamma\in\mathbb{K}_{\infty}$, the map from $\mb{K}_{\infty}$ to $\mb{C}^{\times}$ given by  $\alpha\mapsto E(\gamma\alpha)$ is an additive character on the additive group $\mathbb{K}_{\infty}$.
\par Unless indicated otherwise, all summations other than those in the shape $\sum_{a=\mu}^{\nu}$, where $\mu,\nu\in\mb{Z}\cup\left\{\pm\infty\right\}$, are over polynomials in $\mb{A}$ or tuples of elements in $\mb{A}$. The set $\mb{T}$ defined in \eqref{eq:torus} is a compact additive subgroup of $\mb{K}_{\infty}$, and thus possesses a unique Haar measure $\mr{d}\alpha$. The Haar measure on $\mb{T}$ is normalised so that $\int_{\mb{T}}\mr{d}\alpha=1$. From Lemma 1(f) of \cite{kubota}, we also have the standard orthogonality property that
\begin{equation}
\int_{\mathbb{T}}E(\alpha x)\;\mathrm{d}\alpha=\begin{cases}
1,\qquad & \text{for }x=0,\\
0,\qquad & \text{for }x\in\mb{A}\backslash\left\{0\right\}.\end{cases}\label{eq:orthogonality}\end{equation}
Hence the set $\mb{T}$ defined in \eqref{eq:torus} can be regarded as the function field analogue of $\mb{R}/\mb{Z}$. Equation (3) of \cite{kubota} provides the useful relation
\begin{equation}
\int_{\size{\alpha}<\wh{m}^{-1}}\;\mr{d}\alpha=\wh{m}^{-1},\qquad\text{for all }m\in\mathbb{N}.\label{eq:order less than -m}\end{equation}
We can extend the Haar measure on $\mb{T}$ to a (unique) translation-invariant measure on $\mb{K}_{\infty}$ by countable additivity. Indeed, let $m$ be a positive integer. The set of all $\alpha\in\mb{K}_{\infty}$ with $\size{\alpha}<\wh{m}$ is the disjoint union of the translates $a_0+a_1t+...+a_{m-1}t^{m-1}+\mb{T}$ of $\mb{T}$, with $a_0,...,a_{m-1}$ ranging over $\mb{F}_q$. Using this together with the normalisation $\int_{\mb{T}}\mr{d}\alpha=1$, we have
\begin{equation}\label{eq:extended measure}
\int_{\size{\alpha}<\wh{m}}\mr{d}\alpha=\sum_{ \size{x}<\wh{m}}\int_{\mb{T}+x}\mr{d}\alpha=\sum_{a_0,...,a_{m-1}\in\mb{F}_q}\int_{\mb{T}}\mr{d}\alpha=\wh{m}.\end{equation}
The measures on $\mb{T}$ and $\mb{K}_{\infty}$ also extend easily to product measures on the respective $D$-fold Cartesian products $\mb{T}^D$ and $\mb{K}_{\infty}^D$, for any positive integer $D$.
\par We examine the weak approximation property by introducing the quantities $\bs{\xi}$, $\bs{b}$, $h$ and $N$. Let $\bs{\xi}\in\mb{T}^s$. Let $h$ be a monic polynomial in $\mb{A}$. Let $\bs{b}\in\mb{A}^s$ with $\size{\bs{b}}<\size{h}$. Take $N$ to be a non-negative integer that is bounded in terms of $\size{h}$. Write $W$ for the affine algebraic set given by the equations
\begin{equation}\label{eq:F(hx+b)=0}
\bs{F}(h\bs{x}+\bs{b})=\bs{0}.\end{equation}
Moreover, let $\mc{B}$ be the hypercube defined by
\begin{equation}\label{eq:B}
\mc{B}=\mc{B}_N=\left\{\bs{\alpha}\in\mb{T}^s:\size{\bs{\alpha}-\bs{\xi}}<\wh{N}^{-1}\right\}.\end{equation}
We call $\wh{N}^{-1}$ the \emph{sidelength} of the hypercube $\mc{B}$. Until \S9, we regard $h$ and $\bs{b}$ as fixed and suppress mention of the dependence of $\mc{B}$ on $N$.
\par Let $n$ be a monic polynomial in $\mb{A}$ with degree $P$. When $m$ is a monic polynomial in $\mb{A}$ and $\mc{E}\subseteq\mb{T}^s$, we write $\bs{x}\in m\mc{E}$ to mean that $m^{-1}\bs{x}\in\mc{E}$. In addition, when this notation appears in a summation, the $\bs{x}$ involved are taken to be integral points in $m\mc{E}$. When $m$ is as above, $\bs{\alpha}\in\mathbb{K}_{\infty}^R$ and  $\mc{E}\subseteq\mb{T}^s$, define the generating function
\begin{equation}
T(\bs{\alpha};m,\mc{E})=\sum_{\bs{x}\in m\mc{E}}E(\bs{\alpha}\cdot\bs{F}(h\bs{x}+\bs{b})).\label{eq:T(alpha)}\end{equation}
In particular, write $T(\bs{\alpha})=T(\bs{\alpha};n,\mc{B})$. For any measurable subset $Z$ of $\mb{K}_{\infty}^R$, write
\begin{equation}\label{eq:rho(P;E)}
\rho_{h,\bs{b}}(n;Z)=\int_{Z}T(\bs{\alpha})\;\mr{d}\bs{\alpha}.\end{equation}
It follows from \eqref{eq:orthogonality} that the number $\rho_{h,\bs{b}}(n)$ of integral solutions $\bs{x}$ of \eqref{eq:F(hx+b)=0} with $\bs{x}\in n\mc{B}$ is given by
\begin{equation}\label{eq:rho(P)}
\rho_{h,\bs{b}}(n)=\rho_{h,\bs{b}}(n;\mb{T}^R).\end{equation}
\par Before our analysis of the quantity $\rho_{h,\bs{b}}(n)$, a few comments on our notation are in order. Recall the Vinogradov notations $O(\cdot)$, $\ll$, $\gg$ and $\asymp$. With the exception of Theorem \ref{thm:asymptotic formula}, the implicit constants in each application of these symbols depend at most on $s$, $d$, $R$, $h$, $\bs{b}$ and the coefficients of the forms $\bs{F}$ of interest. The symbol $\varpi$, with or without subscripts, always denotes a monic irreducible polynomial in $\mb{A}$. For any non-zero  $x\in\mb{A}$, we write $\varpi^l\|x$ to mean that $\varpi^l\mid x$ but $\varpi^{l+1}\nmid x$. The summation notation $\dagg\sum$ is reserved for sums over \emph{monic} polynomials in $\mb{A}$. For any field $L$, we denote its algebraic closure by $\overline{L}$. Given any polynomial $H\in L[x_1,...,x_s]$, the symbol $H^*$ denotes the form given by the sum of the terms of highest degree in $H(x_1,...,x_s)$. For any $R$-tuple $\bs{H}=(H_1,...,H_R)$ of polynomials in $L[x_1,...,x_s]$, the symbol $X_{\bs{H}}$ denotes the affine algebraic set over $\overline{L}$ given by the equations $H_i^*(x_1,...,x_s)=0$, where $i=1,...,R$. For any positive integer $D$, write $I_D$ for the $D$-dimensional identity matrix. Finally, it is worthwhile to comment on the abundance of vector notations for vectors of various dimensions arising in this paper. The interpretation of the dimension of each vector involved depends purely on the context in which it lies.
\section{Weyl's inequality}
The aim of this section is to derive a generic upper bound for exponential sums of the form given in \eqref{eq:T(alpha)}. Throughout this section, we fix a subset $\mc{E}\subseteq\mb{T}^s$ and a monic polynomial $m\in\mb{A}$ with degree $Q$, which we think of as tending to infinity. We also abbreviate the exponential sum $T(\bs{\alpha};m,\mc{E})$ to $U(\bs{\alpha})$, and write $\bs{G}(\bs{x})=\bs{F}(h\bs{x}+\bs{b})$. Then
\begin{equation}\label{eq:U(alpha)}
U(\bs{\alpha})=\sum_{\bs{x}\in m\mc{E}}E(\bs{\alpha}\cdot\bs{F}(h\bs{x}+\bs{b}))=\sum_{\bs{x}\in m\mc{E}}E(\bs{\alpha}\cdot\bs{G}(\bs{x})).\end{equation}
The degree-$d$ part $G_i^*$ of each polynomial $G_i$ can be expressed as
\begin{equation}\label{eq:G_i}
G_i^*(\bs{x})=\sum_{j_1,...,j_d=1}^{s}a_{j_1,...,j_d}^{(i)}x_{j_1}...x_{j_d},\end{equation}
where the coefficients $a_{j_1,...,j_d}^{(i)}$ are in $\mb{A}$ and symmetric with respect to permutations on their subscripts, for each $i=1,...,R$. Without loss of generality, we assume that for each such $i$, the coefficients $a_{j_1,...,j_d}^{(i)}$ are coprime. It is also convenient to introduce their associated multilinear forms
\begin{equation}\label{eq:Psi_j}
\Psi_j^{(i)}(\bs{x}^{(1)},...,\bs{x}^{(d-1)})=d!\sum_{j_1,...,j_{d-1}=1}^{s}a_{j_1,...,j_{d-1},j}^{(i)}x_{j_1}^{(1)}...x_{j_{d-1}}^{(d-1)}
\end{equation}
for $i=1,...,R$ and $j=1,...,s$. We denote by $\bs{\Psi}=\bs{\Psi}(\bs{x}^{(1)},...,\bs{x}^{(d-1)})$ the matrix $(\Psi_j^{(i)})$, with $i=1,...,R$ and $j=1,...,s$.
\par The following lemma, which we do not prove, is the function field analogue of Lemma 2.1 of \cite{birch}.
\begin{lem}\label{lem:weyl differencing}
For every $\bs{\alpha}\in\mb{K}_{\infty}^R$, we have
\[|U(\bs{\alpha})|^{2^{d-1}}\leq\;\wh{Q}^{(2^{d-1}-d)s}\sum_{\size{\bs{x}^{(1)}},...,\size{\bs{x}^{(d-1)}}<\wh{Q}}\prod_{j=1}^{s}|\Upsilon_j|,\]
where $\Upsilon_j=\Upsilon_j(\bs{x}^{(1)},...,\bs{x}^{(d-1)};\bs{\alpha})$ is given by
\begin{equation}\label{eq:Upsilon}
\Upsilon_j=\sum_{\size{x}<\wh{Q}}E\Big(\sum_{i=1}^{R}\alpha_i\Psi_j^{(i)}(\bs{x}^{(1)},...,\bs{x}^{(d-1)})x\Big).\end{equation}
\begin{proof}
The argument used in obtaining the required estimate is standard in circle method literature. We refer the reader to the proof of Lemma 11.1 of \cite{schmidt'} for the underlying argument. 
\end{proof}
\end{lem}
Note that when $p\leq d$, the factor $d!$ forces the multilinear forms $\Psi_{j}^{(i)}$ to vanish identically. Hence, the derivation of a non-trivial upper bound for $U(\bs{\alpha})$ inevitably imposes the condition that $p>d$.
\par Before proceeding further, we quote the following useful lemma.
\begin{lem}\label{lem:linear exp sums}
Let $\gamma\in\mb{K}_{\infty}$, and let $c$ be a non-negative integer. Then
\[\sum_{\size{x}<\wh{c}}E(\gamma x)=\begin{cases}
\wh{c},&\qquad\text{if }\sz{\gamma}<\wh{c}^{-1},\\
0,&\qquad\text{otherwise.}\end{cases}\]
\begin{proof}
This is Lemma 7 of \cite{kubota}.
\end{proof}
\end{lem}
Let $\eta$ be a parameter in the interval $(0,1)$ to be specified. Recall the notation $\sz{\cdot}$ given in \eqref{eq:sz of tuple}.  For any $v=0,1,...,d-1$, let $N_{\eta}^{(v)}(\bs{\alpha})$ denote the number of $(d-1)$-tuples $(\bs{x}^{(1)},...,\bs{x}^{(d-1)})\in(\mb{A}^s)^{d-1}$ with
\begin{equation}\label{eq:size of h}
\size{\bs{x}^{(1)}},...,\size{\bs{x}^{(v)}}<\wh{Q}^{\eta}\qquad\text{and}\qquad\size{\bs{x}^{(v+1)}},...,\size{\bs{x}^{(d-1)}}<\wh{Q},\end{equation}
such that
\begin{equation}\label{eq:size of Phi_j}
\sz{\bs{\alpha}\bs{\Psi}(\bs{x}^{(1)},...,\bs{x}^{(d-1)})}<\wh{Q}^{-v-1+v\eta}.\end{equation}
The following corollary majorises $U(\bs{\alpha})$ in terms of  the quantity $N_{\eta}^{(0)}(\bs{\alpha})$.
\begin{cor}\label{lem:T(alpha) in terms of N(alpha)}
Whenever $p>d$, we have
\[|U(\bs{\alpha})|^{2^{d-1}}\leq\wh{Q}^{(2^{d-1}-d+1)s}N_{\eta}^{(0)}(\bs{\alpha}).\]
\begin{proof}
By the previous lemma, the sum defining $\Upsilon_j$ in \eqref{eq:Upsilon} is equal to $\wh{Q}$ when
\[\sz{\sum_{i=1}^{R}\alpha_i\Psi_{j_d}^{(i)}(\bs{x}^{(1)},...,\bs{x}^{(d-1)})}<\wh{Q}^{-1},\]
and zero otherwise. Using the definition of $N_{\eta}^{(0)}(\bs{\alpha})$, we obtain the required result by putting this back into \eqref{eq:Upsilon}.
\end{proof}
\end{cor}
Our next goal is to provide the tools for bounding the quantity $N_{\eta}^{(0)}(\bs{\alpha})$ in terms of the number $N_{\eta}^{(d-1)}(\bs{\alpha})$. The ultimate aim of this manoeuvre is to choose $\eta$ in such a manner as to further relate integral solutions of \eqref{eq:size of Phi_j}, with $v=d-1$, to singular integral points on $X_{\bs{G}}$. The tools in question are provided by results from the geometry of numbers. We refer the reader to Appendix B for the basic notions and properties of lattices defined over $\mb{K}_{\infty}$.
\par In our later deliberations, we take $\Lambda$ to be the $2s$-dimensional lattice given by the equation $\bs{x}=\Lambda\bs{u}$, where $\bs{u}$ runs over all the elements of $\mb{A}^{2s}$,
\begin{equation}\label{eq:lattice used}
\Lambda=
\begin{pmatrix}
t^{-b}I_s & 0\\
t^b\gamma & t^bI_s\end{pmatrix},\end{equation}
$b$ is an integer and $\gamma=(\gamma_{i,j})_{1\leq i,j\leq s}$ is a symmetric matrix with entries in $\mb{K}_{\infty}$. Let $\mr{M}$ be the adjoint lattice of $\Lambda$, with points $\bs{y}=\mr{M}\bs{v}$, where $\bs{v}$ runs over all the points in $\mb{A}^{2s}$. Then the underlying matrix of $\mr{M}$ is
\[\mr{M}=\begin{pmatrix}
t^bI_s & -t^b\gamma\\
0 & t^{-b}I_s\end{pmatrix}.\]
We can make $\mr{M}$ essentially the same as $\Lambda$ by first changing the signs of $v_{s+1},...,v_{2s}$, then changing those of $y_{s+1},...,y_{2s}$, then interchanging the two sets of variables $v_{1},...,v_s$ and $v_{s+1},...,v_{2s}$, and finally interchanging the two sets of variables $y_1,...,y_s$ and $y_{s+1},...,y_{2s}$. Hence the two lattices $\Lambda$ and $M$ can be regarded as the same lattice. Write $\wh{R}_1,...,\wh{R}_{2s}$ for the successive minima of $\Lambda$. Lemma \ref{lem:adjoint lattices} thus implies that
\begin{equation}\label{eq:same lattice}
\wh{R}_{\nu}\wh{R}_{2s-\nu+1}=1\qquad\text{for }\nu=1,...,s.\end{equation}
\par The following lemma concerns the number of integral solutions of a system of linear inequalities defined over $\mb{K}_{\infty}$.
\begin{lem}\label{lem:comparison}
For $i=1,...,s$, let $L_i(\bs{u})$ be the linear forms in $u_1,...,u_s$ given by the equations
\begin{equation}\label{eq:linear forms}
L_i(\bs{u})=\gamma_{i,1}u_1+...+\gamma_{i,s}u_s,\end{equation}
where $\gamma_{i,j}\in\mb{K}_{\infty}$ with $\gamma_{i,j}=\gamma_{j,i}$ for all $i,j=1,...,s$. When $a,Z\in\mb{R}$, let $M(a;Z)$ denote the number of $u_1,...,u_{2s}\in\mb{A}$ for which
\[\size{u_1},...,\size{u_{s}}<\wh{a}\wh{Z}\]
and
\[\size{L_i(\bs{u})+u_{s+i}}<\wh{Z}/\wh{a}\qquad\text{for all }i=1,...,s.\]
Then whenever $\wh{Z}_1\leq\wh{Z}_2\leq 1$, we have
\[\frac{M(a;Z_1)}{M(a;Z_2)}\gg\Big(\frac{\wh{Z}_1}{\wh{Z}_2}\Big)^s.\]
\begin{proof}
Take $\gamma=(\gamma_{i,j})_{1\leq i,j\leq s}$, and let  $\Lambda$ be the lattice given by \eqref{eq:lattice used} with $b=\lceil a\rceil$. On comparing \eqref{eq:linear forms} with \eqref{eq:lattice used}, and using the inequalities $b-1<a\leq b$, we can bound $M(a;Z)$ above and below by the inequalities
\[\#\left\{\bs{x}\in\Lambda:\size{\bs{x}}<q^{-1}\wh{Z}\right\}\leq M(a;Z)\leq\#\left\{\bs{x}\in\Lambda:\size{\bs{x}}<q\wh{Z}\right\}.\]
Put $R_0=-\infty$, so that $\wh{R}_0=0$. Since $\wh{Z}_1\leq\wh{Z}_2$, we have
\begin{equation}\label{eq:range of Z_1}
\wh{R}_{\mu}\leq q^{-1}\wh{Z}_1<\wh{R}_{\mu+1}\end{equation}
and
\[
\wh{R}_{\nu}\leq q\wh{Z}_2<\wh{R}_{\nu+1},\]
for some $\mu,\nu\in\left\{0,...,2s\right\}$ with $\mu\leq\nu$, as a result of our assumption that $\wh{Z}_1\leq\wh{Z}_2$. By Lemma \ref{lem:lattice points in a cube}, we have
\begin{equation}\label{eq:M(a;Z_1)}
M(a;Z_1)\geq\#\left\{\bs{x}\in\Lambda:\size{\bs{x}}<q^{-1}\wh{Z}_1\right\}=\prod_{\omega=1}^{\mu}(q^{-1}\wh{Z}_1/\wh{R}_{\omega})\end{equation}
and
\begin{equation}\label{eq:M(a;Z_2)}
M(a;Z_2)\leq\#\left\{\bs{x}\in\Lambda:\size{\bs{x}}<q\wh{Z}_2\right\}=\prod_{\omega=1}^{\nu}(q\wh{Z}_2/\wh{R}_{\omega}).\end{equation}
Dividing \eqref{eq:M(a;Z_1)} by \eqref{eq:M(a;Z_2)} gives
\begin{equation}\label{eq:M(a;Z_1)/M(a;Z_2)}
\frac{M(a;Z_1)}{M(a;Z_2)}\geq q^{-\mu-\nu}\Big(\prod_{\omega=\mu+1}^{\nu}\frac{\wh{R}_{\omega}}{\wh{Z}_1}\Big)\Big(\frac{\wh{Z}_1}{\wh{Z}_2}\Big)^{\nu}\gg\Big(\prod_{\omega=\mu+1}^{\nu}\frac{\wh{R}_{\omega}}{\wh{Z}_1}\Big)\Big(\frac{\wh{Z}_1}{\wh{Z}_2}\Big)^{\nu}.\end{equation}
\par Suppose $\nu\leq s$. Using the second inequality in \eqref{eq:range of Z_1}, the assumption $\wh{Z}_1\leq\wh{Z}_2$ and the fact that $\wh{R}_{\mu+1}\leq...\leq\wh{R}_{\nu}$, we get
\begin{align*}
\frac{M(a;Z_1)}{M(a;Z_2)}\gg\;&\prod_{\omega=\mu+1}^{\nu}\frac{\wh{R}_{\omega}}{\wh{R}_{\mu+1}}\Big(\frac{\wh{Z}_1}{\wh{Z}_2}\Big)^s\geq \Big(\frac{\wh{Z}_1}{\wh{Z}_2}\Big)^s.\end{align*}
We have thus proved the lemma in the case where $\nu\leq s$. 
\par Suppose $\mu\leq s<\nu$. We can rewrite \eqref{eq:M(a;Z_1)/M(a;Z_2)} in the form
\begin{equation}\label{eq:Z_1/Z_2 when nu>s}
\frac{M(a;Z_1)}{M(a;Z_2)}\gg\Big(\frac{\wh{Z}_1^{\mu}(\wh{R}_{\mu+1}...\wh{R}_s)}{\wh{Z}_2^s}\Big)\Big(\frac{\wh{R}_{s+1}...\wh{R}_{\nu}}{\wh{Z}_2^{\nu-s}}\Big).\end{equation}
From \eqref{eq:same lattice}, we have $\wh{R}_s\wh{R}_{s+1}=1$. Since $\wh{R}_s\leq\wh{R}_{s+1}$, it follows that $\wh{R}_s\leq1\leq\wh{R}_{s+1}$. Since $\nu>s\geq\mu$, we also have
\begin{equation}\label{eq:order of suc min}
\wh{R}_{\nu}\geq...\geq\wh{R}_{s+1}\geq1\qquad\text{and}\qquad\wh{R}_s\geq...\geq\wh{R}_{\mu+1}\geq\wh{Z}_1.\end{equation}
This, together with \eqref{eq:Z_1/Z_2 when nu>s} and the assumption that $\wh{Z}_2\leq1$, gives rise to
\[\frac{M(a;Z_1)}{M(a;Z_2)}\gg \frac{\wh{Z}_1^{\mu}(\wh{Z}_1^{s-\mu})}{\wh{Z}_2^s}=\Big(\frac{\wh{Z}_1}{\wh{Z}_2}\Big)^s,\]
as required. 
\par Finally let $s<\mu\leq\nu$. We first rewrite \eqref{eq:M(a;Z_1)/M(a;Z_2)} as
\[\frac{M(a;Z_1)}{M(a;Z_2)}\gg\Big(\frac{\wh{Z}_1}{\wh{Z}_2}\Big)^s\frac{\wh{Z}_1^{\mu-s}\wh{R}_{\mu+1}...\wh{R}_{\nu}}{\wh{Z}_2^{\nu-s}}.\]
By the assumption that $\wh{Z}_2\leq1$, this gives
\begin{equation}\label{eq:simplify ratio}
\frac{M(a;Z_1)}{M(a;Z_2)}\gg\Big(\frac{\wh{Z}_1}{\wh{Z}_2}\Big)^s\wh{Z}_1^{\mu-s}\wh{R}_{\mu+1}...\wh{R}_{\nu}.\end{equation}
Since now $\mu>s$, the left inequality in \eqref{eq:range of Z_1} further implies that
\[\wh{Z}_1^{\mu-s}\geq\wh{R}_{\mu}^{\mu-s}\geq\wh{R}_{\mu}...\wh{R}_{s+1}.\]
Also, since $\nu>s$, the first assertion in \eqref{eq:order of suc min} still applies in the present case. Using this information in \eqref{eq:simplify ratio} thus gives
\[\frac{M(a;Z_1)}{M(a;Z_2)}\gg\Big(\frac{\wh{Z}_1}{\wh{Z}_2}\Big)^s\wh{R}_{s+1}...\wh{R}_{\mu}\wh{R}_{\mu+1}...\wh{R}_{\nu}\gg\Big(\frac{\wh{Z}_1}{\wh{Z}_2}\Big)^s,\]
as required.
\end{proof}
\end{lem}
Recall the quantity $N_{\eta}^{(v)}(\bs{\alpha})$ defined before Corollary \ref{lem:T(alpha) in terms of N(alpha)}. The following lemma re-expresses the upper bound in Corollary \ref{lem:T(alpha) in terms of N(alpha)} by majorising the quantity $N_{\eta}^{(0)}(\bs{\alpha})$ in terms of $N_{\eta}^{(d-1)}(\bs{\alpha})$, via successive applications of Lemma \ref{lem:comparison}.
\begin{lem}\label{lem:T(alpha) in terms of N_theta(alpha)}
For any $\bs{\alpha}\in\mb{K}_{\infty}^R$ and $\eta\in[0,1]$, we have
\[|U(\bs{\alpha})|^{2^{d-1}}\ll(\wh{Q}^s)^{2^{d-1}-\eta(d-1)}N_{\eta}^{(d-1)}(\bs{\alpha}).\]
\begin{proof}
Courtesy of Corollary \ref{lem:T(alpha) in terms of N(alpha)}, the conclusion of this lemma holds if we establish for each $\bs{\alpha}\in\mb{K}_{\infty}^R$ that
\[
N_{\eta}^{(d-1)}(\bs{\alpha})\gg\wh{Q}^{s(d-1)(\eta-1)}N_{\eta}^{(0)}(\bs{\alpha}).\]
We first show that for each $v=1,...,d-1$, we have
\begin{equation}\label{eq:Ni in terms of Ni-1}
N_{\eta}^{(v)}(\bs{\alpha})\gg\wh{Q}^{s(\eta-1)}N_{\eta}^{(v-1)}(\bs{\alpha}).\end{equation}
\par Fix such a $v$ together with $\bs{x^{(u)}}\in\mb{A}^s$ ($1\leq u\leq d-1,u\neq v$) which satisfy \eqref{eq:size of h}. We apply Lemma \ref{lem:comparison} with 
\[L_j(\bs{x}^{(v)})=\sum_{i=1}^{R}\alpha_i\Psi_j^{(i)}(\bs{x}^{(1)},...,\bs{x}^{(v)},...,\bs{x}^{(d-1)})\qquad\text{for each }j=1,...,s.\] 
This application is justified as long as these $L_j(\bs{x}^{(v)})$ satisfy the symmetry property as alluded to in that lemma. More precisely, we need to verify that for every $j,j_v=1,...,s$, the coefficient of $x_{j_v}^{(v)}$ in $L_j(\bs{x}^{(v)})$ is unchanged by interchanging $j$ with $j_v$. Indeed, from \eqref{eq:Psi_j}, we see that
\begin{align*}
L_j(\bs{x}^{(v)})=\;&d!\sum_{i=1}^{R}\alpha_i\sum_{j_1,...,j_{d-1}=1}^{s}a_{j_1,...,j_{d-1},j}^{(i)}x_{j_1}^{(1)}...x_{j_{d-1}}^{(d-1)}\\
=\;&d!\sum_{i=1}^{R}\alpha_i\sum_{j_v=1}^{s}x_{j_v}^{(v)}\sum_{\substack{1\leq j_u\leq s\\(1\leq u\leq d-1,\\u\neq v)}}a_{j_1,...,j_{d-1},j}^{(i)}x_{j_1}^{(1)}...x_{j_{v-1}}^{(v-1)}x_{j_{v+1}}^{(v+1)}...x_{j_{d-1}}^{(d-1)}.\end{align*}
The coefficient of $x_{j_v}^{(v)}$ in the expansion for $L_j(\bs{x}^{(v)})$ is thus
\[d!\sum_{i=1}^{R}\alpha_i\sum_{\substack{1\leq j_u\leq s\\(1\leq u\leq d-1,\\u\neq v)}}a_{j_1,...,j_{d-1},j}^{(i)}x_{j_1}^{(1)}...x_{j_{v-1}}^{(v-1)}x_{j_{v+1}}^{(v+1)}...x_{j_{d-1}}^{(d-1)}.\]
Since the coefficients $a_{j_1,...,j_{d-1},j}^{(i)}$ are symmetric with respect to any permutation on their indices, it follows that the coefficient of $x_{j_v}^{(v)}$ in the expansion for $L_j(\bs{x}^{(v)})$ above is symmetric in the indices $j$ and $j_v$, as desired. We can therefore estimate the number $M_1$ of $\bs{x}^{(v)}\in\mb{A}^s$ with $\size{\bs{x}^{(v)}}<\wh{Q}^{\eta}$ which satisfy \eqref{eq:size of Phi_j}, in terms of the number $M_2$ of $\bs{x}^{(v)}\in\mb{A}^s$ with $\size{\bs{x}^{(v)}}<\wh{Q}$ which also satisfy \eqref{eq:size of Phi_j}, but with $v$ replaced by $v-1$ in the last quoted relation. In the notation of that lemma, we have $M_1=M(a;Z_1)$ and $M_2=M(a;Z_2)$, with
\[\wh{a}=\wh{Q}^{(v+1-(v-1)\eta)/2},\;\wh{Z}_1=\wh{Q}^{(v+1)(\eta-1)/2}\text{ and }\wh{Z}_2=\wh{Q}^{(v-1)(\eta-1)/2}.\]
Since $\eta\leq1$, we have $\wh{Z}_1\leq\wh{Z}_2\leq1$, so Lemma \ref{lem:comparison} is applicable. An application of that lemma thus yields
\[M_1\gg\Big(\frac{\wh{Q}^{(v+1)(\eta-1)/2}}{\wh{Q}^{(v-1)(\eta-1)/2}}\Big)^sM_2=\wh{Q}^{s(\eta-1)}M_2.\]
Now summing over all the $\bs{x^{(u)}}$ ($1\leq u\leq d-1,u\neq v$) which satisfy \eqref{eq:size of h}, and recalling the definitions of $M_1$ and $M_2$ above, we obtain \eqref{eq:Ni in terms of Ni-1}. Iterating on \eqref{eq:Ni in terms of Ni-1} thus leads to
\[N_{\eta}^{(d-1)}(\bs{\alpha})\gg\wh{Q}^{s(d-1)(\eta-1)}N_{\eta}^{(0)}(\bs{\alpha}),\]
as required.
\end{proof}
\end{lem}
Our remaining task is to relate the tuples counted by $N_{\eta}^{(d-1)}(\bs{\alpha})$ to the integral singular points on $X_{\bs{G}}$. Henceforth, the notation $\bs{a}/g\in\mb{K}^R$ means that $g$ is a monic polynomial in $\mb{A}$, and $\bs{a}\in\mb{A}^R$ such that $(a_1,...,a_R,g)=1$. We make use of the following lemma, which is a higher-dimensional generalisation of Lemma 2.3 of \cite{heathbrowncubic14} to function fields. It provides conditions under which we can infer from a system of linear inequalities that the underlying matrix does not have full rank.
\begin{lem}\label{lem:inequalities into equations}
Let $M\in\mb{R}$, and $\bs{\alpha}=\bs{a}/g+\bs{\beta}$, where $\bs{a}/g\in\mb{K}^R$ and $\bs{\beta}\in\mb{T}^R$. Suppose $\size{g\bs{\beta}}<\wh{M}^{-1}$, and $\Phi=(\phi_{i,j})$ is an $R\times s$-matrix with entries in $\mb{A}$ such that $\size{\phi_{i,j}}<\wh{M}$ for all $i=1,...,R$ and $j=1,...,s$. Moreover, assume that $\sz{\bs{\alpha}\Phi}<\wh{Y}^{-1}$, where $Y\in\mb{R}$ such that $\wh{Y}>\size{g}$. Then the rank of $\Phi$ (mod $g$) is less than $R$. In addition, if either
\[\size{g}\geq\wh{M}^R\qquad\text{or}\qquad\size{g\bs{\alpha}-\bs{a}}\geq\wh{Y}^{-1}\wh{M}^{R-1},\]
then $\mr{rank}\;\Phi<R$.
\begin{proof}
We first prove the first assertion of the lemma. The assumptions $\size{g\bs{\beta}}<\wh{M}^{-1}$ and $\size{\phi_{i,j}}<\wh{M}$ ($1\leq j\leq s,1\leq i\leq R$) imply that
\begin{align}
\size{\bs{\beta}\Phi}=\;&\max_{1\leq j\leq s}\size{\sum_{i=1}^{R}\beta_i\phi_{i,j}}\leq\max_{1\leq j\leq s}\max_{1\leq i\leq R}\size{\beta_i}\size{\phi_{i,j}}\nonumber\\
<\;&(\size{g}^{-1}\wh{M}^{-1})\wh{M}=\size{g}^{-1}.\label{eq:beta.phi small}\end{align}
Therefore $\bs{\beta}\Phi\in\mb{T}^s$. We then deduce that
\[\left\{\bs{a}\Phi/g\right\}_0=\left\{\bs{\alpha}\Phi-\bs{\beta}\Phi\right\}_0=\left\{\bs{\alpha}\Phi\right\}_0-\bs{\beta}\Phi.\]
Using the assumptions that $\sz{\bs{\alpha}\Phi}<\wh{Y}^{-1}$ and $\size{g}<\wh{Y}$ in conjunction with \eqref{eq:beta.phi small}, we get
\[\sz{\bs{a}\Phi/g}\leq\ma{\sz{\bs{\alpha}\Phi},\size{\bs{\beta}\Phi}}<\ma{\wh{Y}^{-1},\size{g}^{-1}}=\size{g}^{-1}.\]
This necessarily leads to $g\mid\bs{a}\Phi$. Since not all of the $a_i$ are divisible by $g$, it follows from $g|\bs{a}\Phi$ that $\bs{a}$ is a non-trivial solution to the homogeneous linear equation
\[\bs{a}\Phi\equiv\bs{0}\qquad\md{g}.\]
The matrix $\Phi$ (mod $g$) must therefore have less than the full rank $R$. This concludes the proof of the first assertion of the lemma. 
\par We then establish the remaining assertion. Suppose $\rk{\Phi}=R$. Then one of its $R\times R$-determinants does not vanish. To ease the notation, we assume without loss of generality that this non-vanishing determinant is the leading minor of $\Phi$. Let $H$ denote the $R\times R$-matrix formed from the top left corner of $\Phi$. The first assertion of the lemma implies that every $R\times R$-determinant of $\Phi$, and in particular $\det H$, is divisible by $g$ as a result. Also, every entry of $\Phi$ has size less than $\wh{M}$, so $\size{\det H}<\wh{M}^R$. These two observations lead to
\[\size{g}\leq\size{\det H}<\wh{M}^R.\]
Our work in the last paragraph shows that $g\mid\bs{a}\Phi$, and in particular $g$ divides the first $R$ coordinates of the $s$-dimensional vector $\bs{a}\Phi$. But these coordinates form the vector $\bs{a}H$. On recalling that $\bs{\alpha}=\bs{a}/g+\bs{\beta}$, we can write
\[\bs{\alpha}H=\bs{A}+\bs{\Upsilon},\]
where $\bs{A}=\bs{a}H/g\in\mb{A}^R$ and $\bs{\Upsilon}=\bs{\beta}H\in\mb{K}_{\infty}^R$. Hence $g\bs{\alpha}-\bs{a}$ satisfies the linear equation
\[(g\bs{\alpha}-\bs{a})H=g\bs{\Upsilon}.\]
By our assumption that $\det H\neq0$, we can express $g\bs{\alpha}-\bs{a}$ in terms of $H$, $g$ and $\bs{\Upsilon}$ by Cram\'{e}r's rule. Since $g\mid \det H$ and each entry of $H$ has size less than $\wh{M}$, an application of Cram\'{e}r's rule  leads to the bound
\[\size{g\bs{\alpha}-\bs{a}}<\size{g/(\det H)}\size{\bs{\Upsilon}}\wh{M}^{R-1}\leq\size{\bs{\Upsilon}}\wh{M}^{R-1}.\]
Finally, the earlier deductions that $\bs{\beta}\Phi\in\mb{T}^s$ and $g|\bs{a}\Phi$, together with the assumption $\sz{\bs{\alpha}\Phi}<\wh{Y}^{-1}$, implies that 
\[\size{\bs{\Upsilon}}=\size{\bs{\beta}H}\leq\size{\bs{\beta}\Phi}=\sz{\bs{\beta}\Phi}=\sz{\bs{\alpha}\Phi}<\wh{Y}^{-1}.\]
We have therefore shown that, if $\rk{\Phi}=R$, then $\size{g}<\wh{M}^R$ and $\size{g\bs{\alpha}-\bs{a}}<\wh{Y}^{-1}\wh{M}^{R-1}$, contrary to the assumptions made in the latter half of the statement of the lemma.
\end{proof}
\end{lem}
Let $\bs{\alpha}=\bs{a}/g+\bs{\beta}$, where $\bs{a}/g\in\mb{K}^R$, and $\bs{\beta}\in\mb{T}^R$. Our intermediate goal is to bound $N_{\eta}^{(d-1)}(\bs{\alpha})$ above by the number $U_{\eta}$ of $(d-1)$-tuples $(\bs{x}^{(1)},...,\bs{x}^{(d-1)})\in\mb{A}^{s(d-1)}$ satisfying \eqref{eq:size of h} with $v=d-1$, such that
\begin{equation}\label{eq:less than full rank}
\rk{\bs{\Psi}(\bs{x}^{(1)},...,\bs{x}^{(d-1)})}<R.\end{equation}
To this end, we choose $\eta$ carefully so as to render Lemma  \ref{lem:inequalities into equations} applicable. Write $\wh{\kappa}$ for the maximum absolute value of the coefficients of the multilinear forms $\Psi_j^{(i)}$. We then apply that lemma with
\[\wh{M}=\wh{\kappa}\wh{Q}^{(d-1)\eta},\qquad\wh{Y}=\wh{Q}^{d-(d-1)\eta},\]
and $\phi_{i,j}=\Psi_j^{(i)}(\bs{x}^{(1)},...,\bs{x}^{(d-1)})$, where $(\bs{x}^{(1)},...,\bs{x}^{(d-1)})$ varies over those $(d-1)$-tuples that are counted by $N_{\eta}^{(d-1)}(\bs{\alpha})$. An examination of the hypotheses of Lemmata \ref{lem:T(alpha) in terms of N_theta(alpha)} and \ref{lem:inequalities into equations} reveals that they are both applicable if $\eta$ satisfies $\eta\leq1$,
\[\size{g\bs{\beta}}<\wh{\kappa}^{-1}\wh{Q}^{-(d-1)\eta},\qquad\size{g}<\wh{Q}^{d-(d-1)\eta},\]
and
\[\size{g}\geq\mi{\wh{\kappa}^R\wh{Q}^{R(d-1)\eta},\wh{\kappa}^{R-1}\wh{Q}^{R(d-1)\eta}(\wh{Q}^d\size{\bs{\beta}})^{-1}}.\]
The bound given by Lemma \ref{lem:T(alpha) in terms of N_theta(alpha)} is optimised by making $\eta$ as large as possible. We thus choose $\eta$ so that
\begin{align}\label{eq:choice for eta}
&\wh{Q}^{R(d-1)\eta}\nonumber\\
=\;&q^{-1}\mi{\wh{Q}^{R(d-1)},\wh{\kappa}^{-R}\size{g\bs{\beta}}^{-R},\wh{Q}^{Rd}\size{g}^{-R},\size{g}\wh{\kappa}^{-R}\ma{1,\wh{\kappa}\wh{Q}^d\size{\bs{\beta}}}}.\end{align}
An application of Lemma \ref{lem:inequalities into equations}, then yields
\begin{equation}\label{eq:N_eta in terms of U_eta}
N_{\eta}^{(d-1)}(\bs{\alpha})\leq U_{\eta}.\end{equation}
\par Finally we relate the points counted by $U_{\eta}$ to the singular points of the affine algebraic set $X_{\bs{G}}$. Let $\mc{U}$ be the set in affine $s(d-1)$-space given by the relation \eqref{eq:less than full rank}, and let $\mc{D}$ be the diagonal in affine $s(d-1)$-space given by the equations
\[\bs{x}^{(1)}=...=\bs{x}^{(d-1)}.\]
From \eqref{eq:G_i} and \eqref{eq:Psi_j}, we deduce that when $i=1,...,R$ and $j=1,...,s$, we have
\[(d-1)!\partial_jG_i^*(\bs{x})=d!\sum_{j_1,...,j_{d-1}=1}^{s}a_{j_1,...,j_{d-1},j}^{(i)}x_{j_1}...x_{j_{d-1}}=\Psi_j^{(i)}(\bs{x},...,\bs{x}).\]
So when $p>d$, we have $\mc{U}\cap\mc{D}=\Sing X_{\bs{G}}$. Let $\sigma$ temporarily denote the dimension of the singular locus of $X_{\bs{G}}$. By the affine dimension theorem (see Proposition I.7.1 of \cite{hartshorne}), we get
\[\sigma=\dim(\mc{U}\cap\mc{D})\geq\dim\mc{U}+\dim\mc{D}-s(d-1)=\dim\mc{U}-s(d-2).\]
Thus $\dim\mc{U}\leq\sigma+s(d-2)$. The argument which led to Lemma 3.1 of \cite{birch} implies that
\begin{align*}
U_{\eta}=\;&\#\left\{(\bs{x}^{(1)},...,\bs{x}^{(d-1)})\in\mc{U}:\bs{x}^{(u)}\in\mb{A}^s,\;\size{\bs{x}^{(u)}}<\wh{Q}^{\eta}\text{ for }u=1,...,d-1\right\}\\
\ll\;&(\wh{Q}^{\eta})^{\sigma+s(d-2)}.\end{align*}
This, together with \eqref{eq:N_eta in terms of U_eta} and Lemma \ref{lem:T(alpha) in terms of N_theta(alpha)}, leads to the estimate
\begin{align*}
|U(\bs{\alpha})|^{2^{d-1}}\ll\;&(\wh{Q}^s)^{2^{d-1}-\eta(d-1)}U_{\eta}\\
\ll\;&(\wh{Q}^s)^{2^{d-1}-\eta(d-1)}(\wh{Q}^{\eta})^{\sigma+s(d-2)}\\
=\;&\wh{Q}^{2^{d-1}s}(\wh{Q}^{\eta})^{\sigma-s}.\end{align*}
Define the real number $L$ by the equation
\begin{equation}\label{eq:K'}
s-\sigma=2^{d-1}L.\end{equation}
The above inequality can then be rewritten as
\[U(\bs{\alpha})\ll\wh{Q}^{s-L\eta}=\wh{Q}^s(\wh{Q}^{R(d-1)\eta})^{-L/(R(d-1))}.\]
The choice of $\eta$ in \eqref{eq:choice for eta} gives 
\[U(\bs{\alpha})\ll\wh{Q}^s\mi{\wh{Q}^{R(d-1)},\size{g\bs{\beta}}^{-R},\wh{Q}^{Rd}\size{g}^{-R},\size{g}(1+\wh{Q}^d\size{\bs{\beta}})}^{-L/(R(d-1))}.\]
On recalling that $\bs{\beta}=\bs{\alpha}-\bs{a}/g$, we obtain  the following lemma, commonly known as \emph{Weyl's inequality} in circle method literature.
\begin{lem}\label{lem:weyl's inequality}
Let $\mc{E}\subseteq\mb{T}^s$. Set $\bs{G}(\bs{x})=\bs{F}(h\bs{x}+\bs{b})$. Let $L$ be given by \eqref{eq:K'}. Let $\bs{\alpha}\in\mb{K}_{\infty}^R$, $\bs{a}/g\in\mb{K}^R$, and take $m$ to be a monic polynomial in $\mb{A}$ with large degree $Q$. Suppose $p>d$. Then
\begin{align*}
&U(\bs{\alpha})\\
\ll\;&\wh{Q}^s\Big(\frac{\wh{Q}^R+\size{g}^R+\wh{Q}^{Rd}\size{g\bs{\alpha}-\bs{a}}^R}{\wh{Q}^{Rd}}+\frac{1}{\size{g}+\wh{Q}^d\size{g\bs{\alpha}-\bs{a}}}\Big)^{L/(R(d-1))}.\end{align*}
\end{lem}
\section{Dissection into major and minor arcs}
The key step in counting solutions of the simultaneous equations \eqref{eq:X} using the circle method is to dissect the torus $\mb{T}^R$ into \emph{major arcs}, and the complementary set, called the \emph{minor arcs}. Informally speaking, the major arcs consist of those $R$-tuples in $\mb{K}_{\infty}^R$ that are simultaneously well approximated by $R$-tuples of rationals in $\mb{K}^R$ whose common denominator has small degree. The remaining step in our analysis is to show that the contribution from the minor arcs to the underlying integral in \eqref{eq:rho(P)} is negligible compared to that from the major arcs. This is shown in the remaining sections of this paper.
\par We need two slightly different types of major arcs in our work. When $\theta\in(0,1)$, and $\bs{a}/g\in\mb{K}^R$ with $\size{\bs{a}}<\size{g}<\wh{P}^{R(d-1)\theta}$, write
\begin{equation}\label{eq:M(g,a)}
\mf{M}(g,\bs{a};\theta)=\left\{\bs{\alpha}\in\mb{T}^R:\size{g\bs{\alpha}-\bs{a}}<\wh{\Delta}^{-1}\size{h}^{-d}\wh{P}^{-d+R(d-1)\theta}\right\}.\end{equation}
Let $\mf{M}(\theta)$ denote the union of all these sets over the $R$-tuples of rationals $\bs{a}/g$ described above. Also write $\mf{m}(\theta)$ for the complement of $\mf{M}(\theta)$ in $\mb{T}^R$. The sets $\mf{M}(\theta)$, for various values of $\theta$, are used as a tool for the pruning technique employed in obtaining a satisfactory minor arc estimate in the next section. Meanwhile, let $\theta_0$ be a parameter in $(0,1)$ to be specified in the following section. When $\bs{a}/g\in\mb{K}^R$ with $\size{\bs{a}}<\size{g}<\wh{P}^{R(d-1)\theta_0}$, write
\begin{equation}\label{eq:N(g,a)}
\mf{N}(g,\bs{a})=\left\{\bs{\alpha}\in\mb{T}^R:\size{g\bs{\alpha}-\bs{a}}<\wh{\Delta}^{-1}\size{h}^{-d}\size{g}\wh{P}^{-d+R(d-1)\theta_0}\right\}.\end{equation}
Let $\mf{N}$ denote the union of all the sets $\mf{N}(g,\bs{a})$, over the $\bs{a}/g$ described above. Also write $\mf{n}=\mb{T}^R\backslash\mf{N}$. The sets $\mf{N}(g,\bs{a})$ are the major arcs employed in obtaining an asymptotic formula for the major arc contribution to the quantity in \eqref{eq:rho(P)}.
\par With this purpose in mind, we require a mild condition on $\theta_0$ under which the sets $\mf{N}(g,\bs{a})$ are disjoint. This is provided by the following lemma.
\begin{lem}\label{eq:disjointness}
Let $\theta_0$ be a positive number with $\theta_0\leq d/(3R(d-1))$. Then whenever $\bs{a}^{(1)}/g_1,\bs{a}^{(2)}/g_2$ are distinct in $\mb{K}^R$, the sets $\mf{N}(g_1,\bs{a}^{(1)})$ and $\mf{N}(g_2,\bs{a}^{(2)})$ are disjoint.
\begin{proof}
Suppose the conclusion of the lemma were false. Then we could find two distinct $R$-tuples $\bs{a}^{(1)}/g_1,\bs{a}^{(2)}/g_2\in\mb{K}^R$, with $\size{\bs{a}^{(l)}}<\size{g_l}<\wh{P}^{R(d-1)\theta_0}$ for $l=1,2$, such that there would exist $\bs{\alpha}\in\mf{N}(g_1,\bs{a}^{(1)})\cap\mf{N}(g_2,\bs{a}^{(2)})$. The ultrametric property of the absolute value $\size{\cdot}$ on $\mb{K}_{\infty}^R$, together with \eqref{eq:N(g,a)}, would then give
\begin{align*}
\wh{P}^{-2R(d-1)\theta_0}<\size{g_1g_2}^{-1}\leq\;&\size{\bs{a}^{(1)}/g_1-\bs{a}^{(2)}/g_2}\\
\leq\;&\ma{\size{\bs{\alpha}-\bs{a}^{(1)}/g_1},\size{\bs{\alpha}-\bs{a}^{(2)}/g_2}}\\
<\;&\wh{P}^{-d+R(d-1)\theta_0}.\end{align*}
This would lead to $1<\wh{P}^{-d+3R(d-1)\theta_0}$, which is a contradiction under our hypothesis on $\theta_0$.
\end{proof}
\end{lem}
We also need a crude upper bound on the measure of each of the sets $\mf{M}(g,a;\theta)$. From \eqref{eq:M(g,a)} and \eqref{eq:order less than -m}, we deduce that
\begin{align}\label{eq:measure of M(g,a)}
\vol\mf{M}(\theta)\leq\;&\dagg\sum_{\size{g}<\wh{P}^{R(d-1)\theta}}\sum_{\substack{\size{\bs{a}}<\size{g}\\(a_1,...,a_R,g)=1}}\vol\mf{M}(g,\bs{a};\theta)\nonumber\\
\leq\;&\dagg\sum_{\size{g}<\wh{P}^{R(d-1)\theta}}\size{g}^R\prod_{i=1}^{R}\Big(\int_{\size{\beta_i}<\wh{\Delta}^{-1}\size{h}^{-d}\size{g}^{-1}\wh{P}^{-d+R(d-1)\theta}}\mr{d}\beta_i\Big)\nonumber\\
=\;&\dagg\sum_{\size{g}<\wh{P}^{R(d-1)\theta}}\size{g}^R\prod_{i=1}^{R}(\wh{\Delta}^{-1}\size{h}^{-d}\size{g}^{-1}\wh{P}^{-d+R(d-1)\theta})\nonumber\\
\leq\;&\wh{\Delta}^{-R}\size{h}^{-Rd}\wh{P}^{-Rd+R(R+1)(d-1)\theta}.\end{align}
\section{Minor arc estimate}
This section is dedicated to bounding the minor arc contribution $\rho_{h,\bs{b}}(P;\mf{n})$ using the pruning argument employed by Birch in the proof of Lemma 4.4 of \cite{birch}. Let $K$ be the real number defined by the relation
\begin{equation}\label{eq:K}
s-\dim\Sing X=2^{d-1}K.\end{equation}
We need the following function field analogue of the multidimensional Dirichlet approximation lemma.
\begin{lem}\label{lem:dirichlet}
Let $R$ be a positive integer, $Y>0$ and $\bs{\alpha}\in\mb{K}_{\infty}^R$. Then there exists an $R$-tuple $\bs{a}/g\in\mb{K}^R$ with $\size{g}\leq\wh{Y}$ such that $\size{g\bs{\alpha}-\bs{a}}^R<\wh{Y}^{-1}$.
\begin{proof}
Take $N=\lfloor Y/R\rfloor$, so that $NR\leq Y$. The assertion of the lemma then follows once we find a polynomial $g\in\mb{A}$ in the shape
\begin{equation}\label{eq:find g}
g=g_0+g_1t+...+g_{NR}t^{NR},\end{equation}
such that 
\begin{equation}\label{eq:aim of dirichlet}
\sz{g\bs{\alpha_i}}=\size{\left\{g\alpha_i\right\}_0}<\wh{N}^{-1}\qquad\text{ for all }i=1,...,R.\end{equation}
The last condition requires that the coefficients of $t^{-l}$ ($1\leq l\leq N$) in the expansions of $\left\{g\alpha_i\right\}_0$ vanish for all $i$. This gives rise to  exactly $NR$ equations in the $NR+1$ variables $g_0,...,g_{NR}$. Further, a short computation reveals that these equations are actually homogeneous and linear. This guarantees the existence of a non-zero $(NR+1)$-tuple $(g_0,...,g_{NR})$, with the property that the corresponding polynomial $g$ in \eqref{eq:find g} satisfies \eqref{eq:aim of dirichlet}. This completes the proof of this lemma.
\end{proof}
\end{lem}
We then derive an upper bound for the generating function $T(\bs{\alpha})$ using Lemma \ref{lem:weyl's inequality}, in the case when $\bs{\alpha}\notin\mf{M}(\theta)$, for various $\theta$.
\begin{lem}\label{lem:T(alpha) on minor arcs}
Let $\bs{\alpha}\in\mb{K}_{\infty}^R$ and $0<\theta\leq d/((R+1)(d-1))$ such that $\bs{\alpha}\notin\mf{M}(\theta)$. Then
\[T(\bs{\alpha})\ll\wh{P}^{s-K\theta}.\]
\begin{proof}
Observe that for each $i=1,...,R$, the degree-$d$ part of the polynomial $F_i(h\bs{x}+\bs{b})$ in $\bs{x}$ is $F_i(h\bs{x})=h^dF_i(\bs{x})$. On recalling the equations defining $X$ in \eqref{eq:X}, we see that the algebraic set generated by these degree-$d$ parts coincide with $X$. On considering \eqref{eq:K'}, we thus take $L=K$ in applying Lemma \ref{lem:weyl's inequality} to the generating function
\[T(\bs{\alpha})=\sum_{\bs{x}\in n\mc{B}}E(\bs{\alpha}\cdot\bs{F}(h\bs{x}+\bs{b})).\]
\par Let $\bs{\alpha}\notin\mf{M}(\theta)$, and $Y$ be a positive parameter to be chosen. By Lemma \ref{lem:dirichlet}, we can find an $R$-tuple $\bs{a}/g\in\mb{K}^R$ with $\size{g}\leq\wh{Y}$, such that
\[\size{g\bs{\alpha}-\bs{a}}^R<\wh{Y}^{-1}.\]
The assumption on $\bs{\alpha}$, together with the definition \eqref{eq:M(g,a)} of the sets $\mf{M}(g,\bs{a};\theta)$, implies that
\begin{equation}\label{eq:minor arc info}
\size{g}+\wh{P}^d\size{g\bs{\alpha}-\bs{a}}\gg\wh{P}^{R(d-1)\theta}.\end{equation}
Meanwhile, it follows from the information obtained from Lemma \ref{lem:dirichlet} that
\[\wh{P}^R+\size{g}^R+\wh{P}^{Rd}\size{g\bs{\alpha}-\bs{a}}^R\ll\wh{P}^R+\wh{Y}^R+\wh{P}^{Rd}\wh{Y}^{-1}.\]
To optimise this bound, we choose $\wh{Y}=\wh{P}^{Rd/(R+1)}$. This yields
\[\wh{P}^R+\size{g}^R+\wh{P}^{Rd}\size{g\bs{\alpha}-\bs{a}}^R\ll\wh{Y}^R=(\wh{P}^{Rd})^{R/(R+1)}.\]
Applying Lemma \ref{lem:weyl's inequality} with $L=K$, $m=n$ and $Q=P$, and using the last estimate along with \eqref{eq:minor arc info} in the process, we obtain
\begin{align*}
&T(\bs{\alpha})\\
\ll\;&\wh{P}^s\Big(\frac{\wh{P}^R+\size{g}^R+\wh{P}^{Rd}\size{g\bs{\alpha}-\bs{a}}^R}{\wh{P}^{Rd}}+\frac{1}{\size{g}+\wh{P}^d\size{g\bs{\alpha}-\bs{a}}}\Big)^{K/(R(d-1))}\\
\ll\;&\wh{P}^s\Big((\wh{P}^{Rd})^{R/(R+1)-1}+\wh{P}^{-R(d-1)\theta}\Big)^{K/R(d-1)}\\
\ll\;&\wh{P}^{s-K\theta}.\end{align*}
Here the last line follows from our assumption on $\theta$. This completes the proof of the lemma.
\end{proof}
\end{lem}
Evidently the bound given in this lemma is sharpest when $\theta=d/((R+1)(d-1))$. With this in mind, we devise our pruning argument for the minor arc contribution $\rho_{h,\bs{b}}(n;\mf{n})$ as follows. Let $T$ be a sufficiently large integer, and $\theta_0,...,\theta_T$ be positive numbers such that
\begin{equation}\label{eq:thetas}
1/(2R(d-1))=\theta_0<\theta_1<...<\theta_T=d/((R+1)(d-1)),\end{equation}
and
\begin{equation}\label{eq:gaps between thetas}
\theta_{r+1}-\theta_r<(1-R(R+1)(d-1)/K)\theta_0/2\qquad\text{for all }r=0,1,...,T-1.\end{equation}
As elaborated in the next section, the choice for $\theta_0$ in \eqref{eq:thetas} is essentially the largest possible value for which several key results in our subsequent major arc analysis hold. In the meantime, choosing $\theta_0$ to be as large as possible in this respect offers a superior estimate for $\rho_{h,\bs{b}}(n;\mf{n})$, as revealed shortly. The condition \eqref{eq:gaps between thetas} serves a similar purpose on the minor arcs $\mf{n}$, and it is satisfied as long as $T$ is large enough. The plan now is to decompose $\mf{n}$ as the disjoint union of sets
\[\mf{n}=(\mf{M}(\theta_0)\backslash\mf{N})\cup\Big(\bigcup_{r=0}^{T-1}(\mf{M}(\theta_{r+1})\backslash\mf{M}(\theta_r))\Big)\cup\mf{m}(\theta_T).\]
This together with \eqref{eq:rho(P;E)} gives
\begin{align}\label{eq:pruning}
&\rho_{h,\bs{b}}(n;\mf{n})\nonumber\\
=\;&\rho_{h,\bs{b}}(n;\mf{M}(\theta_0)\backslash\mf{N})+\sum_{r=0}^{T-1}\rho_{h,\bs{b}}(n;\mf{M}(\theta_{r+1})\backslash\mf{M}(\theta_r))+\rho_{h,\bs{b}}(n;\mf{m}(\theta_T)).\end{align}
The contribution from each term above can be bounded using \eqref{eq:measure of M(g,a)} and Lemma \ref{lem:T(alpha) on minor arcs} in conjunction. To ensure that the bound obtained on each segment is satisfactory, we assume that
\begin{equation}\label{eq:assumption on K}
K>R(R+1)(d-1),\end{equation}
and choose the $\theta_r$ ($0\leq r\leq T-1$) such that
\begin{equation}\label{eq:delta_r}
\delta_r=K\theta_r-R(R+1)(d-1)\theta_{r+1}>0\qquad\text{for }r=0,...,T-1.\end{equation}
Our choice of $\theta_T$ in \eqref{eq:thetas}, together with \eqref{eq:assumption on K}, implies that the quantity
\begin{equation}\label{eq:delta_T}
\delta_T=K\theta_T-Rd\end{equation}
is positive. Finally write
\begin{equation}\label{eq:delta}
\delta=\min_{0\leq r\leq T}\delta_r.\end{equation}
The desired minor arc estimate is given in the following lemma.
\begin{lem}\label{lem:minor arc estimate}
Under the assumptions \eqref{eq:assumption on K} and $p>d$, we have
\[\rho_{h,\bs{b}}(n;\mf{n})\ll\wh{P}^{s-Rd-\delta},\]
where
\[\delta>(K-R(R+1)(d-1))/(4R(d-1)).\]
\begin{proof}
We bound the contribution from each term in \eqref{eq:pruning} in turn. From \eqref{eq:M(g,a)} and \eqref{eq:N(g,a)}, we see that $\mf{M}(g,\bs{a};\theta_0)\subseteq\mf{N}(g,\bs{a})$ for all $\bs{a}/g\in\mb{K}^R$ for which both sets of arcs are defined. Hence $\mf{M}(\theta_0)\backslash\mf{N}$ is empty, and  $\rho_{h,\bs{b}}(n;\mf{M}(\theta_0)\backslash\mf{N})=0$ as a result. On recalling \eqref{eq:rho(P;E)}, and using \eqref{eq:measure of M(g,a)}, Lemma \ref{lem:T(alpha) on minor arcs} and \eqref{eq:delta_r} in turn, we obtain
\begin{align*}
|\rho_{h,\bs{b}}(n;\mf{M}(\theta_{r+1})\backslash\mf{M}(\theta_r))|\leq\;&\vol\mf{M}(\theta_{r+1})\sup_{\bs{\alpha}\notin\mf{M}(\theta_r)}|T(\bs{\alpha})|\\
\ll\;&\wh{P}^{-Rd+R(R+1)(d-1)\theta_{r+1}}\wh{P}^{s-K\theta_r}\\
=\;&\wh{P}^{s-Rd-\delta_r},\end{align*}
when $0\leq r\leq T-1$. Finally, on applying \eqref{eq:rho(P;E)}, Lemma \ref{lem:T(alpha) on minor arcs} and \eqref{eq:delta_T} together, we get
\[|\rho_{h,\bs{b}}(n;\mf{m}(\theta_T))|\leq\sup_{\bs{\alpha}\in\mf{m}(\theta_T)}|T(\bs{\alpha})|\ll\wh{P}^{s-K\theta_T}=\wh{P}^{s-Rd-\delta_T}.\]
The desired estimate for $\rho_{h,\bs{b}}(n;\mf{n})$ thus follows by substituting these three estimates back into \eqref{eq:pruning}, and recalling \eqref{eq:delta}.
\par It remains to establish the lower bound for $\delta$ asserted by the lemma. Using \eqref{eq:delta_T} together with the choice of $\theta_T$ in \eqref{eq:thetas}, we can rewrite $\delta_T$ as
\[\delta_T=Kd/((R+1)(d-1))-Rd=(K-R(R+1)(d-1))\theta_T.\]
Furthermore, it follows from \eqref{eq:delta_r}, \eqref{eq:thetas} and \eqref{eq:gaps between thetas} that when $0\leq r\leq T-1$, we get
\begin{align*}
\delta_r=\;&-K(\theta_{r+1}-\theta_r)+(K-R(R+1)(d-1))\theta_{r+1}\\
>\;&-K(1-R(R+1)(d-1)/K)\theta_0/2+(K-R(R+1)(d-1))\theta_0\\
=\;&(K-R(R+1)(d-1))\theta_0/2.\end{align*}
On recalling the definition \eqref{eq:delta} of $\delta$ as well as the choice of $\theta_0$ in \eqref{eq:thetas}, we therefore arrive at the lower bound
\[\delta>(K-R(R+1)(d-1))\theta_0/2=(K-R(R+1)(d-1))/(4R(d-1)),\]
as required.
\end{proof}
\end{lem}
\section{Major arc analysis}
Now we turn our attention towards obtaining an asymptotic formula for the major arc contribution $\rho_{h,\bs{b}}(n;\mf{N})$. The standard approach here is to express the generating function $T(\bs{\alpha})$ in terms of analytically well-behaved functions, making explicit use of the assumption that $\bs{\alpha}\in\mf{N}$ in the process. Before we proceed, we introduce the necessary notations for the aforementioned functions. When $g\in\mb{A}$ and $\bs{a}\in\mb{A}^R$, write
\begin{equation}\label{eq:S(g,a)}
S(g,\bs{a})=\sum_{\size{\bs{x}}<\size{g}}E(\bs{a}\cdot\bs{F}(h\bs{x}+\bs{b})/g).\end{equation}
When $\bs{\gamma}\in\mb{K}_{\infty}^R$ and $\mc{E}\subseteq\mb{T}^s$, write
\begin{equation}\label{eq:I}
I(\bs{\gamma};\mc{E})=\int_{\mc{E}}E(\bs{\gamma}\cdot\bs{F}(\bs{\sigma}))\;\mr{d}\bs{\sigma}.\end{equation}
Also write $\wh{\Delta}$ for the maximum absolute value of the coefficients of the forms $F_i$ ($i=1,...,R$).
\par The following lemma gives the desired decomposition of the generating function $U(\bs{\alpha})$, given in \eqref{eq:U(alpha)}, under rather generic assumptions on $\bs{\alpha}$.
\begin{lem}\label{lem:approximate T(alpha)}
Let $\bs{\alpha}\in\mb{K}_{\infty}^R$. Take $\mc{E}$ to be a hypercube in $\mb{T}^s$ with sidelength $\wh{M}^{-1}$, for some non-negative integer $M$. Let $m$ be a monic polynomial in $\mb{A}$ with degree $Q$. Suppose $\bs{a}/g\in\mb{K}^R$ such that $\size{g}<\wh{Q}\wh{M}^{-1}$ and
\begin{equation}\label{eq:approximation condition on beta}
\size{g\bs{\alpha}-\bs{a}}<\wh{\Delta}^{-1}\size{h}^{-d}\wh{Q}^{1-d}.\end{equation}
Then with $\bs{\beta}=\bs{\alpha}-\bs{a}/g$, we have
\[U(\bs{\alpha})=\wh{Q}^s\size{g}^{-s}S(g,\bs{a})I(h^dm^d\bs{\beta};\mc{E}).\]
\begin{proof}
To express $U(\bs{\alpha})$ as asserted in the statement of the lemma, we first examine the summation condition $\bs{x}\in m\mc{E}$ on integral $s$-tuples $\bs{x}$. We express the hypercube $\mc{E}$ as
\begin{equation}\label{eq:hypercube E}
\mc{E}=\left\{\bs{\gamma}\in\mb{T}^s:\size{\bs{\gamma}-\bs{\zeta}}<\wh{M}^{-1}\right\},\end{equation}
for some $\bs{\zeta}\in\mb{T}^s$. The summation condition on $\bs{x}$ is then equivalent to
\[\size{\bs{x}-m\bs{\zeta}}<\wh{Q}\wh{M}^{-1}.\]
In terms of the truncation operation defined in \eqref{eq:integral part}, this can be rewritten further as
\[\lfloor\bs{x}\rfloor_{Q-M}=\bs{\zeta}^*,\]
where
\begin{equation}\label{eq:zeta*}
\bs{\zeta}^*=\lfloor m\bs{\zeta}\rfloor_{Q-M}.\end{equation}
Hence $\bs{x}\in\mb{A}^s\cap m\mc{E}$ if and only if $\bs{x}=\bs{w}+\bs{\zeta}^*$, for some unique $\bs{w}\in\mb{A}^s$ with $\size{\bs{w}}<\wh{Q}\wh{M}^{-1}$. Since $\mb{K}_{\infty}^s$ is equipped with the non-archimedean absolute value given in \eqref{eq:size of tuple}, every point in $\mc{E}$ is a centre of $\mc{E}$. Without loss of generality, we can choose our centre $\bs{\zeta}$ of $\mc{E}$ in such a way that no coordinate of $\mb{\zeta}^*$ equals zero. We can then express $U(\bs{\alpha})$ as
\begin{align}
U(\bs{\alpha})=\;&\sum_{\size{\bs{w}}<\wh{Q}\wh{M}^{-1}}E(\bs{\alpha}\cdot\bs{F}(h(\bs{w}+\bs{\zeta}^*)+\bs{b}))\label{eq:pre-rewrite T(alpha)}\\
=\;&\sum_{\size{\bs{w}}<\wh{Q}\wh{M}^{-1}}E(\bs{\alpha}\cdot\bs{F}(h\bs{w}+\wt{\bs{\zeta}})),\label{eq:rewrite T(alpha)}\end{align}
where
\begin{equation}\label{eq:zeta tilde}
\widetilde{\bs{\zeta}}=h\bs{\zeta}^*+\bs{b}.\end{equation}
\par Since $\wh{Q}\wh{M}^{-1}>\size{g}$ by assumption, we can express each $\bs{w}$ occuring in the sum \eqref{eq:rewrite T(alpha)} uniquely in the shape $\bs{w}=g\bs{y}+\bs{z}$, where $\bs{y},\bs{z}\in\mb{A}^s$ with $\size{\bs{y}}<\wh{Q}\wh{M}^{-1}\size{g}^{-1}$ and $\size{\bs{z}}<\size{g}$. It follows from \eqref{eq:rewrite T(alpha)} that
\begin{align}
&U(\bs{\alpha})\nonumber\\
=\;&\sum_{\size{\bs{z}}<\size{g}}\sum_{\size{\bs{y}}<\wh{Q}\wh{M}^{-1}\size{g}^{-1}}E(\bs{\beta}\cdot\bs{F}(h(g\bs{y}+\bs{z})+\wt{\bs{\zeta}}))E(\bs{a}\cdot\bs{F}(h(g\bs{y}+\bs{z})+\wt{\bs{\zeta}})/g)\nonumber\\
=\;&\sum_{\size{\bs{z}}<\size{g}}E(\bs{a}\cdot\bs{F}(h\bs{z}+\wt{\bs{\zeta}})/g)\sum_{\size{\bs{y}}<\wh{Q}\wh{M}^{-1}\size{g}^{-1}}E(\bs{\beta}\cdot\bs{F}(h(g\bs{y}+\bs{z})+\wt{\bs{\zeta}})).\label{eq:simplify T(alpha)}\end{align}
Here we have used the simple observation that $\bs{F}(h(g\bs{y}+\bs{z})+\wt{\bs{\zeta}})\equiv\bs{F}(h\bs{z}+\wt{\bs{\zeta}})$ (mod $g$).
\par We aim to remove the $\bs{z}$-dependence in the inner sum in \eqref{eq:simplify T(alpha)}. To this end, we seek to establish that
\[E(\bs{\beta}\cdot\bs{F}(h(g\bs{y}+\bs{z})+\wt{\bs{\zeta}}))=E(\bs{\beta}\cdot\bs{F}(hg\bs{y}+\wt{\bs{\zeta}}))\]
whenever $\bs{z},\bs{y}\in\mb{A}^s$ with $\size{\bs{z}}<\size{g}$ and  $\size{\bs{y}}<\wh{Q}\wh{M}^{-1}\size{g}^{-1}$. From \eqref{defn:res} and \eqref{defn:e}, it suffices to show that
\begin{equation}
\ord(\bs{\beta}\cdot(\bs{F}(h(g\bs{y}+\bs{z})+\wt{\bs{\zeta}})-\bs{F}(hg\bs{y}+\wt{\bs{\zeta}})))<-1\label{eq:order}\end{equation}
whenever $\size{\bs{y}}<\wh{Q}\wh{M}^{-1}\size{g}^{-1}$ and $\size{\bs{z}}<\size{g}$. Since each $F_i$ is a form of degree $d$ with coefficients having degree at most  $\Delta$, it follows that
\begin{align}\label{eq:order estimate I}
&\ord(F_i(hg\bs{y}+h\bs{z}+\wt{\bs{\zeta}})-F_i(hg\bs{y}+\wt{\bs{\zeta}}))\nonumber\\
\leq\;&\Delta+\max_{1\leq u\leq s}\ord(hz_u)+(d-1)\max_{1\leq u\leq s}\ord(hgy_u+\wt{\zeta_u}).\end{align}
On recalling \eqref{eq:zeta tilde}, we see that
\[hgy_u+\wt{\zeta_u}=h(gy_u+\zeta_u^*)+b_u.\]
Further, since $\size{gy_u}<\wh{Q}\wh{M}^{-1}\leq\size{\zeta_u^*}$ for all $u=1,...,s$, we see that $gy_u+\zeta_u^*\neq0$ for each $u$. This together with $\size{\bs{b}}<\size{h}$ gives
\[\ord(hgy_u+\wt{\zeta_u})\leq\ord(h(gy_u+\zeta_u^*)).\]
Using the assumption $\size{\bs{y}}<\wh{Q}\wh{M}^{-1}\size{g}^{-1}$, we obtain $\size{g\bs{y}+\bs{\zeta}^*}<\wh{Q}$. It follows that
\begin{equation}\label{eq:order estimate a}
\text{ord}(hgy_u+\wt{\zeta_u})\leq\ord h+(Q-1),\end{equation}
for all $u=1,...,s$. Also, since $\size{\bs{z}}<\size{g}$, we have
\begin{equation}\label{eq:order estimate b}
\ord(hz_u)\leq\ord h+(\ord g-1),\end{equation}
for any such $u$. 
Putting \eqref{eq:order estimate a} and \eqref{eq:order estimate b} into \eqref{eq:order estimate I} thus gives
\begin{align*}
&\ord(F_i(hg\bs{y}+h\bs{z}+\wt{\bs{\zeta}})-F_i(hg\bs{y}+\wt{\bs{\zeta}}))\\
\leq\;&\Delta+\ord h+(\ord g-1)+(d-1)(\ord h+Q-1)\\
=\;&\Delta+d\;\ord h+(d-1)Q+\ord g-d.\end{align*}
The assumption \eqref{eq:approximation condition on beta} also implies that for any such $i$, we have
\[\ord\beta_i\leq-\Delta-d\;\ord h-\ord g-(d-1)Q-1\qquad\text{for all }i=1,...,R.\]
The last two inequalities combined give rise to the upper bound
\[
\ord(\bs{\beta}\cdot(\bs{F}(h(g\bs{y}+\bs{z})+\wt{\bs{\zeta}})-\bs{F}(hg\bs{y}+\wt{\bs{\zeta}})))\leq-d-1<-1.\]
The desired relation \eqref{eq:order} thus follows.
\par Having settled the claim \eqref{eq:order}, it follows from \eqref{eq:simplify T(alpha)} that
\begin{align}\label{eq:factorise T(alpha)}
&U(\bs{\alpha})\nonumber\\
=\;&\Big(\sum_{\size{\bs{z}}<\size{g}}E(\bs{a}\cdot\bs{F}(h\bs{z}+\wt{\bs{\zeta}})/g)\Big)\Big(\sum_{\size{\bs{y}}<\wh{Q}\wh{M}^{-1}\size{g}^{-1}}E(\bs{\beta}\cdot\bs{F}(hg\bs{y}+\wt{\bs{\zeta}}))\Big).\end{align}
Recalling \eqref{eq:zeta tilde}, we can write $h\bs{z}+\wt{\bs{\zeta}}=h\bs{z}+h\bs{\zeta}^*+\bs{b}$. From \eqref{eq:zeta*} and our choice of $\bs{\zeta}$, we see that $\size{\bs{\zeta}^*}\geq\wh{Q}\wh{M}^{-1}>\size{g}$. Hence we can write $\bs{\zeta}^*=g\bs{\eta}+\bs{\phi}$, where $\bs{\eta},\bs{\phi}\in\mb{A}^s$ with $\size{\bs{\phi}}<\size{g}$. Thus
\[h\bs{z}+\wt{\bs{\zeta}}=h\bs{z}+h(g\bs{\eta}+\bs{\phi})+\bs{b}\equiv h(\bs{z}+\bs{\phi})+\bs{b}\;\md{g}.\]
As $\bs{z}$ runs through all residue classes (mod $g$), so does $\bs{z}+\bs{\phi}$. On applying \eqref{eq:S(g,a)}, the first sum in \eqref{eq:factorise T(alpha)} can therefore be simplified as
\begin{equation}\label{eq:first factor in factorisation}
\sum_{\size{\bs{z}}<\size{g}}E(\bs{a}\cdot\bs{F}(h\bs{z}+\wt{\bs{\zeta}})/g)=\sum_{\size{\bs{z}}<\size{g}}E(\bs{a}\cdot\bs{F}(h\bs{z}+\bs{b})/g)
=S(g,\bs{a}).\end{equation}
To simplify the second sum in \eqref{eq:factorise T(alpha)}, we average over all $\bs{z}\in\mb{A}^s$ with $\size{\bs{z}}<\size{g}$, use the claim before \eqref{eq:order} again, put back $\bs{w}=g\bs{y}+\bs{z}$, and use \eqref{eq:zeta tilde}. This gives
\begin{align}\label{eq:T(beta;F,Q)}
&\sum_{\size{\bs{y}}<\wh{Q}\wh{M}^{-1}\size{g}^{-1}}E(\bs{\beta}\cdot\bs{F}(hg\bs{y}+\wt{\bs{\zeta}}))\nonumber\\
=\;&\size{g}^{-s}\sum_{\size{\bs{z}}<\size{g}}\sum_{\size{\bs{y}}<\wh{Q}\wh{M}^{-1}\size{g}^{-1}}E(\bs{\beta}\cdot\bs{F}(h(g\bs{y}+\bs{z})+\wt{\bs{\zeta}}))\nonumber\\
=\;&\size{g}^{-s}\sum_{\size{\bs{w}}<\wh{Q}\wh{M}^{-1}}E(\bs{\beta}\cdot\bs{F}(h\bs{w}+\wt{\bs{\zeta}}))\nonumber\\
=\;&\size{g}^{-s}\sum_{\size{\bs{w}}<\wh{Q}\wh{M}^{-1}}E(\bs{\beta}\cdot\bs{F}(h(\bs{w}+\bs{\zeta}^*)+\bs{b}))\nonumber\\
=\;&\size{g}^{-s}U(\bs{\beta}).\end{align}
Here the last two lines follows from \eqref{eq:rewrite T(alpha)} and \eqref{eq:pre-rewrite T(alpha)} respectively, with $\bs{\alpha}=\bs{\beta}$ in the cited equations.
\par To complete the proof of this lemma, it therefore remains to show that
\begin{equation}\label{eq:sum=integral}
U(\bs{\beta})=\wh{Q}^sI(h^dm^{d}\bs{\beta};\mc{E}).\end{equation}
To achieve this, we use a similar argument to the one employed in establishing \eqref{eq:order}. Our first step is to show that for any $\bs{\sigma}\in\mb{T}^s$ and any integral $s$-tuple $\bs{w}$ with $\size{\bs{w}}<\wh{Q}\wh{M}^{-1}$, we have
\begin{equation}\label{eq:order II}
\size{\bs{\beta}\cdot(\bs{F}(h(\bs{w}+\bs{\zeta}^*+\bs{\sigma})+\bs{b})-\bs{F}(h(\bs{w}+\bs{\zeta}^*)+\bs{b}))}<q^{-1}.\end{equation}
This suffices to give the relation
\[E(\bs{\beta}\cdot\bs{F}(h(\bs{w}+\bs{\zeta}^*+\bs{\sigma})+\bs{b}))=E(\bs{\beta}\cdot\bs{F}(h(\bs{w}+\bs{\zeta}^*)+\bs{b}))\]
for all such $\bs{\sigma}$ and $\bs{w}$. To prove \eqref{eq:order II}, recall that each $F_i$ is a form of degree $d$ with coefficients having absolute value at most $\wh{\Delta}$. This implies that
\[
\size{F_i(h(\bs{w}+\bs{\zeta}^*+\bs{\sigma})+\bs{b})-F_i(h(\bs{w}+\bs{\zeta}^*)+\bs{b})}
\leq\wh{\Delta}\size{h\bs{\sigma}}\size{h(\bs{w}+\bs{\zeta}^*)+\bs{b}}^{d-1}.\]
On recalling \eqref{eq:zeta*} and that $\size{\bs{b}}<\size{h}$ and $0<\size{\bs{w}+\bs{\zeta}^*}<\wh{Q}$, we see that for any $\bs{w}$ and $\bs{\sigma}$ involved in \eqref{eq:order II}, we have the estimate
\begin{align*}
\size{F_i(h(\bs{w}+\bs{\zeta}^*+\bs{\sigma})+\bs{b})-F_i(h(\bs{w}+\bs{\zeta}^*)+\bs{b})}
\leq\;&\wh{\Delta}\size{h\bs{\sigma}}\size{h(\bs{w}+\bs{\zeta}^*)}^{d-1}\\
\leq\;&\wh{\Delta}(\size{h}q^{-1})(\size{h}q^{-1}\wh{Q})^{d-1}\\
=\;&\wh{\Delta}\size{h}^dq^{-d}\wh{Q}^{d-1}.\end{align*}
The condition \eqref{eq:approximation condition on beta} also gives $\size{\beta_i}<\wh{\Delta}^{-1}\size{h}^{-d}\wh{Q}^{1-d}$ for each $i=1,...,R$. On combining these two degree estimates, we get
\begin{align*}
&\size{\beta_i(F_i(h(\bs{w}+\bs{\zeta}^*+\bs{\sigma})+\bs{b})-F_i(h(\bs{w}+\bs{\zeta}^*)+\bs{b}))}\\
\leq\;&(\wh{\Delta}^{-1}\size{h}^{-d}\wh{Q}^{1-d}q^{-1})(\wh{\Delta}\size{h}^dq^{-d}\wh{Q}^{d-1})\\
=\;&q^{-d-1}<q^{-1}.\end{align*}
This yields \eqref{eq:order II}. Integrating both sides of \eqref{eq:pre-rewrite T(alpha)} with respect to all $\bs{\sigma}\in\mb{T}^s$, and using the normalisation on $\mb{T}$ and the relation displayed after \eqref{eq:order II}, we can thus express $U(\bs{\beta})$ as
\[U(\bs{\beta})=\int_{\mb{T}^s}\sum_{\size{\bs{w}}<\wh{Q}\wh{M}^{-1}}E(\bs{\beta}\cdot\bs{F}(h(\bs{w}+\bs{\sigma}+\bs{\zeta}^*)+\bs{b}))\;\mr{d}\bs{\sigma}.\]
As $\bs{w}$ runs through all integral $s$-tuples with $\size{\bs{w}}<\wh{Q}\wh{M}^{-1}$, and $\bs{\sigma}$ runs through all points in $\mb{T}^s$, the $s$-tuple $\bs{\omega}=\bs{w}+\bs{\sigma}$ runs through all real points with $\size{\bs{\omega}}<\wh{Q}\wh{M}^{-1}$. This leads to
\begin{equation}\label{eq:T(beta) as integral}
U(\bs{\beta})=\int_{\size{\bs{\omega}}<\wh{Q}\wh{M}^{-1}}E(\bs{\beta}\cdot\bs{F}(h(\bs{\omega}+\bs{\zeta}^*)+\bs{b}))\;\mr{d}\bs{\omega}.\end{equation}
\par We plan to make the change of variables $\bs{\omega}=m\bs{\sigma}$ in the integral above. To accomplish this using Theorem \ref{thm:change of variables}, we need to verify that the function
\[\bs{\omega}\mapsto E(\bs{\beta}\cdot\bs{F}(h(\bs{\omega}+\bs{\zeta}^*)+\bs{b}))\]
is continuous when $\size{\bs{\omega}}<\wh{Q}\wh{M}^{-1}$. Indeed, note that this function takes only a discrete set of values in $\mb{C}$. Hence, to check the required continuity condition, it suffices to show that
\[E(\bs{\beta}\cdot\bs{F}(h(\bs{\omega}+\bs{\zeta}^*)+\bs{b}))=E(\bs{\beta}\cdot\bs{F}(h(\bs{\omega}'+\bs{\zeta}^*)+\bs{b}))\]
when $\bs{\omega},\bs{\omega}'\in\mb{K}_{\infty}^s$ such that $\size{\bs{\omega}-\bs{\omega}'}$ is sufficiently small. On recalling \eqref{defn:e} and \eqref{defn:res}, we can reduce this task further to showing that
\begin{equation}\label{eq:order estimate}
\size{\bs{\beta}\cdot(\bs{F}(h(\bs{\omega}+\bs{\zeta}^*)+\bs{b})-\bs{F}(h(\bs{\omega}'+\bs{\zeta}^*+\bs{b})))}<q^{-1},\end{equation}
when $\size{\bs{\omega}},\size{\bs{\omega}'}<\wh{Q}\wh{M}^{-1}$ and $\size{\bs{\omega}-\bs{\omega}'}$ is small enough. Observe that
\[(h(\bs{\omega}+\bs{\zeta}^*)+\bs{b})-(h(\bs{\omega}'+\bs{\zeta}^*)+\bs{b})=h(\bs{\omega}-\bs{\omega}').\]
Since each form $F_i$ has degree $d$ with coefficients having absolute value at most $\wh{\Delta}$, it follows from this observation together with the inequalities $\size{\bs{b}}<\size{h}$ and $\size{\bs{\omega}'+\bs{\zeta}^*}>0$ that
\begin{align*}
\size{F_i(h(\bs{\omega}+\bs{\zeta}^*)+\bs{b})-F_i(h(\bs{\omega}'+\bs{\zeta}^*)+\bs{b})}
\leq\;&\wh{\Delta}\size{h(\bs{\omega}'+\bs{\zeta}^*)+\bs{b}}^{d-1}\size{h(\bs{\omega}-\bs{\omega}')}\\
\leq\;&\wh{\Delta}\size{h(\bs{\omega}'+\bs{\zeta}^*)}^{d-1}\size{h(\bs{\omega}-\bs{\omega}')}.\end{align*}
Equation \eqref{eq:zeta*} gives $\size{\bs{\omega}'+\bs{\zeta}^*}<\wh{Q}$. This gives
\begin{align*}
\size{F_i(h(\bs{\omega}+\bs{\zeta}^*)+\bs{b})-F_i(h(\bs{\omega}'+\bs{\zeta}^*)+\bs{b})}\leq\;&\wh{\Delta}(\size{h}^{d-1}(\wh{Q}q^{-1})^{d-1})\size{h}\size{\bs{\omega}-\bs{\omega}'}\\
=\;&\wh{\Delta}\size{h}^d\wh{Q}^{d-1}q^{-(d-1)}\size{\bs{\omega}-\bs{\omega}'}.\end{align*}
Using the assumption \eqref{eq:approximation condition on beta} again, we get
\begin{align*}
&\size{\beta_i(F_i(h(\bs{\omega}+\bs{\zeta}^*)+\bs{b})-F_i(h(\bs{\omega}'+\bs{\zeta}^*)+\bs{b}))}\\
\leq\;&(\wh{\Delta}^{-1}\size{h}^{-d}q^{-1}\wh{Q}^{1-d})(\wh{\Delta}\size{h}^d\wh{Q}^{d-1}q^{-(d-1)}\size{\bs{\omega}-\bs{\omega}'})\\
=\;&q^{-d}\size{\bs{\omega}-\bs{\omega}'},\end{align*}
for each $i=1,...,R$. Hence \eqref{eq:order estimate} follows as long as $\size{\bs{\omega}-\bs{\omega}'}<q^{d-1}$. The integrand in \eqref{eq:T(beta) as integral} is indeed a  continuous function of $\bs{\omega}$, when $\size{\bs{\omega}}<\wh{Q}\wh{M}^{-1}$.
\par Hence we can apply Theorem \ref{thm:change of variables} in \eqref{eq:T(beta) as integral} with the substitution $\bs{\omega}=m\bs{\lambda}$. On recalling \eqref{eq:zeta*} and the assumption that $\size{m}=\wh{Q}$, this yields
\[U(\bs{\beta})=\wh{Q}^s\int_{\size{\bs{\lambda}}<\wh{M}^{-1}}E(\bs{\beta}\cdot\bs{F}(h(m\bs{\lambda}+\lfloor m\bs{\zeta}\rfloor_{Q-M})+\bs{b}))\;\mr{d}\bs{\lambda}.\]
Since $\bs{F}$ is homogeneous of degree $d$, we can rewrite the above as
\begin{equation}\label{eq:extract m}
U(\bs{\beta})=\wh{Q}^s\int_{\bs{\lambda}<\wh{M}^{-1}}E(m^d\bs{\beta}\cdot\bs{F}(h(\bs{\lambda}+\bs{\mu})+\bs{c}))\;\mr{d}\bs{\lambda},\end{equation}
where
\begin{equation}\label{eq:choice of c}
\bs{\mu}=m^{-1}\lfloor m\bs{\zeta}\rfloor_{Q-M}\qquad\text{and}\qquad\bs{c}=m^{-1}\bs{b}.\end{equation}
Note that
\[\bs{\mu}=m^{-1}(m\bs{\zeta}-\left\{m\bs{\zeta}\right\}_{Q-M})=\bs{\zeta}-m^{-1}\left\{m\bs{\zeta}\right\}_{Q-M}.\]
Equation \eqref{eq:hypercube E} and the assumption that  $\size{m}=\wh{Q}$ imply that  $\bs{\lambda}+\bs{\mu}\in\mc{E}$ whenever $\size{\bs{\lambda}}<\wh{M}^{-1}$. Putting $\bs{\sigma}=\bs{\lambda}+\bs{\mu}$ in \eqref{eq:extract m}, and using the translation-invariance of the measure, we obtain
\begin{equation}\label{eq:remove c}
U(\bs{\beta})=\wh{Q}^s\int_{\mc{E}}E(m^d\bs{\beta}\cdot\bs{F}(h\bs{\sigma}+\bs{c}))\;\mr{d}\bs{\sigma}.\end{equation}
\par Our final step in proving \eqref{eq:sum=integral} is to remove the $\bs{c}$ in the integrand above. From \eqref{defn:e} and \eqref{defn:res}, this is valid as long as
\[\size{m^d\bs{\beta}\cdot(\bs{F}(h\bs{\sigma}+\bs{c})-\bs{F}(h\bs{\sigma}))}<q^{-1}\]
for all $\bs{\sigma}\in\mc{E}$. Arguing as in the proof of \eqref{eq:order estimate}, we see that when $i=1,...,R$, we have
\begin{align*}
\size{m^d\beta_i(F_i(h\bs{\sigma}+\bs{c})-F_i(h\bs{\sigma}))}=\;&\size{m}^d\size{\beta_i}\size{F_i(h\bs{\sigma}+\bs{c})-F_i(h\bs{\sigma})}\\
\leq\;&\size{m}^d(\wh{\Delta}^{-1}\wh{Q}^{1-d}\size{h}^{-d})(\wh{\Delta}\size{\bs{c}}\size{h\bs{\sigma}}^{d-1}).\end{align*}
On recalling the definition of $\bs{c}$ in \eqref{eq:choice of c}, together with the assumptions that $\size{m}=\wh{Q}$, $\size{\bs{b}}<\size{h}$ and $\size{\bs{\sigma}}<1$, we get
\begin{align*}
\size{m^d\beta_i(F_i(h\bs{\sigma}+\bs{c})-F_i(h\bs{\sigma}))}\leq\;&\size{m}^d(\wh{\Delta}^{-1}\wh{Q}^{1-d}\size{h}^{-d})(\wh{\Delta}\size{m^{-1}\bs{b}}\size{h}^{d-1}\size{\bs{\sigma}}^{d-1})\\
=\;&(\size{\bs{b}}\size{h}^{-1})\size{\bs{\sigma}}^{d-1}\\
\leq\;&q^{-1}(q^{-1})^{d-1}<q^{-1}\end{align*}
for all $i=1,...,R$, as required. Removing the $\bs{c}$ from the integral in \eqref{eq:remove c}, and once again using the homogeneity of the forms $\bs{F}$ of degree $d$, we therefore obtain
\[
U(\bs{\beta})=\wh{Q}^s\int_{\mc{E}}E(m^d\bs{\beta}\cdot\bs{F}(h\bs{\sigma}))\;\mr{d}\bs{\sigma}=\wh{Q}^s\int_{\mc{E}}E(h^dm^d\bs{\beta}\cdot\bs{F}(\bs{\sigma}))\;\mr{d}\bs{\sigma}.\]
On recalling \eqref{eq:I}, this gives \eqref{eq:sum=integral}, as desired.
\par On substituting \eqref{eq:sum=integral} into \eqref{eq:T(beta;F,Q)}, and putting the resulting equation together with \eqref{eq:first factor in factorisation} into \eqref{eq:factorise T(alpha)}, the equality required by the lemma thus follows.
\end{proof}
\end{lem}
Recall \eqref{eq:N(g,a)}, the choice of $\theta_0$ in \eqref{eq:thetas} and the assumption that $N\ll1$. From these we deduce that when $\bs{\alpha}\in\mf{N}(g,\bs{a})$, for some $\bs{a}/g\in\mb{K}^R$ for which the set $\mf{N}(g,\bs{a})$ is defined, we have
\[\size{g}<\wh{P}^{R(d-1)\theta_0}=\wh{P}^{1/2}<\wh{P}\wh{N}^{-1},\]
as long as $P$ is large enough. Also, equations \eqref{eq:N(g,a)} and \eqref{eq:thetas} give
\begin{align*}
\size{g\bs{\alpha}-\bs{a}}<\;&\wh{\Delta}^{-1}\size{h}^{-d}\size{g}\wh{P}^{-d+R(d-1)\theta_0}<\wh{\Delta}^{-1}\size{h}^{-d}\wh{P}^{-d+2R(d-1)\theta_0}\\
=\;&\wh{\Delta}^{-1}\size{h}^{-d}\wh{P}^{1-d}.\end{align*}
Hence every $\bs{\alpha}\in\mf{N}(g,\bs{a})$ satisfies the conditions of Lemma \ref{lem:approximate T(alpha)} with $\mc{E}=\mc{B}$, $M=N$, $m=n$ and $Q=P$. Using that lemma provides the formula
\begin{equation}\label{eq:T(alpha) on N(g,a)}
T(\bs{\alpha})=\wh{P}^s\size{g}^{-s}S(g,\bs{a})I(h^dn^d\bs{\beta};\mc{B})\qquad\text{for all }\bs{\alpha}\in\mf{N}(g,\bs{a}),\end{equation}
where, as before, $\bs{\beta}=\bs{\alpha}-\bs{a}/g$. 
\par For any monic $g\in\mb{A}$, write
\begin{equation}\label{eq:A(g)}
\mc{A}(g)=\sum_{\substack{\size{\bs{a}}<\size{g}\\(a_1,...,a_R,g)=1}}\size{g}^{-s}S(g,\bs{a}).\end{equation}
When $Y$ is a real number, set
\begin{equation}\label{eq:S(Y)}
\mf{S}(Y)=\dagg\sum_{\size{g}<\wh{Y}}\mc{A}(g)\end{equation}
and
\begin{equation}\label{eq:first J(Y)}
J(Y)=\int_{\size{\bs{\beta}}<\wh{\Delta}^{-1}\size{h}^{-d}\wh{Y}}I(h^dn^d\bs{\beta};\mc{B})\;\mr{d}\bs{\beta}.\end{equation}
Now integrate both sides of \eqref{eq:T(alpha) on N(g,a)} with respect to all $\bs{\alpha}\in\mf{N}(g,\bs{a})$, use the translation-invariance of the measure on $\mb{T}^R$, sum over all $\bs{a}$ and $g$ for which $\mf{N}(g,\bs{a})$ is defined, and finally use Lemma \ref{eq:disjointness} and the definitions \eqref{eq:A(g)}, \eqref{eq:S(Y)} and \eqref{eq:first J(Y)}. By the choice of $\theta_0$ in \eqref{eq:thetas}, we arrive at the formula
\begin{align}\label{eq:approximate asymptotic for rho(P;N)}
\rho_{h,\bs{b}}(n;\mf{N})=\;&\wh{P}^{s}\mf{S}(R(d-1)\theta_0P)J((-d+R(d-1)\theta_0)P)\nonumber\\
=\;&\wh{P}^{s}\mf{S}(P/2)J((-d+1/2)P).\end{align}
\par With the asymptotic formula given in Theorem \ref{thm:asymptotic formula} in mind, the final aim of this section is to extract a factor of $\wh{P}^{-Rd}$ from the quantity $J((-d+1/2)P)$. For any real number $Y$, let
\begin{equation}\label{eq:J(Y)}
\mc{J}(Y)=\int_{\size{\bs{\gamma}}<\wh{\Delta}^{-1}\wh{Y}}I(\bs{\gamma};\mc{B})\;\mr{d}\bs{\gamma}.\end{equation}
We wish to make the substitution $\bs{\gamma}=h^dn^d\bs{\beta}$ in the integral
\begin{equation}\label{eq:final change of variables}
J((-d+1/2)P)=\int_{\size{\bs{\beta}}<\wh{\Delta}^{-1}\size{h}^{-d}\wh{P}^{-d+1/2}}I(h^dn^d\bs{\beta};\mc{B})\;\mr{d}\bs{\beta}.\end{equation}
On recalling that $\size{n}=\wh{P}$, we see that this involves an application of Theorem \ref{thm:change of variables}, which requires that the function $\bs{\gamma}\mapsto I(\bs{\gamma};\mc{B})$ be continuous when $\size{\bs{\gamma}}<\wh{\Delta}^{-1}\wh{P}^{1/2}$. Again from \eqref{defn:e} and \eqref{eq:I}, we see that this function takes a discrete set of values. So our requirement here is satisfied if we have
\[I(\bs{\gamma};\mc{B})=I(\bs{\gamma'};\mc{B})\qquad\text{when }\size{\bs{\gamma}-\bs{\gamma'}}\text{ is small enough.}\]
On recalling the definition \eqref{eq:I} of $I(\bs{\gamma};\mc{B})$, it suffices to show that
\begin{equation}\label{eq:e uniformly continuous II}
E(\bs{\gamma}\cdot\bs{F}(\bs{\sigma}))=E(\bs{\gamma'}\cdot\bs{F}(\bs{\sigma}))\qquad\text{for all }\bs{\sigma}\in\mc{B},\end{equation}
whenever $\size{\bs{\gamma}-\bs{\gamma'}}$ is sufficiently small. From \eqref{defn:e} and \eqref{defn:res}, we see that the displayed equality here holds once we have
\[\size{\bs{\gamma}-\bs{\gamma'}}<(\max_{\bs{\sigma}\in\mc{B}}\size{\bs{F}(\bs{\sigma})}q)^{-1}.\]
From \eqref{eq:B} and \eqref{eq:torus}, the hypercube $\mc{B}$ lies in the compact set $\mb{T}^s$. Thus the maximum above definitely exists. In conclusion, when $\bs{\gamma}$ and $\bs{\gamma'}$ satisfy the above inequality, the desired equality \eqref{eq:e uniformly continuous II} holds, which implies that the function $\bs{\gamma}\mapsto I(\bs{\gamma};\mc{B})$ is indeed continuous when $\size{\bs{\gamma}}<\wh{\Delta}^{-1}\wh{P}^{1/2}$.
\par Making the desired substitution $\bs{\gamma}=h^dn^d\bs{\beta}$ in \eqref{eq:final change of variables}, and using Theorem \ref{thm:change of variables} along with \eqref{eq:J(Y)}, we obtain
\[J((-d+1/2)P)=\size{h}^{-Rd}\wh{P}^{-Rd}\mc{J}(P/2).\]
Putting this back into \eqref{eq:approximate asymptotic for rho(P;N)}, we obtain the simplified formula
\begin{equation}\label{eq:approx asymptotic for rho(P;N)}
\rho_{h,\bs{b}}(n;\mf{N})=\size{h}^{-Rd}\wh{P}^{s-Rd}\mf{S}(P/2)\mc{J}(P/2).\end{equation}
\section{Singular series}
This section is devoted to the investigation of the factor $\mf{S}(P/2)$ that appears in \eqref{eq:approx asymptotic for rho(P;N)}. To make the analysis of this factor more transparent, we extend the $g$-sum which defines it in \eqref{eq:S(Y)} to the infinite sum
\begin{equation}\label{eq:S}
\mf{S}=\dagg\sum_{g\in\mb{A}}\mc{A}(g).\end{equation} 
This quantity is known as the singular series associated to the system \eqref{eq:F(hx+b)=0} of equations. We need to know whether the extension of sums alluded to above is valid, as well as the size of the error incurred in performing this extension. Moreover, to make the resulting asymptotic for $\rho_{h,\bs{b}}(n;\mf{N})$ useful, we also require $\mf{S}$ to be positive. Equations \eqref{eq:S} and \eqref{eq:A(g)} suggest that this information can be obtained by investigating the complete exponential sums $S(g,\bs{a})$. The first two issues can be resolved via an application of Lemma \ref{lem:weyl's inequality}.
On recalling \eqref{eq:T(alpha)} and \eqref{eq:S(g,a)}, note that $S(g,\bs{a})=T(\bs{a}/g;g,\mb{T}^s)$. That lemma therefore yields the upper bound
\[S(g,\bs{a})\ll\size{g}^s\Big(\frac{\size{g}^R}{\size{g}^{Rd}}+\frac{1}{\size{g}}\Big)^{K/(R(d-1))}\ll\size{g}^{s-K/(R(d-1))}.\]
Putting this into \eqref{eq:A(g)}, we therefore obtain
\begin{equation}\label{eq:upper bound for A(g)}
\mc{A}(g)\ll\size{g}^R\size{g}^{-s}\size{g}^{s-K/(R(d-1))}=\size{g}^{R-K/(R(d-1))}.\end{equation}
Armed with this estimate, we are now able to confirm the convergence of $\mf{S}$.
\begin{lem}\label{lem:S converges}
Under the assumptions \eqref{eq:assumption on K} and $p>d$, we have
\[\mf{S}-\mf{S}(Y)\ll\wh{Y}^{R+1-K/(R(d-1))}\qquad\text{for all real }Y.\]
In particular, the singular series $\mf{S}$ converges absolutely.
\begin{proof}
On recalling \eqref{eq:S(Y)} and \eqref{eq:S}, and using \eqref{eq:upper bound for A(g)}, we obtain
\begin{align*}
\mf{S}-\mf{S}(Y)=\;&\dagg\sum_{\size{g}\geq\wh{Y}}\mc{A}(g)\\
\ll\;&\dagg\sum_{\size{g}\geq\wh{Y}}\size{g}^{R-K/(R(d-1))}\\
=\;&\sum_{l\geq Y}(q^l)^{R-K/(R(d-1))}\dagg\sum_{\size{g}=\wh{l}}1\\
=\;&\sum_{l\geq Y}(q^l)^{R+1-K/(R(d-1))}\\
\ll\;&\wh{Y}^{R+1-K/(R(d-1))}.\end{align*}
Here the last line follows from \eqref{eq:assumption on K}. This proves the first assertion of the lemma. The second assertion is an immediate consequence of this and \eqref{eq:assumption on K}.
\end{proof}
\end{lem}
It thus remains to establish the positivity of $\mf{S}$. To this end, we first need the following standard factorisation property of the complete exponential sums $S(g,\bs{a})$.
\begin{lem}\label{lem:factorisation of S(g,a)}
Let $\bs{a}_1/g_1,\bs{a}_2/g_2\in\mb{K}^R$ such that $(g_1,g_2)=1$. Then
\[S(g_1,\bs{a}_1)S(g_2,\bs{a}_2)=S(g_1g_2,g_1\bs{a}_2+g_2\bs{a}_1).\]
\begin{proof}
Since $(g_1,g_2)=1$, every $s$-tuple of residue classes (mod $g_1g_2$) can be written uniquely as $g_1\bs{x}_2+g_2\bs{x}_1$, where $\bs{x}_u$ is a complete residue class (mod $g_u$), for each $u=1,2$. Using this in \eqref{eq:S(g,a)} then yields
\begin{align*}
&S(g_1g_2,g_1\bs{a}_2+g_2\bs{a}_1)\\
=\;&\sum_{\size{\bs{x}_1}<\size{g_1}}\sum_{\size{\bs{x}_2}<\size{g_2}}E\Big(\frac{(g_1\bs{a}_2+g_2\bs{a}_1)\cdot\bs{F}(h(g_1\bs{x}_2+g_2\bs{x}_1)+\bs{b})}{g_1g_2}\Big).\end{align*}
Expanding out each polynomial $F_i(h(g_1\bs{x}_2+g_2\bs{x}_1)+\bs{b})$ ($i=1,...,R$), we  see that
\[F_i(h(g_1\bs{x}_2+g_2\bs{x}_1)+\bs{b})\equiv\begin{cases}
F_i(h(g_2\bs{x}_1)+\bs{b})&\md{g_1}\\
F_i(h(g_1\bs{x}_2)+\bs{b})&\md{g_2}.\end{cases}\]
Hence
\begin{align*}
&S(g_1g_2,g_1\bs{a}_2+g_2\bs{a}_1)\\
=\;&\sum_{\size{\bs{x}_1}<\size{g_1}}\sum_{\size{\bs{x}_2}<\size{g_2}}E(\bs{a}_1\cdot\bs{F}(h(g_2\bs{x}_1)+\bs{b})/g_1)E(\bs{a}_2\cdot\bs{F}(h(g_1\bs{x}_2)+\bs{b})/g_2).\end{align*}
Finally, since $(g_1,g_2)=1$, as $\bs{x}_1$ (resp. $\bs{x}_2$) runs through all $s$-tuples of residue classes (mod $g_1$) (resp. mod $g_2$), so does $g_2\bs{x}_1$ (resp. $g_1\bs{x}_2$). On recalling \eqref{eq:S(g,a)}, we obtain
\begin{align*}
&S(g_1g_2,g_1\bs{a}_2+g_2\bs{a}_1)\\
=\;&\Big(\sum_{\size{\bs{x}_1}<\size{g_1}}E(\bs{a}_1\cdot\bs{F}(h\bs{x}_1+\bs{b})/g_1)\Big)\Big(\sum_{\size{\bs{x}_2}<\size{g_2}}E(\bs{a}_2\cdot\bs{F}(h\bs{x}_2+\bs{b})/g_2)\Big)\\
=\;&S(g_1,\bs{a}_1)S(g_2,\bs{a}_2),\end{align*}
as required.
\end{proof}
\end{lem}
For any $\varpi$, write
\begin{equation}\label{eq:Omega}
\Omega(\varpi)=1+\sum_{\nu=1}^{\infty}\mc{A}(\varpi^{\nu}).\end{equation}
Lemma \ref{lem:factorisation of S(g,a)} is then the ignition spark for the following Euler product representation of $\mf{S}$.
\begin{lem}\label{lem:factorise S}
We can write
\[\mf{S}=\prod_{\varpi}\Omega(\varpi).\]
Further, there exists a constant $Z=Z(s,R,d,\bs{F},h,\bs{b})$ for which
\[\prod_{\size{\varpi}\geq\wh{Z}}\Omega(\varpi)\asymp1.\]
\begin{proof}
On recalling the definition \eqref{eq:A(g)} of the arithmetic function $\mc{A}$, it follows from Lemma \ref{lem:factorisation of S(g,a)} that $\mc{A}$ is multiplicative. The first assertion of this lemma thus follows by using this observation together with \eqref{eq:Omega} on \eqref{eq:S}. To prove the second assertion, note that the upper bound part follows from the first assertion together with Lemma \ref{lem:S converges}. It remains to prove the lower bound part. An application of the upper bound \eqref{eq:upper bound for A(g)} on \eqref{eq:Omega} gives
\[\Omega(\varpi)-1=\sum_{\nu=1}^{\infty}\mc{A}(\varpi^{\nu})\ll\sum_{\nu=1}^{\infty}\size{\varpi^{\nu}}^{R-K/(R(d-1))}\ll\size{\varpi}^{R-K/(R(d-1))}.\]
There thus exists a constant $Z=Z(s,d,R,\bs{F},h,\bs{b})$ for which
\[\sum_{\size{\varpi}\geq\wh{Z}}|\Omega(\varpi)-1|\ll\sum_{\size{\varpi}\geq\wh{Z}}\size{\varpi}^{R-K/(R(d-1))}.\]
By the polynomial prime number theorem (see Theorem 2.2 of \cite{rosen}), the number of monic irreducible polynomials in $\mb{A}$ with degree $l$ is $O(q^l/l)$. It thus follows from the last displayed inequality and \eqref{eq:assumption on K} that
\begin{align*}
\sum_{\size{\varpi}\geq\wh{Z}}|\Omega(\varpi)-1|\ll\;&\sum_{l\geq Z}(q^l)^{R-K/(R(d-1))}\sum_{\size{\varpi}=\wh{l}}1\\
\ll\;&\sum_{l\geq Z}(q^l)^{R+1-K/(R(d-1))}/l\ll1.\end{align*}
On extracting logarithms, we find that
\[\prod_{\size{\varpi}\geq\wh{Z}}\Omega(\varpi)\gg1.\]
The lower bound part in the second assertion of the lemma duly follows.
\end{proof}
\end{lem}
Hence, to establish that $\mf{S}>0$, it suffices to show that $\Omega(\varpi)>0$ for all $\varpi$ with $\size{\varpi}\ll1$. There is a standard interpretation of $\Omega(\varpi)$ as the density of $\varpi$-adic solutions of the simultaneous equations in \eqref{eq:F(hx+b)=0}. To make this precise, we have to introduce more notation. For any monic $g\in\mb{A}$, let $M(g)$ denote the number of solutions of the simultaneous congruences
\begin{equation}\label{eq:M(g)}
\bs{F}(h\bs{x}+\bs{b})\equiv\bs{0}\qquad\md{g}.\end{equation}
We aim to relate $\Omega(\varpi)$ to $M(\varpi^k)$ for large integers $k$. The first step to achieving this is the following lemma concerning orthogonality relations for complete exponential sums.
\begin{lem}\label{lem:orthogonality for complete exp sums}
Let $y,g\in\mb{A}$ with $g\neq0$. Then
\[\sum_{\size{x}<\size{g}}E(xy/g)=\begin{cases}
\size{g},&\qquad\text{if }g\mid y,\\
0,&\qquad\text{otherwise.}\end{cases}\]
\begin{proof}
This is Lemma 1(g) of \cite{kubota}.
\end{proof}
\end{lem}
Our objective, alluded to before the last lemma, is made exact by the following lemma.
\begin{lem}\label{lem:local solutions}
For any $\varpi$ and positive integer $k$, we have
\[M(\varpi^k)=\size{\varpi}^{k(s-R)}\sum_{\nu=0}^{k}\mc{A}(\varpi^{\nu}).\]
In particular, we have
\[\Omega(\varpi)=\lim_{k\rightarrow\infty}\size{\varpi}^{-k(s-R)}M(\varpi^k).\]
\begin{proof}
Applying Lemma \ref{lem:orthogonality for complete exp sums} on the counting function $M(\varpi^k)$, we obtain
\[M(\varpi^k)=\size{\varpi}^{-Rk}\sum_{\size{\bs{u}}<\size{\varpi}^k}\sum_{\size{\bs{x}}<\size{\varpi}^k}E(\bs{u}\cdot\bs{F}(h\bs{x}+\bs{b})/\varpi^k).\]
Now put $\bs{u}=\varpi^{k-\nu}\bs{w}$, where $\nu=0,1,...,k$ and $\bs{w}\in\mb{A}^R$ with the property that $(w_1,...,w_R,\varpi)=1$. This gives
\[M(\varpi^k)=\size{\varpi}^{-Rk}\sum_{\nu=0}^{k}\sum_{\substack{\size{\bs{w}}<\size{\varpi}^{\nu}\\(w_1,...,w_R,\varpi)=1}}\sum_{\size{\bs{x}}<\size{\varpi}^k}E(\bs{w}\cdot\bs{F}(h\bs{x}+\bs{b})/\varpi^{\nu}).\]
The innermost sum here is periodic (mod $\varpi^{\nu}$), hence
\[M(\varpi^k)=\size{\varpi}^{-Rk}\sum_{\nu=0}^{k}\size{\varpi}^{(k-\nu)s}\sum_{\substack{\size{\bs{w}}<\size{\varpi}^{\nu}\\(w_1,...,w_R,\varpi)=1}}\sum_{\size{\bs{x}}<\size{\varpi}^{\nu}}E(\bs{w}\cdot\bs{F}(h\bs{x}+\bs{b})/\varpi^{\nu}).\]
On recalling \eqref{eq:S(g,a)} and \eqref{eq:A(g)}, we obtain
\begin{align*}
M(\varpi^k)=\;&\size{\varpi}^{k(s-R)}\sum_{\nu=0}^{k}\size{\varpi}^{-\nu s}\sum_{\substack{\size{\bs{w}}<\size{\varpi}^{\nu}\\(w_1,...,w_R,\varpi)=1}}S(\varpi^{\nu},\bs{w})\\
=\;&\size{\varpi}^{k(s-R)}\sum_{\nu=0}^{k}\mc{A}(\varpi^{\nu}).\end{align*}
This proves the first assertion of the lemma. The second assertion follows immediately from rearranging this equality, and recalling the definition \eqref{eq:Omega} of $\Omega(\varpi)$.
\end{proof}
\end{lem}
The second assertion of the above lemma implies that, in order to have $\Omega(\varpi)>0$, it suffices to show that
\begin{equation}\label{eq:many local solutions}
M(\varpi^k)\gg\size{\varpi}^{k(s-R)}\qquad\text{for all large }k,\end{equation}
where the implied constant here could depend on $\varpi$, but is independent of $k$. The plan is to show that, as long as the simultaneous equations \eqref{eq:F(hx+b)=0} have a non-singular $\varpi$-adic solution for every $\varpi$, then the desired inequality \eqref{eq:many local solutions} holds. When $\bs{x}\in\mb{K}_{\infty}^s$, write $\mc{F}(\bs{x})$ for the matrix with $(i,j)$-entry $(\partial/\partial x_j)(F_i(h\bs{x}+\bs{b}))$, for any $i=1,...,R$ and $j=1,...,s$. The following lemma brings us one step closer to the desired relation \eqref{eq:many local solutions} via a Hensel lifting argument.
\begin{lem}\label{lem:smooth local solutions}
Suppose that for every $\varpi$, there exists a positive integer $l$ along with an $s$-tuple $\bs{x}$ (mod $\varpi^{2l-1}$) for which
\begin{equation}\label{eq:F(x)=0 mod pi^2l-1}
\bs{F}(h\bs{x}+\bs{b})\equiv\bs{0}\qquad\md{\varpi^{2l-1}}\end{equation}
and
\begin{equation}\label{eq:X most singular mod pi^l-1}
\mr{rank}\;\mc{F}(\bs{x})<R\;\md{\varpi^{l-1}}\;\text{but }\mr{rank}\;\mc{F}(\bs{x})=R\;\md{\varpi^l}.\end{equation}
Then for every positive integer $\nu$, there exists at least $\size{\varpi}^{\nu(s-R)}$ tuples $\bs{y}$ $\md{\varpi^{2l-1+\nu}}$ such that
\begin{equation}\label{eq:F(x)=0 mod pi^2l-1+nu}
\bs{F}(h\bs{y}+\bs{b})\equiv\bs{0}\qquad\md{\varpi^{2l-1+\nu}}\end{equation}
and
\begin{equation}\label{eq:X most singular mod pi^l-1 II}
\mr{rank}\;\mc{F}(\bs{y})<R\;\md{\varpi^{l-1}}\;\text{but }\mr{rank}\;\mc{F}(\bs{y})=R\;\md{\varpi^l}.\end{equation}
\begin{proof}
This proof is based on the argument which led to Hensel's lemma, also used in the proof of Lemma 17.1 of \cite{davenport}. We proceed by induction on $\nu$. When $\nu=1$, let $\bs{v}$ be an $s$-tuple (mod $\varpi$) to be specified. By the multidimensional Taylor's theorem, we have
\[\bs{F}(h(\bs{x}+\varpi^l\bs{v})+\bs{b})\equiv\bs{F}(h\bs{x}+\bs{b})+\nabla\bs{F}(h\bs{x}+\bs{b})\cdot\varpi^lh\bs{v}\;\md{\varpi^{2l}}.\]
On recalling the definition of the matrix $\mc{F}(\bs{x})$, we can rewrite the above as
\begin{equation}\label{eq:taylor when nu=1}
\bs{F}(h(\bs{x}+\varpi^l\bs{v})+\bs{b})\equiv\bs{F}(h\bs{x}+\bs{b})+\mc{F}(\bs{x})\cdot\varpi^l\bs{v}\;\md{\varpi^{2l}}.\end{equation}
We want $\bs{v}$ to satisfy the system of congruences
\[\bs{F}(h(\bs{x}+\varpi^l\bs{v})+\bs{b})\equiv\bs{0}\;\md{\varpi^{2l}},\]
so that $\bs{y}=\bs{x}+\varpi^l\bs{v}$ is a solution of the congruences in \eqref{eq:F(x)=0 mod pi^2l-1+nu} with $\nu=1$. From \eqref{eq:taylor when nu=1}, the last displayed congruence relation is equivalent to
\[\bs{F}(h\bs{x}+\bs{b})+\mc{F}(\bs{x})\cdot\varpi^l\bs{v}\equiv\bs{0}\;\md{\varpi^{2l}}.\]
The assumptions \eqref{eq:F(x)=0 mod pi^2l-1} and \eqref{eq:X most singular mod pi^l-1} allow us to divide through the whole congruence by $\varpi^{2l-1}$ and hence rewrite it as
\[\varpi^{-(2l-1)}\bs{F}(h\bs{x}+\bs{b})+\varpi^{-(l-1)}\mc{F}(\bs{x})\cdot\bs{v}\equiv\bs{0}\;\md{\varpi}.\]
The latter half of \eqref{eq:X most singular mod pi^l-1} implies that the rank of the matrix $\varpi^{-(l-1)}\mc{F}(\bs{x})$ (mod $\varpi$) is $R$. Hence the system of linear  congruences above has at least $\size{\varpi}^{s-R}$ solutions in $\bs{v}$ (mod $\varpi$). For each of these tuples $\bs{v}$, the tuple $\bs{y}=\bs{x}+\varpi^l\bs{v}$ satisfies both \eqref{eq:F(x)=0 mod pi^2l-1+nu} and \eqref{eq:X most singular mod pi^l-1 II} in the current case for $\nu$. This finishes the proof of the base case of our induction.
\par Suppose the induction hypothesis holds for some positive integer $\nu$. This provides at least $\size{\varpi}^{\nu(s-R)}$ solutions of \eqref{eq:F(x)=0 mod pi^2l-1+nu} which satisfy \eqref{eq:X most singular mod pi^l-1 II} as well. Let $\bs{y}$ be one of these solutions. Again let $\bs{v}$ be an $s$-tuple (mod $\varpi$) to be chosen. Using the multidimensional Taylor's theorem once more, we obtain
\begin{align}\label{eq:taylor for nu}
\bs{F}(h(\bs{y}+\varpi^{l+\nu}\bs{v})+\bs{b})\equiv\;&\bs{F}(h\bs{y}+\bs{b})+\nabla\bs{F}(h\bs{y}+\bs{b})\cdot\varpi^{l+\nu}h\bs{v}\nonumber\\
=\;&\bs{F}(h\bs{y}+\bs{b})+\mc{F}(\bs{y})\cdot\varpi^{l+\nu}\bs{v}\;\md{\varpi^{2l+\nu}}.\end{align}
As in the case where $\nu=1$, we plan to choose $\bs{v}$ such that $\bs{F}(h(\bs{y}+\varpi^{l+\nu}\bs{v})+\bs{b})\equiv\bs{0}$ (mod $\varpi^{2l+\nu}$). Using \eqref{eq:taylor for nu}, we can rewrite this desired relation as
\[\bs{F}(h\bs{y}+\bs{b})+\varpi^{l+\nu}\mc{F}(\bs{y})\cdot\bs{v}\equiv\bs{0}\;\md{\varpi^{2l+\nu}}.\]
The conditions \eqref{eq:F(x)=0 mod pi^2l-1+nu} and \eqref{eq:X most singular mod pi^l-1 II} satisfied by $\bs{y}$ allow us to divide through this congruence by $\varpi^{2l-1+\nu}$, and hence rewrite it as
\[\varpi^{-(2l-1+\nu)}\bs{F}(h\bs{y}+\bs{b})+\varpi^{-(l-1)}\mc{F}(\bs{y})\cdot\bs{v}\equiv\bs{0}\;\md{\varpi}.\]
The latter half of \eqref{eq:X most singular mod pi^l-1 II} implies, as before, that this system of congruences has at least $\size{\varpi}^{s-R}$ solutions in $\bs{v}$ (mod $\varpi$). So altogether there are at least $\size{\varpi}^{\nu(s-R)}\size{\varpi}^{s-R}=\size{\varpi}^{(\nu+1)(s-R)}$ tuples $\bs{z}=\bs{y}+\varpi^{l+\nu}\bs{v}$ that satisfy \eqref{eq:F(x)=0 mod pi^2l-1+nu} and \eqref{eq:X most singular mod pi^l-1 II}, with $\bs{y}$ and $\nu$ replaced by $\bs{z}$ and $\nu+1$ respectively. This completes the proof of the inductive step. The lemma thus follows by induction on $\nu$.
\end{proof}
\end{lem}
We are now ready to establish that $\mf{S}>0$.
\begin{cor}\label{lem:S>0}
Suppose that the affine algebraic set $W$ has a non-singular $\varpi$-adic point for every $\varpi$. Then $\mf{S}>0$.
\begin{proof}
Courtesy of \eqref{eq:F(hx+b)=0}, our assumption here implies the existence of a $\varpi$-adic $s$-tuple $\bs{y}$ for which  $\rk\mc{F}(\bs{y})=R$ over $\mb{K}_{\varpi}$. There thus exists a positive integer $l$ together with an $s$-tuple $\bs{x}$ (mod $\varpi^{2l-1}$) for which the conditions \eqref{eq:F(x)=0 mod pi^2l-1} and \eqref{eq:X most singular mod pi^l-1} of Lemma \ref{lem:smooth local solutions} are satisfied. Recalling that $M(\varpi^{2l-1+\nu})$ counts the number of solutions of the congruences \eqref{eq:M(g)} with $g=\varpi^{2l-1+\nu}$, and using Lemma \ref{lem:smooth local solutions}, we obtain
\[M(\varpi^{2l-1+\nu})\geq\size{\varpi}^{\nu(s-R)}=C\size{\varpi}^{(2l-1+\nu)(s-R)},\]
where $C=\size{\varpi}^{(2l-1)(R-s)}$. Evidently $C$ is independent of $\nu$. We have thus established the desired inequality \eqref{eq:many local solutions}, with $k=2l-1+\nu$. Lemma \ref{lem:local solutions} therefore gives $\Omega(\varpi)>0$ for every $\varpi$ with $\size{\varpi}<\wh{Z}$, where $Z$ is given in Lemma \ref{lem:factorise S}.That lemma then leads to the conclusion that $\mf{S}>0$.
\end{proof}
\end{cor}
\section{Singular integral}
Similar to the way that the singular series is treated in the previous section, our aim now is to extend the finite integral $\mc{J}(P/2)$ in the formula \eqref{eq:approx asymptotic for rho(P;N)} to the \emph{complete singular integral}
\begin{equation}\label{eq:J}
\mc{J}=\int_{\mb{K}_{\infty}^R}I(\bs{\gamma};\mc{B})\;\mr{d}\bs{\gamma}=\lim_{Y\rightarrow\infty}\mc{J}(Y).\end{equation}
We therefore need to know whether this integral converges absolutely, and the size of the error incurred by replacing $\mc{J}(P/2)$ by $\mc{J}$ in \eqref{eq:approx asymptotic for rho(P;N)}. To make the resulting asymptotic formula for $\rho_{h,\bs{b}}(n;\mf{N})$ useful, we also require that $\mc{J}>0$, under certain conditions on the geometry of $X$. The first two pieces of information are obtained via an estimate for the exponential integral $I(\bs{\gamma};\mc{B})$ given in \eqref{eq:I}. This is provided by the following lemma.
\begin{lem}\label{lem:upper bound for I}
Let $\bs{\gamma}\in\mb{K}_{\infty}^R$ and $M\in\mb{N}$. Let $\mc{C}$ be any hypercube in $\mb{T}^s$ with sidelength $\wh{M}^{-1}$. Then provided that $p>d$, we have
\[I(\bs{\gamma};\mc{C})\ll\wh{M}^{-s}(1+\wh{M}^{-d}\size{\bs{\gamma}})^{-K/(R(d-1))}.\]
\begin{proof}
From \eqref{eq:I}, we see that when $\size{\bs{\gamma}}\leq\wh{M}^{d}$, the estimate given by the lemma is trivial. Hence we may suppose that
\begin{equation}\label{eq:WLOG on gamma}
\size{\bs{\gamma}}>\wh{M}^d\geq1.\end{equation}
Let $Q$ be a large parameter to be chosen. The plan is to apply Lemma \ref{lem:approximate T(alpha)} to relate $I(\bs{\gamma};\mc{C})$ to the generating function
\[V(\bs{\alpha})=\sum_{\bs{x}\in t^Q\mc{C}}E(\bs{\alpha}\cdot\bs{F}(\bs{x})),\]
for some $\bs{\alpha}$ to be specified, and then use Lemma \ref{lem:weyl's inequality} to obtain an upper bound for the generating function.
\par Provided that $Q$ satisfies the conditions
\[\wh{Q}>\wh{M}\;\text{and }\size{\bs{\alpha}}<\wh{\Delta}^{-1}\wh{Q}^{1-d},\]
Lemma \ref{lem:approximate T(alpha)} is applicable with $\bs{a}=0$, $g=1$, $h=1$, $m=t^Q$ and $\bs{b}=\bs{0}$. From \eqref{eq:S(g,a)}, we have $S(1,\bs{0})=1$, so Lemma \ref{lem:approximate T(alpha)} yields
\begin{equation}\label{eq:use T(alpha) for I}
V(\bs{\alpha})=\wh{Q}^sI(t^{Qd}\bs{\alpha};\mc{C}).\end{equation}
Choose $\bs{\alpha}=t^{-Qd}\bs{\gamma}$. The second condition on $Q$ above is thus satisfied as long as $\wh{Q}>\wh{\Delta}\size{\bs{\gamma}}$. Meanwhile, write $\mc{E}=t^M\mc{C}$. Then $\mc{E}$ is a hypercube in $\mb{T}^s$. An application of Lemma \ref{lem:weyl's inequality} with $L=K$, $h=g=1$, $\bs{b}=\bs{a}=\bs{0}$ and $Q$ replaced by $Q-M$, gives the estimate
\begin{align*}
&V(\bs{\alpha})\\
=\;&\sum_{\bs{x}\in t^{Q-M}\mc{E}}E(\bs{\alpha}\cdot\bs{F}(\bs{x}))\\
\ll\;&(\wh{Q}\wh{M}^{-1})^s\Big(\frac{(\wh{Q}\wh{M}^{-1})^R+(\wh{Q}\wh{M}^{-1})^{Rd}\size{\bs{\alpha}}^R}{(\wh{Q}\wh{M}^{-1})^{Rd}}+\frac{1}{1+(\wh{Q}\wh{M}^{-1})^{d}\size{\bs{\alpha}}}\Big)^{K/(R(d-1))}\\
=\;&\wh{Q}^s\wh{M}^{-s}\Big((\wh{Q}\wh{M}^{-1})^{R(1-d)}+\wh{Q}^{-Rd}\size{\bs{\gamma}}^R+\frac{1}{1+\wh{M}^{-d}\size{\bs{\gamma}}}\Big)^{K/(R(d-1))}\\
\ll\;&\wh{Q}^s\wh{M}^{-s}\Big((\wh{Q}\wh{M}^{-1})^{R(1-d)}+\wh{Q}^{-Rd}\size{\bs{\gamma}}^R+\wh{M}^d\size{\bs{\gamma}}^{-1}\Big)^{K/(R(d-1))}.\end{align*}
Here the last two lines follow from our choice of $\bs{\alpha}$ and the assumption \eqref{eq:WLOG on gamma} respectively. On combining this with \eqref{eq:use T(alpha) for I}, we get
\[I(\bs{\gamma};\mc{C})\ll\wh{M}^{-s}\Big((\wh{Q}\wh{M}^{-1})^{R(1-d)}+\wh{Q}^{-Rd}\size{\bs{\gamma}}^R+\wh{M}^d\size{\bs{\gamma}}^{-1}\Big)^{K/(R(d-1))}.\]
Note that the desired upper bound for $I(\bs{\gamma};\mc{C})$ follows if the last term in the bracket above dominates. This is indeed the case if we choose $Q$ so large that
\[\wh{Q}>\ma{\wh{M},\wh{\Delta}\size{\bs{\gamma}},\size{\bs{\gamma}}^{(R+1)/(Rd)}\wh{M}^{-1/R},\wh{M}^{1-d/(R(d-1))}\size{\bs{\gamma}}^{1/(R(d-1))}}.\]
This completes the proof of the lemma.
\end{proof}
\end{lem}
We are now in a position to resolve the convergence issue of the singular integral.
\begin{cor}\label{lem:J converges}
Let $Y$ be any real number with $Y\geq dN+\Delta$. Then when $K>R^2(d-1)$, we have
\[\mc{J}-\mc{J}(Y)\ll\wh{Y}^{R-K/(R(d-1))}.\]
In particular, the singular integral $\mc{J}$ converges absolutely. Further, we have $\mc{J}\ll1$.
\begin{proof}
Recall from \eqref{eq:B} that $\mc{B}$ is a hypercube with sidelength $\wh{N}^{-1}$, with $N$ absolutely bounded. By Lemma \ref{lem:upper bound for I} followed by $R$ applications of \eqref{eq:order less than -m}, we get
\begin{align*}
\mc{J}-\mc{J}(Y)\ll\;&\int_{\size{\bs{\gamma}}\geq \wh{\Delta}^{-1}\wh{Y}}|I(\bs{\gamma};\mc{B})|\;\mr{d}\bs{\gamma}\\
\ll\;&\int_{\size{\bs{\gamma}}\geq \wh{\Delta}^{-1}\wh{Y}}\size{\bs{\gamma}}^{-K/(R(d-1))}\;\mr{d}\bs{\gamma}\\
=\;&\sum_{l\geq Y-\Delta}(q^l)^{-K/(R(d-1))}\int_{\size{\bs{\gamma}}=\wh{l}}\mr{d}\bs{\gamma}\\
\ll\;&\sum_{l\geq Y-\Delta}(q^l)^{-K/(R(d-1))}(q^l)^R\\
\ll\;&\wh{Y}^{R-K/(R(d-1))}.\end{align*}
Here the last line follows from the assumption that $K>R^2(d-1)$. This proves the first assertion of the corollary. The second assertion is an immediate consequence of this. Meanwhile, a parallel argument to the one given above yields
\begin{align*}
\mc{J}(Y)\ll\;&\wh{N}^{-s}\int_{\size{\bs{\gamma}}<\wh{\Delta}^{-1}\wh{Y}}(1+\wh{N}^{-d}\size{\bs{\gamma}})^{-K/(R(d-1))}\;\mr{d}\bs{\gamma}\\
\ll\;&\int_{\size{\bs{\gamma}}<\wh{N}^d}1\;\mr{d}\bs{\gamma}+\int_{\wh{N}^d\leq\size{\bs{\gamma}}<\wh{Y}}\size{\bs{\gamma}}^{-K/(R(d-1))}\;\mr{d}\bs{\gamma}\\
=\;&\wh{N}^d+\sum_{dN\leq l<Y}(q^l)^{-K/(R(d-1))}\int_{\size{\bs{\gamma}}=\wh{l}}\;\mr{d}\bs{\gamma}\\
\ll\;&1+\sum_{dN\leq l<Y}(q^l)^{R-K/(R(d-1))}\\
\ll\;&1.\end{align*}
This establishes the third assertion of the corollary.
\end{proof}
\end{cor}
It therefore remains to establish that $\mc{J}>0$. On considering the relation \eqref{eq:J}, it suffices to show that $\mc{J}(Y)\gg1$ for large $Y$. Recalling the definitions \eqref{eq:J(Y)} and \eqref{eq:I} of $\mc{J}(Y)$ and $I(\bs{\gamma};\mc{B})$ respectively, and applying Fubini's theorem, we obtain
\begin{align*}
\mc{J}(Y)=\;&\int_{\size{\bs{\gamma}}<\wh{\Delta}^{-1}\wh{Y}}I(\bs{\gamma};\mc{B})\;\mr{d}\bs{\gamma}\\
=\;&\int_{\size{\bs{\gamma}}<\wh{\Delta}^{-1}\wh{Y}}\int_{\mc{B}}E(\bs{\gamma}\cdot\bs{F}(\bs{\sigma}))\;\mr{d}\bs{\sigma}\mr{d}\bs{\gamma}\\
=\;&\int_{\mc{B}}\int_{\size{\bs{\gamma}}<\wh{\Delta}^{-1}\wh{Y}}E(\bs{\gamma}\cdot\bs{F}(\bs{\sigma}))\;\mr{d}\bs{\gamma}\mr{d}\bs{\sigma}.\end{align*}
With $R$ applications of Lemma 1(f) of \cite{kubota}, the $\bs{\gamma}$-integral above gives
\begin{equation}\label{eq:geometric singular integral}
\mc{J}(Y)=\wh{\Delta}^{-R}\wh{Y}^RM_{\mc{B}}(Y),\end{equation}
where $M_{\mc{B}}(Y)=\vol\mc{M}_{\mc{B}}(Y)$, with
\begin{equation}\label{eq:M_B(Y)}
\mc{M}_{\mc{B}}(Y)=\left\{\bs{\sigma}\in\mc{B}:\size{\bs{F}(\bs{\sigma})}<\wh{\Delta}\wh{Y}^{-1}\right\}.\end{equation}
To have the relation $\mc{J}(Y)\gg1$, it is thus sufficient to prove that
\begin{equation}\label{eq:M_B(Y) large}
\mc{M}_{\mc{B}}(Y)\gg\wh{Y}^{-R}\qquad\text{for large }Y.\end{equation}
To this end, write
\[\mc{M}_{\mc{B}}=\bigcap_{Y\geq\Delta}\mc{M}_{\mc{B}}(Y).\]
Equations \eqref{eq:M_B(Y)} and \eqref{eq:X} provide the relation
\begin{equation}\label{eq:M_B}
\mc{M}_{\mc{B}}=\left\{\bs{\sigma}\in\mc{B}:\bs{F}(\bs{\sigma})=\bs{0}\right\}=\mc{B}\cap X.\end{equation}
The following lemma provides a sufficient geometric condition under which \eqref{eq:M_B(Y) large} is true.
\begin{lem}\label{lem:J>0}
Provided that $\dim\mc{M}_{\mc{B}}\geq s-R$, we have $\mc{J}>0$.
\begin{proof}
We use the same argument used in the proof of Lemma 2 of \cite{schmidt}. Write $\mc{B}(Y)$ for the set of all points in $\mc{B}$ which are at least a distance $\wh{Y}^{-1}$ away from the boundary of $\mc{B}$. Set $\mc{M}'=\mc{M}_{\mc{B}}\cap\mc{B}(Y)$. Choose $\mc{M}''$ to be a component of $\mc{M}'$ with precise dimension $s-R$. Choose another subset $\mc{M}'''\subseteq\mc{M}''$ which can be parametrised by $s-R$ coordinates. Without loss of generality, we can assume these to be the first $s-R$ coordinates $\sigma_1,...,\sigma_{s-R}$. This means that we can find a non-empty open set $\mc{O}\subseteq\mb{T}^{s-R}$ together with a map $\bs{h}=(h_1,...,h_R):\mc{O}\rightarrow\mb{T}^R$ such that
\begin{equation}(\sigma_1,...,\sigma_{s-R},\bs{h}(\sigma_1,...,\sigma_{s-R}))\in\mc{M}^{'''}\qquad\text{when }(\sigma_1,...,\sigma_{s-R})\in\mc{O}.\label{eq:in M'''}\end{equation}
Since $\mc{M}'''\subseteq\mc{M}_{\mc{B}}$, equation \eqref{eq:M_B} gives
\begin{equation}\label{eq:F maps O to 0}
\bs{F}(\sigma_1,...,\sigma_{s-R},\bs{h}(\sigma_1,...,\sigma_{s-R}))=\bs{0}\qquad\text{when }(\sigma_1,...,\sigma_{s-R})\in\mc{O}.\end{equation}
For any positive integer $T$, write $S_T$ for the set of $s$-tuples $(\sigma_1,...,\sigma_s)\in\mb{T}^s$ with $(\sigma_1,...,\sigma_{s-R})\in\mc{O}$ such that
\[\size{\sigma_{s-R+i}-h_i(\sigma_1,...,\sigma_{s-R})}<\wh{T}^{-1}\;\text{for all }i=1,...,R.\]
Our choice of $\mc{M}'''$ indicates that all points of $\mc{M}'''$ are in $\mc{B}$ and at least a distance $\wh{Y}^{-1}$ away from the boundary of $\mc{B}$. It thus follows from \eqref{eq:in M'''} and the inequality above that $S_T\subseteq\mc{B}$ whenever $T>Y$. Further, observe that
\begin{align*}
&\size{(\sigma_1,...,\sigma_{s-R},\bs{h}(\sigma_1,...,\sigma_{s-R}))-(\sigma_1,...,\sigma_s)}\\
=\;&\max_{1\leq i\leq R}\size{\sigma_{s-R+i}-h_i(\sigma_1,...,\sigma_{s-R})}.\end{align*}
On recalling \eqref{eq:F maps O to 0} and the fact that the forms $\bs{F}$ have degree $d$ with coefficients having absolute value at most $\wh{\Delta}$, we deduce that when $(\sigma_1,...,\sigma_{s-R})\in\mc{O}$, we have
\begin{align*}
\size{\bs{F}(\sigma_1,...,\sigma_s)}=\;&\size{\bs{F}(\sigma_1,...,\sigma_s)-\bs{F}(\sigma_1,...,\sigma_{s-R},\bs{h}(\sigma_1,...,\sigma_{s-R}))}\\
\leq\;&\wh{\Delta}\size{(\sigma_1,...,\sigma_s)}^{d-1}\max_{1\leq i\leq R}\size{\sigma_{s-R+i}-h_i(\sigma_1,...,\sigma_{s-R})}.\end{align*}
In particular, when $(\sigma_1,...,\sigma_s)\in S_{Y}$, we get
\[\size{\bs{F}(\sigma_1,...,\sigma_s)}<\wh{\Delta}\wh{Y}^{-1}.\]
From \eqref{eq:M_B(Y)}, the above inequality implies that  $S_{Y}\subseteq\mc{M}_{\mc{B}}(Y)$. It is apparent from the definition of $S_{Y}$ that $\vol S_{Y}\gg\wh{Y}^{-R}$. The last inclusion relation thus implies that $\vol\mc{M}_{\mc{B}}(Y)\gg\wh{Y}^{-R}$. This is exactly \eqref{eq:M_B(Y) large}. The conclusion of the lemma thus follows by inserting this estimate into \eqref{eq:geometric singular integral}, and recalling the relation \eqref{eq:J}.
\end{proof}
\end{lem}
With this result, we only have to show that $\dim\mc{M}_{\mc{B}}\geq s-R$. But we have chosen our box $\mc{B}$ to be centred at a non-singular real point on $X$. By the Implicit Function Theorem in arbitrary fields (see Theorem 2.2.1 in \cite{igusa}), the set $\mc{M}_{\mc{B}}$ of real points in $X\cap\mc{B}$ has dimension at least $s-R$, as required. On combining this with Lemma \ref{lem:J>0}, we conclude that $\mc{J}>0$, under the condition that $X$ contains a non-singular real point.
\section{Conclusion}
We are now in a position to settle the outstanding theorems in this paper. We first derive an asymptotic formula for the number $\rho_{h,\bs{b}}(n)$ of integral solutions $\bs{x}$ to \eqref{eq:F(hx+b)=0} with $\bs{x}\in n\mc{B}_N$. Applying Lemma \ref{lem:S converges} and Corollary \ref{lem:J converges} in \eqref{eq:approx asymptotic for rho(P;N)} with $Y=P/2$, we obtain
\[
\rho_{h,\bs{b}}(n;\mf{N})=\size{h}^{-Rd}\wh{P}^{s-Rd}(\mf{S}\mc{J}+O((\wh{P}^{1/2})^{R+1-K/(R(d-1))})).\]
On recalling that $\mb{T}^R=\mf{N}\cup\mf{n}$, and using Lemma \ref{lem:minor arc estimate}, we obtain
\begin{align}\label{eq:generic asymptotic}
\rho_{h,\bs{b}}(n)=\;&\rho_{h,\bs{b}}(n;\mf{N})+\rho_{h,\bs{b}}(n;\mf{n})\nonumber\\
=\;&\size{h}^{-Rd}\wh{P}^{s-Rd}(\mf{S}\mc{J}+O(\wh{P}^{-(K/(R(d-1))-R-1)/4})),\end{align}
for every monic polynomials $n,h\in\mb{A}$ with $\size{n}=\wh{P}$, and every integral $s$-tuple $\bs{b}$ with $\size{\bs{b}}<\size{h}$.
\par Now we make different choices for $h$, $\bs{b}$, $n$ and $N$ in deriving Theorems \ref{thm:asymptotic formula} and \ref{thm:weak approx}. To establish Theorem \ref{thm:asymptotic formula}, we take $h=1$, $\bs{b}=\bs{0}$, $n=t^P$, $N=0$ and $\bs{\xi}$ to be an arbitrary point in $\mb{T}^s$. From \eqref{eq:B}, we see that $\mc{B}_0=\mb{T}^s$. Meanwhile, on considering the respective underlying equations \eqref{eq:F(hx+b)=0} and \eqref{eq:X} of the algebraic sets $W$ and $X$ respectively, we have $W=X$. Thus $\rho_{1,\bs{0}}(n)=\rho(P)$, where $\rho(P)$ is defined just before Theorem \ref{thm:asymptotic formula}. On recalling \eqref{eq:assumption on k} and \eqref{eq:K}, equation \eqref{eq:generic asymptotic} thus gives
\[\rho(P)=\wh{P}^{s-Rd}(\mf{S}\mc{J}+O(\wh{P}^{-k/(2^{d+1}R(d-1))})).\]
This is the asymptotic formula required in Theorem \ref{thm:asymptotic formula}. Further, Corollary \ref{lem:S>0} and the discussion after Lemma \ref{lem:J>0} imply that the quantities $\mf{S}$ and $\mc{J}$ are indeed positive under the hypotheses of Theorem \ref{thm:asymptotic formula}. The proof of that theorem is thus complete.
\par Before proving Theorem \ref{thm:weak approx}, we introduce a few more notations. For every $\varpi$, let $\size{\cdot}_{\varpi}$ be the $\varpi$-adic absolute value on $\mb{K}$. More precisely, when $x\in\mb{A}\backslash\left\{0\right\}$, we put $\size{x}_{\varpi}=\size{\varpi}^{-l}$, where $l$ is the non-negative integer such that $\varpi^l\|x$. By convention, put $\size{0}_{\varpi}=0$. We extend this absolute value to $\mb{K}$ by the relation
\[\size{x/y}_{\varpi}=\size{x}_{\varpi}/\size{y}_{\varpi}\qquad(x,y\in\mb{A},y\neq 0).\]
Let $\mb{A}_{\varpi}$ and $\mb{K}_{\varpi}$ denote respectively the completions of $\mb{A}$ and $\mb{K}$ with respect to this absolute value.
\par The proof of Theorem \ref{thm:weak approx} is largely based on the argument which led to Corollary 1 of \cite{skinner}. Note that \eqref{eq:generic asymptotic} has the following immediate consequence.
\begin{lem}\label{lem:weak approx}
Suppose $\bs{F}$ has a non-singular real zero $\bs{\xi}$. Let $h\in\mb{A}$, $\bs{b}\in\mb{A}^s$, with $h$ monic and $\size{\bs{b}}<\size{h}$, such that for every $\varpi$, there exists a non-singular solution $\bs{x}$ of the equations $\bs{F}(h\bs{x}+\bs{b})=\bs{0}$, with $\bs{x}\in\mb{A}_{\varpi}^s$. Then for every positive integer $N$, there exists a $\bs{\xi}_0\in\mb{A}^s$ such that
\begin{equation}\label{eq:xi_0}
\bs{F}(h\bs{\xi}_0+\bs{b})=\bs{0},\end{equation}
and there exists a positive constant $C=C(s,d,R,\bs{F},h,\bs{b},N)$ such that whenever $n$ is a monic polynomial in $\mb{A}$ with $\size{n}\geq\wh{C}$, we have
\begin{equation}\label{eq:xi_0 close to xi}
\size{n^{-1}\bs{\xi}_0-\bs{\xi}}<\wh{N}^{-1}.\end{equation}
\begin{proof}
From \eqref{eq:B} and the definition of $\rho_{h,\bs{b}}(n)$ given before \eqref{eq:rho(P)}, we see that $\rho_{h,\bs{b}}(n)$ counts the number of integral points $\bs{\xi}_0$ which satisfy \eqref{eq:xi_0} and \eqref{eq:xi_0 close to xi}. Hence, under the conditions that $\mf{S}>0$ and $\mc{J}>0$, the existence of such $\bs{\xi}_0$ is confirmed by \eqref{eq:generic asymptotic} when $\size{n}=\wh{P}\geq\wh{C}$, for some positive constant $C=C(s,d,R,\bs{F},h,\bs{b},N)$. But our present hypothesis on the existence of non-singular $\varpi$-adic solutions for every $\varpi$ renders Corollary \ref{lem:S>0} applicable. That corollary gives  $\mf{S}>0$. Also, now that $X$ has a non-singular real point, the conclusion of \S8 implies that $\mc{J}>0$. This completes the proof of the lemma.
\end{proof}
\end{lem}
We now apply the above lemma to establish Theorem \ref{thm:weak approx}, which we reformulate below.
\begin{lem}\label{lem:weak approx II}
Suppose the hypotheses of Theorem \ref{thm:weak approx} hold. Let $\bs{\zeta}_{\infty}\in X(\mb{K}_{\infty})$. Fix a collection of distinct monic irreducible polynomials $\varpi_1,...,\varpi_r$. For each $i=1,...,r$, let $\bs{\zeta}_i\in X(\mb{K}_{\varpi_i})$. Let $N_{\infty},N_1,...,N_r$ be positive integers. Then there exists a $\bs{\zeta}\in X(\mb{K})$ such that
\[\size{\bs{\zeta}-\bs{\zeta}_{\infty}}<\wh{N}_{\infty}^{-1}\]
and
\[\size{\bs{\zeta}-\bs{\zeta}_i}_{\varpi_i}\leq\size{\varpi_i}^{-N_i}\qquad\text{for all }i=1,...,r.\]
\begin{proof}
Let $h=\varpi_1^{N_1}...\varpi_r^{N_r}$ and $\bs{\xi}=h^{-1}\bs{\zeta}_{\infty}$. By the homogeneity of $\bs{F}$, we can assume without loss of generality that $\bs{\zeta}_i\in\mb{A}_{\varpi_i}^s$ for each $i=1,...,r$. Let $\bs{b}_i\in\mb{A}^s$ such that $\bs{\zeta}_i\equiv\bs{b}_i\;\md{\varpi_i^{N_i}}$. By the Chinese Remainder Theorem, we can find a $\bs{b}\in\mb{A}^s$ such that $\bs{b}\equiv\bs{b}_i\;\md{\varpi_i^{N_i}}$, for each $i$. 
%
\par We first verify that, with the above choices of $h$ and $\bs{b}$, the $\varpi$-adic solubility condition of Lemma \ref{lem:weak approx} is satisfied for every $\varpi$. First consider the case where $\varpi\notin\left\{\varpi_1,...,\varpi_r\right\}$. The hypotheses of Theorem \ref{thm:weak approx} imply the existence of a $\bs{\zeta}^{(\varpi)}\in X(\mb{K}_{\varpi})$. By homogeneity again, we can assume that $\bs{\zeta}^{(\varpi)}\in X(\mb{A}_{\varpi})$. Take $\bs{x}^{(\varpi)}=(\bs{\zeta}^{(\varpi)}-\bs{b})/h$. Then $\bs{F}(h\bs{x}^{(\varpi)}+\bs{b})=\bs{0}$. Since $(h,\varpi)=1$, it follows that $\bs{x}^{(\varpi)}\in\mb{A}_{\varpi}^s$. Further, now that $X$ is non-singular, it follows that $\bs{x}^{(\varpi)}$ is a non-singular solution of the equations $\bs{F}(h\bs{x}+\bs{b})=\bs{0}$. Now consider the monic irreducibles amongst $\varpi_1,...,\varpi_r$. Fix $i\in\left\{1,...,r\right\}$. As in the previous case, take $\bs{x}_i=(\bs{\zeta}_i-\bs{b})/h$. Since $\varpi_i^{N_i}\|h$ and $\bs{\zeta}_i\equiv\bs{b}\;\md{\varpi_i^{N_i}}$, it follows that $\bs{x}_i\in\mb{A}_{\varpi_i}$. The non-singularity of $X$ once again implies that $\bs{x}_i$ is a non-singular solution of $\bs{F}(h\bs{x}+\bs{b})=\bs{0}$.
We have thus shown that for every $\varpi$, the $\varpi$-adic solubility condition of Lemma \ref{lem:weak approx} is indeed satisfied.
\par We apply that lemma with $\wh{N}=\wh{N}_{\infty}\size{h}$. This gives an integral $s$-tuple $\bs{\xi}_0$ such that
\begin{equation}\label{eq:xi_0 solution}
\bs{F}(h\bs{\xi}_0+\bs{b})=\bs{0},\end{equation}
together with a positive constant $C=C(s,d,R,\bs{F},h,\bs{b})$ such that whenever $n$ is a monic polynomial in $\mb{A}$ with $\size{n}\geq\wh{C}$, we have
\[\size{n^{-1}\bs{\xi}_0-\bs{\xi}}<\wh{N}_{\infty}^{-1}\size{h}^{-1}.\]
Recalling our choice of $\bs{\xi}$ at the beginning of this proof, we multiply through the last relation by $h$ and obtain \begin{equation}\label{eq:xi_0 close to zeta}
\size{hn^{-1}\bs{\xi}_0-\bs{\zeta}_{\infty}}<\wh{N}_{\infty}^{-1}.\end{equation}
Choose $n$ to be a monic polynomial with 
\begin{equation}\label{eq:congruence on n}
n\equiv1\;\md{\varpi_1^{N_1}...\varpi_r^{N_r}},\end{equation}
such that 
\begin{equation}\label{eq:size of n}\size{n}>\ma{\size{\bs{b}}\wh{N}_{\infty},\wh{C}}.\end{equation}
Take $\bs{\zeta}=n^{-1}(h\bs{\xi}_0+\bs{b})$. By the homogeneity of $\bs{F}$ together with \eqref{eq:xi_0 solution}, we have $\bs{\zeta}\in X(\mb{K})$. Using \eqref{eq:xi_0 close to zeta} together with \eqref{eq:size of n}, we also get
\begin{align*}
\size{\bs{\zeta}-\bs{\zeta}_{\infty}}=\;&\size{n^{-1}(h\bs{\xi}_0+\bs{b})-\bs{\zeta}_{\infty}}\\
\leq\;&\ma{\size{hn^{-1}\bs{\xi}_0-\bs{\zeta}_{\infty}},\size{\bs{b}}/\size{n}}<\wh{N}_{\infty}^{-1}.\end{align*}
Finally, our choice of $\bs{\zeta}$ implies that
\begin{align}\label{eq:zeta_0 close to zeta_i}
\size{\bs{\zeta}-\bs{\zeta}_i}_{\varpi_i}\leq\;&\ma{\size{\bs{\zeta}-\bs{\zeta}_i/n}_{\varpi_i},\size{\bs{\zeta}_i/n-\bs{\zeta}_i}_{\varpi_i}}\nonumber\\
=\;&\ma{\size{h\bs{\xi}_0+\bs{b}-\bs{\zeta}_i}_{\varpi_i}/\size{n}_{\varpi_i},\size{\bs{\zeta}_i}_{\varpi_i}\size{n-1}_{\varpi_i}/\size{n}_{\varpi_i}}.\end{align}
Moreover, our choices of $h$ and $\bs{b}$ reveal that $h\equiv0\;\md{\varpi_i^{N_i}}$ and $\bs{b}\equiv\bs{\zeta}_i\;\md{\varpi_i^{N_i}}$. Hence
\[\size{h\bs{\xi}_0+\bs{b}-\bs{\zeta}_i}_{\varpi_i}\leq\size{\varpi_i}^{-N_i}.\]
Meanwhile, relation \eqref{eq:congruence on n} gives $\size{n-1}_{\varpi_i}\leq\size{\varpi_i}^{-N_i}$ and $\size{n}_{\varpi_i}=1$. These observations together with \eqref{eq:zeta_0 close to zeta_i} thus imply that
\[\size{\bs{\zeta}-\bs{\zeta}_i}_{\varpi_i}\leq\;\size{\varpi_i}^{-N_i}\qquad\text{for all }i=1,...,r.\]
This completes the proof of the lemma, and thus that of Theorem \ref{thm:weak approx}.
\end{proof}
\end{lem}
Our task in this paper is therefore complete, subject to the resolution of a number of technical results given in the two appendices below.
\appendix
\section{Integration over function fields}
In this appendix, we state and prove the change-of-variable property for integrals over function fields. We investigate integrals in the form
\[
\int_{\Gamma}f(\bs{\alpha})\;\mr{d}\bs{\alpha},\]
for some suitable function $f:\Gamma\rightarrow\mb{C}$ defined on some $D$-dimensional box $\Gamma$. Henceforth, a $D$-dimensional box refers to the set of real $D$-tuples $\bs{x}$ satisfying constraints of the form
\begin{equation}\label{eq:box}
\size{x_u}<\wh{R}_u\qquad(u=1,...,D),\end{equation}
for some real numbers $R_u$ ($1\leq u\leq D$). The special case where $D=1$ and $f=E\circ\varphi$, where $E$ is defined in \eqref{defn:e} and $\varphi$ is a polynomial map in one variable, is investigated in Proposition I.6 of \cite{car_exp}. The aim here is then to extend that result to rather more general functions defined over a box as in \eqref{eq:box}.
\par For convenience of exposition, we introduce further notations. For any $\zeta\in\mb{K}_{\infty}\backslash\left\{0\right\}$ and any integer $Y$, write $\Theta_Y(\zeta)$ for the coefficient of $t^Y$ in the expansion for $\zeta$. By convention, we define $\Theta_Y(0)=0$ for any integer $Y$. In addition, for every integer $L$, write $\mc{D}_L$ for the set of all $\alpha\in\mb{K}_{\infty}$ such that $\size{\alpha}<\wh{L}$. We ultimately establish the following theorem.
\begin{thm}\label{thm:change of variables}
Let $D$ be a positive integer, and let $\Gamma$ be any box in $\mb{K}_{\infty}^D$. Let $f:\Gamma\rightarrow\mb{C}$ be continuous. Suppose $M\in GL_D(\mb{K}_{\infty})$. Then
\begin{equation}\label{eq:eqn in change of variables}
\int_{\Gamma}f(\bs{\alpha})\;\mr{d}\bs{\alpha}=\size{\det M}\int_{M\bs{\gamma}\in\Gamma}f(M\bs{\gamma})\;\mr{d}\bs{\gamma}.\end{equation}
\end{thm}
Since $GL_D(\mb{K}_{\infty})$ is generated by the elementary matrices, namely permutation, diagonal and unipotent matrices, it suffices to prove Theorem \ref{thm:change of variables} when $M$ is one of these matrices. When $M$ is a permutation matrix, Theorem \ref{thm:change of variables} is equivalent to the property that we can relabel the variables of integration without changing the value of the integral of interest. This is evidently true, so we concentrate on the remaining two cases for $M$. Induction on the dimension $D$ is central to our treatment of these cases. Hence, for diagonal matrices $M$, we first require the following lemma, which deals with the inductive base where $D=1$.
\begin{lem}\label{lem:change of variables}
Let $f:\mc{D}_L\rightarrow\mb{C}$ be a continuous function. Let $L$ be an integer, and $\beta\in\mb{K}_{\infty}\backslash\left\{0\right\}$. Then
\begin{equation}\label{eq:eqn of change of variables lemma}
\int_{\size{\alpha}<\wh{L}}f(\alpha)\;\mr{d}\alpha=\size{\beta}\int_{\size{\gamma}<\wh{L}\size{\beta}^{-1}}f(\beta\gamma)\;\mr{d}\gamma.\end{equation}
\begin{proof}
Since $f$ is continuous on the compact set $\mc{D}_L$, it is uniformly continuous. For the sake of simplicity, when $H,J\in\mb{Z}$ with $-H\leq J-1$, and $a_{-H},...,a_{J-1}\in\mb{F}_q$, write
\[g=g(\bs{a};H,J)=a_{-H}t^{-H}+...+a_{J-1}t^{J-1}.\]
First we prove that for any integer $J$ and any uniformly continuous function $h$ on $\mc{D}_J$, we have
\begin{equation}\label{eq:limit repn of integral}
\int_{\size{\alpha}<\wh{J}}h(\alpha)\;\mr{d}\alpha=\lim_{H\rightarrow\infty}\lambda(H,J),\end{equation}
where
\begin{equation}\label{eq:lambda(H,J)}
\lambda(H,J)=\wh{H}^{-1}\sum_{\substack{a_k\in\mb{F}_q\\(-H\leq k\leq J-1)}}h(g).\end{equation}
For any $\alpha\in\mc{D}_J$, we make use of the substitution $\alpha=g+\gamma$, where the $a_k$ ($-H\leq k\leq J-1$) run through elements of $\mb{F}_q$ and $\gamma$ runs through elements of $\mb{K}_{\infty}$ with $\size{\gamma}<\wh{H}^{-1}$. This yields
\[\int_{\size{\alpha}<\wh{J}}h(\alpha)\;\mr{d}\alpha=\int_{\size{\gamma}<\wh{H}^{-1}}\sum_{\substack{a_k\in\mb{F}_q\\(-H\leq k\leq J-1)}}h(g+\gamma)\;\mr{d}\gamma,\]
for all positive integers $H$. Meanwhile, on recalling \eqref{eq:order less than -m}, equation \eqref{eq:lambda(H,J)} can be rewritten as
\[\lambda(H,J)=\int_{\size{\gamma}<\wh{H}^{-1}}\sum_{\substack{a_k\in\mb{F}_q\\(-H\leq k\leq J-1)}}h(g)\;\mr{d}\gamma.\]
To prove \eqref{eq:limit repn of integral}, it thus suffices to show that, for each fixed choice of $J$ and any $\varepsilon>0$, there exists a positive integer $H_0=H_0(\varepsilon,J)$ such that when $H\in\mb{Z}$ with $H\geq H_0$, we have
\[\Big|\int_{\size{\gamma}<\wh{H}^{-1}}\sum_{\substack{a_k\in\mb{F}_q\\(-H\leq k\leq J-1)}}h(g+\gamma)\;\mr{d}\gamma-\int_{\size{\gamma}<\wh{H}^{-1}}\sum_{\substack{a_k\in\mb{F}_q\\(-H\leq k\leq J-1)}}h(g)\;\mr{d}\gamma\Big|<\varepsilon.\]
There are exactly $q^{J-(-H)}=\wh{J}\wh{H}$ choices for the $a_k$ ($-H\leq k\leq J-1$) in the sum above, whereas the $\gamma$-integral above is over a set of measure $\wh{H}^{-1}$, courtesy of \eqref{eq:order less than -m}. We can therefore reduce the task of establishing the last displayed relation to showing that, for every $\varepsilon>0$, there exists a positive integer $H_0=H_0(\varepsilon,J)$ such that whenever $H\geq H_0$, we have
\[|h(g+\gamma)-h(g)|<\varepsilon\wh{J}^{-1},\]
for any $\gamma\in\mc{D}_{-H}$ and any $a_k\in\mb{F}_q$ ($-H\leq k\leq J-1$). But this is simply a consequence of the assumption that $h$ is uniformly continuous.
\par Write $\wh{M}=\size{\beta}$, for some integer $M$. Using \eqref{eq:limit repn of integral} with $h=f$, $J=L$ and $H=G-M$, we can rewrite the left side of \eqref{eq:eqn of change of variables lemma} as
\begin{align*}
&\int_{\size{\alpha}<\wh{L}}f(\alpha)\;\mr{d}\alpha\\
=\;&\wh{M}\lim_{G\rightarrow\infty}\wh{G}^{-1}\sum_{\substack{a_k\in\mb{F}_q\\(-G+M\leq k\leq L-1)}}f(a_{-G+M}t^{-G+M}+...+a_{L-1}t^{L-1}).\end{align*}
When $c_j\in\mb{F}_q$ ($j\leq L-M-1$), put
\begin{equation}\label{eq:expansion for w}
w=w(\bs{c};L,M)=\sum_{j=-\infty}^{L-M-1}c_jt^j.\end{equation}
In terms of the truncation operation defined in \eqref{eq:integral part}, we have
\[c_{-G}t^{-G}+...+c_{L-M-1}t^{L-M-1}=\lfloor w\rfloor_{-G}.\]
Put $v(\gamma)=f(\beta\gamma)$. Evidently $v$ is still a continuous function on $\mc{D}_{L-M}$. Using \eqref{eq:limit repn of integral} again with $H=G$, $h=v$ and $\wh{J}=\wh{L}\wh{M}^{-1}$, we can rewrite the right side of \eqref{eq:eqn of change of variables lemma} as
\begin{align*}
&\size{\beta}\int_{\size{\gamma}<\wh{L}\size{\beta}^{-1}}f(\beta\gamma)\;\mr{d}\gamma\\
=\;&\wh{M}\lim_{G\rightarrow\infty}\wh{G}^{-1}\sum_{\substack{c_r\in\mb{F}_q\\(-G\leq r\leq L-M-1)}}v(c_{-G}t^{-G}+...+c_{L-M-1}t^{L-M-1})\\
=\;&\wh{M}\lim_{G\rightarrow\infty}\wh{G}^{-1}\sum_{\substack{c_r\in\mb{F}_q\\(-G\leq r\leq L-M-1)}}f(\beta\lfloor w\rfloor_{-G}).\end{align*}
To prove the lemma, it thus remains to show that, when $L$ and $\beta$ are fixed as in the statement of the lemma, and $G$ is sufficiently large in terms of $L$ and $M$, we have
\begin{align}\label{eq:first reduction in sub II}
&\sum_{\substack{a_k\in\mb{F}_q\\(-G+M\leq k\leq L-1)}}f(a_{-G+M}t^{-G+M}+...+a_{L-1}t^{L-1})\nonumber\\
=\;&\sum_{\substack{c_r\in\mb{F}_q\\(-G\leq r\leq L-M-1)}}f(\beta\lfloor w\rfloor_{-G})+o(\wh{G}).\end{align}
\par To this end, let $G$ be a large integer to be specified. 
We show that we can replace $f(\beta\lfloor w\rfloor_{-G})$ by $f(\lfloor\beta w\rfloor_{-G+M})$, allowing for a negligible error in the process. This allows us to work with the finite expansion $\lfloor\beta w\rfloor_{-G+M}$ in $t$ in place of the infinite expansion $\beta\lfloor w\rfloor_{-G}$. This is valid if we can prove that, for any $\varepsilon>0$, there exists a $G_0=G_0(L,M,\varepsilon)>0$ such that whenever $G\geq G_0$ and $w\in\mc{D}_{L-M}$, we have
\begin{equation}\label{eq:replacement}
|f(\beta\lfloor w\rfloor_{-G})-f(\lfloor\beta w\rfloor_{-G+M})|<\varepsilon.\end{equation}
But this follows from the uniform continuity of $f$. Indeed, note that
\begin{align*}
\beta\lfloor w\rfloor_{-G}-\lfloor\beta w\rfloor_{-G+M}=\;&(\beta w-\beta\left\{w\right\}_{-G})-(\beta w-\left\{\beta w\right\}_{-G+M})\\
=\;&\left\{\beta w\right\}_{-G+M}-\beta\left\{ w\right\}_{-G}.\end{align*}
On recalling that $\size{\beta}=\wh{M}$, we thus have
\[\size{\beta\lfloor w\rfloor_{-G}-\lfloor\beta w\rfloor_{-G+M}}<q^{-G+M}=\wh{G}^{-1}\wh{M},\]
which can be made arbitrarily small by taking $G$ to be large. This together with the uniform continuity of $f$ implies \eqref{eq:replacement}, as a result. Since there are exactly $q^{L-M+G}=\wh{L}\wh{M}^{-1}\wh{G}$ elements $c_r\in\mb{F}_q$ ($-G\leq r\leq L-M-1$), relation \eqref{eq:replacement}, with $\varepsilon$ replaced by $\varepsilon\wh{L}^{-1}\wh{M}$ further implies that when $G$ is large in terms of $L$ and $M$, we have
\begin{equation}\label{eq:first reduction in sub III}
\sum_{\substack{c_r\in\mb{F}_q\\(-G\leq r\leq L-M-1)}}f(\beta\lfloor w\rfloor_{-G})=\sum_{\substack{c_r\in\mb{F}_q\\(-G\leq r\leq L-M-1)}}f(\lfloor\beta w\rfloor_{-G+M})+o(\wh{G}).\end{equation}
To establish \eqref{eq:first reduction in sub II}, it thus suffices to show that, when $L$ and $\beta$ are as in the lemma, and $G$ is sufficiently large in terms of $L$ and $M$, we have
\begin{align}\label{eq:second reduction in sub}
&\sum_{\substack{a_k\in\mb{F}_q\\(-G+M\leq k\leq L-1)}}f(a_{-G+M}t^{-G+M}+...+a_{L-1}t^{L-1})\nonumber\\
=\;&\sum_{\substack{c_r\in\mb{F}_q\\(-G\leq r\leq L-M-1)}}f(\lfloor\beta w\rfloor_{-G+M})+o(\wh{G}).\end{align}
\par To prove \eqref{eq:second reduction in sub}, we compare the coefficients of like powers of $t$ between the expansions for $\lfloor\beta w\rfloor_{-G+M}$ and $a_{-G+M}t^{-G+M}+...+a_{L-1}t^{L-1}$. Since $\size{\beta}=\wh{M}$, we can express $\beta$ as
\[\beta=\sum_{j=-\infty}^{M}b_jt^j,\]
where $b_j\in\mb{F}_q$ ($j\leq M$) with $b_M\neq0$. On recalling \eqref{eq:expansion for w}, a modicum of computation reveals that
\begin{align*}
\Theta_{L-1}(\beta w)=\;&b_Mc_{L-M-1},\\
\Theta_{L-2}(\beta w)=\;&b_Mc_{L-M-2}+b_{M-1}c_{L-M-1},\\
=\;&\vdots\\
\Theta_{-G+M}(\beta w)=\;&b_Mc_{-G}+b_{M-1}c_{-G+1}+...+b_{-G+2M-L+1}c_{L-M-1}.\end{align*}
Since $b_M\neq0$, the change of variables from $a_k$ ($-G+M\leq k\leq L-1$) to $c_r$ ($-G\leq r\leq L-M-1$) via the linear equations $a_k=\Theta_k(\beta w)$ ($-G+M\leq k\leq L-1$) is invertible. In other words, as the $a_k$ run through all elements of $\mb{F}_q$, so do the $c_r$. The relation \eqref{eq:second reduction in sub} is thus an equality without the error term. Putting \eqref{eq:second reduction in sub} into \eqref{eq:first reduction in sub III} yields \eqref{eq:first reduction in sub II}. This completes the proof of the lemma.
\end{proof}
\end{lem}
With this lemma, we are now ready to settle Theorem \ref{thm:change of variables}, in the case where $M$ is a diagonal matrix. 
\begin{lem}\label{lem:diagonal sub}
Theorem \ref{thm:change of variables} holds for any invertible diagonal matrix $M$.
\begin{proof}
Suppose the lemma holds for some positive integer $D$. Now let $\Gamma'$ be the box given by \eqref{eq:box}, and let $\Gamma$ be the same box with $D$ replaced by $D+1$. Let $f:\Gamma\rightarrow\mb{C}$ be continuous, and $M\in GL_{D+1}(\mb{K}_{\infty})$. We express $M$ as
\[M=\begin{pmatrix}
\beta_1&0&\cdots&0&0\\
0&\beta_2&\cdots&0&0\\
\vdots&&\ddots&\vdots&\vdots\\
0&0&\cdots&\beta_D&0\\
0&0&\cdots&0&\beta_{D+1}\end{pmatrix},\]
for some $\beta_1,...,\beta_{D+1}\in\mb{K}_{\infty}\backslash\left\{0\right\}$. Write $M'$ for the $D\times D$-matrix obtained from the top left corner of $M$. Then with $\bs{\alpha}'=(\alpha_1,...,\alpha_D)$ and $\bs{\alpha}=(\alpha_1,...,\alpha_{D+1})$, we have
\[\int_{\Gamma}f(\bs{\alpha})\;\mr{d}\bs{\alpha}=\int_{\Gamma'}\Big(\int_{\size{\alpha_{D+1}}<\wh{R}_{D+1}}f(\bs{\alpha}',\alpha_{D+1})\;\mr{d}\alpha_{D+1}\Big)\mr{d}\bs{\alpha}'.\]
Note that for any fixed $\bs{\alpha}'\in\Gamma'$, the function $\alpha_{D+1}\mapsto f(\bs{\alpha}',\alpha_{D+1})$ is  continuous on $\mc{D}_{R_{D+1}}$. We can thus apply Lemma \ref{lem:change of variables} to perform the substitution $\alpha_{D+1}=\beta_{D+1}\gamma_{D+1}$ in the $\alpha_{D+1}$-integral. This, together with an application of Fubini's theorem, yields
\begin{align*}
&\int_{\Gamma}f(\bs{\alpha})\;\mr{d}\bs{\alpha}\\
=\;&\size{\beta_{D+1}}\int_{\Gamma'}\int_{\size{\gamma_{D+1}}<\wh{R}_{D+1}\size{\beta_{D+1}}^{-1}}f(\bs{\alpha}',\beta_{D+1}\gamma_{D+1})\;\mr{d}\gamma_{D+1}\mr{d}\bs{\alpha}'\\
=\;&\size{\beta_{D+1}}\int_{\size{\gamma_{D+1}}<\wh{R}_{D+1}\size{\beta_{D+1}}^{-1}}\Big(\int_{\Gamma'}f(\bs{\alpha}',\beta_{D+1}\gamma_{D+1})\;\mr{d}\bs{\alpha}'\Big)\mr{d}\gamma_{D+1}.\end{align*}
As before, when $\gamma_{D+1}$ is fixed, the function $\bs{\alpha}'\mapsto f(\bs{\alpha}',\beta_{D+1}\gamma_{D+1})$ is continuous on $\Gamma'$. By the induction hypothesis, we can therefore apply the change of variables $\alpha_i=\beta_i\gamma_i$ ($i=1,...,D$) to the $\bs{\alpha}'$-integral above. With $\bs{\gamma}=(\gamma_1,...,\gamma_{D+1})$ and $\bs{\gamma}'=(\gamma_1,...,\gamma_D)$, this gives
\begin{align*}
&\int_{\Gamma}f(\bs{\alpha})\;\mr{d}\bs{\alpha}\\
=\;&\Big(\prod_{i=1}^{D+1}\size{\beta_i}\Big)\int_{\size{\gamma_{D+1}}<\wh{R}_{D+1}\size{\beta_{D+1}}^{-1}}\int_{M'\bs{\gamma}'\in\Gamma'}f(\beta_1\gamma_1,...,\beta_{D+1}\gamma_{D+1})\;\mr{d}\bs{\gamma}'\mr{d}\gamma_{D+1}.\end{align*}
On recalling the definitions of the boxes $\Gamma$ and $\Gamma'$, together with those of the matrices $M$ and $M'$, we simplify this further as
\[\int_{\Gamma}f(\bs{\alpha})\;\mr{d}\bs{\alpha}
=\size{\det M}\int_{M\bs{\gamma}\in\Gamma}f(M\bs{\gamma})\;\mr{d}\bs{\gamma},\]
as required. By induction, Theorem \ref{thm:change of variables} thus holds for every positive integer $D$.
\end{proof}
\end{lem}
It therefore remains to settle Theorem \ref{thm:change of variables} in the case where $M$ is a unipotent matrix. This is accomplished in the following lemma.
\begin{lem}\label{lem:change of variables II}
Theorem \ref{thm:change of variables} holds for any unipotent matrix $M$ with entries in $\mb{K}_{\infty}$.
\begin{proof}
We prove this by induction on $D$. The case where $D=1$ is trivial. Suppose the lemma is true for some positive integer $D$. Let $\Gamma'$ and $\Gamma$ be as in the proof of the previous lemma. Now let $f:\Gamma\rightarrow\mb{C}$ be  continuous. Take $M\in GL_{D+1}(\mb{K}_{\infty})$. Since Theorem \ref{thm:change of variables} holds for permutation matrices, we can multiply $M$ by permutation matrices without altering either $\size{\det M}$ or the value of the integral on the right side of \eqref{eq:eqn in change of variables}. We can thus assume that $M$ is in the form
\[M=\begin{pmatrix}
1&0&\cdots&0&0\\
m_{2,1}&1&\cdots&0&0\\
\vdots&\vdots&\ddots&\vdots&\vdots\\
m_{D,1}&m_{D,2}&\cdots&1&0\\
m_{D+1,1}&m_{D+1,2}&\cdots&m_{D+1,D}&1,\end{pmatrix}\]
where $m_{i,j}\in\mb{K}_{\infty}$ for $1\leq j<i\leq D+1$. Let $M'$ be the $D\times D$-matrix formed from the entries in the top left hand corner of $M$. Take any $\bs{\gamma}\in\mb{K}_{\infty}^{D+1}$, and write $\bs{\gamma}'=(\gamma_1,...,\gamma_D)$. An elementary computation reveals that the first $D$ coordinates of the $(D+1)$-dimensional vector $M\bs{\gamma}$ form the vector $M'\bs{\gamma}'$. Crucially, this $D$-dimensional vector is independent of $\gamma_{D+1}$. Meanwhile, the last coordinate of $M\bs{\gamma}$ is
\[(M\bs{\gamma})_{D+1}=\zeta+\gamma_{D+1},\]
where $\zeta=\zeta(\gamma_1,...,\gamma_D)=m_{D+1,1}\gamma_1+...+m_{D+1,D}\gamma_D$. Hence we can rewrite the integral on the right side of \eqref{eq:eqn in change of variables} as
\[\int_{M\bs{\gamma}\in\Gamma}f(M\bs{\gamma})\;\mr{d}\bs{\gamma}=\int_{M'\bs{\gamma}'\in\Gamma'}\int_{\size{\zeta+\gamma_{D+1}}<\wh{R}_{D+1}}f(M'\bs{\gamma}',\zeta+\gamma_{D+1})\;\mr{d}\gamma_{D+1}\mr{d}\bs{\gamma}'.\]
Putting $\alpha_{D+1}=\gamma_{D+1}+\zeta$ in the innermost integral, and using Fubini's theorem, we can simplify the above as
\begin{align*}
\int_{M\bs{\gamma}\in\Gamma}f(M\bs{\gamma})\;\mr{d}\bs{\gamma}=\;&\int_{M'\bs{\gamma}'\in\Gamma'}\int_{\size{\alpha_{D+1}}<\wh{R}_{D+1}}f(M'\bs{\gamma}',\alpha_{D+1})\;\mr{d}\alpha_{D+1}\mr{d}\bs{\gamma}'\\
=\;&\int_{\size{\alpha_{D+1}}<\wh{R}_{D+1}}\int_{M'\bs{\gamma}'\in\Gamma'}f(M'\bs{\gamma}',\alpha_{D+1})\;\mr{d}\bs{\gamma}'\mr{d}\alpha_{D+1}.\end{align*}
Put $\bs{\alpha}'=(\alpha_1,...,\alpha_D)$. As before, we see that for fixed $\alpha_{D+1}$, the function $\bs{\alpha}'\mapsto f(\bs{\alpha}',\alpha_{D+1}):\Gamma'\rightarrow\mb{C}$ is continuous. Moreover, from the way in which the matrix $M'$ is chosen, we see that $M'$ is a $D\times D$-unipotent matrix. Using the induction hypothesis on the inner integral above, we have
\begin{align*}
\int_{M\bs{\gamma}\in\Gamma}f(M\bs{\gamma})\;\mr{d}\bs{\gamma}=\size{\det M}^{-1}\;&\int_{\size{\alpha_{D+1}}<\wh{R}_{D+1}}\int_{\Gamma'}f(\bs{\alpha}',\alpha_{D+1})\;\mr{d}\bs{\alpha}'\mr{d}\alpha_{D+1}.\end{align*}
On recalling the definitions of the boxes $\Gamma$ and $\Gamma'$, and noting that $\det M'=\det M$, we obtain
\[\size{\det M}\int_{M\bs{\gamma}\in\Gamma}f(M\bs{\gamma})\;\mr{d}\bs{\gamma}=\int_{\Gamma}f(\bs{\alpha})\;\mr{d}\bs{\alpha},\]
as required. This completes the proof of the inductive step, and hence that of the lemma.
\end{proof}
\end{lem}
Theorem \ref{thm:change of variables}, in its full generality, now follows on combining Lemmata \ref{lem:diagonal sub} and \ref{lem:change of variables II} together with the trivial case where $M$ is a permutation matrix.
\par A couple of remarks are in order. First, in the proof of Theorems \ref{thm:asymptotic formula} and \ref{thm:weak approx}, we only use Theorem \ref{thm:change of variables} in the simple case where $f$ is defined in terms of the exponential given in \eqref{defn:e}, and $M=\beta I_D$, for various $\beta\in\mb{K}_{\infty}\backslash\left\{0\right\}$ and positive integers $D$. Second, we have so far only proved the change-of-variable property when the substitutions are linear. It is natural to ask whether non-linear substitutions are permissible. Unfortunately we do not have this luxury in function fields as we do over $\mb{R}$. As an illustration, suppose we wish to make the substitution $\alpha=\gamma^N$, for some positive integer $N$. Then $\ord\alpha$ must be divisible by $N$. So those $\alpha$ with $\ord\alpha$ not divisible by $N$ certainly cannot be substituted in this manner. However, for the purposes of this paper, it suffices to content ourselves with  linear substitutions only.
\section{Geometry of numbers over function fields}\label{sec:geometry of numbers}
In this section, we establish function field analogues of well-known results from the geometry of numbers. It is worth mentioning that Mahler \cite{mahler} has proved the analogues of Minkowski's theorems over any field with a non-Archimedean valuation. But we deem it worthwhile to prove these results over the field $\mb{K}_{\infty}$ in particular, so as to taylor them for the specific purposes of this paper. The approach taken here is adapted from that of Chapter 12 of \cite{davenport}.
\par We define a \emph{lattice} in $\mb{K}_{\infty}^D$ to be the set of points in the form $\bs{x}=\Lambda\bs{u}$, where $\Lambda$ is a $D\times D$-matrix over $\mb{K}_{\infty}$, and $\bs{u}$ is an integral $D$-tuple. By a slight abuse of notation, we also denote the set of such $\bs{x}$ by $\Lambda$. It is apparent that every lattice is an additive subgroup of $\mb{K}_{\infty}^D$, and has an \emph{integral basis} $\bs{x}^{(1)},...,\bs{x}^{(D)}$, given by the columns of the matrix $\Lambda$. Under this basis, every point in this lattice can be expressed as a linear combination of the $\bs{x}^{(\nu)}$ ($1\leq \nu\leq D$) over $\mb{A}$. We define the \emph{determinant} $d(\Lambda)$ of the lattice $\Lambda$ to be $\size{\det\Lambda}$. Here the $\Lambda$ in the latter notation refers to the underlying matrix of the lattice. It is easily seen that the determinant of a lattice is independent of the choice of integral basis. Indeed, let $\left\{\bs{v}^{(\nu)}:\nu=1,...,D\right\}$ and $\left\{\bs{w}^{(\nu)}:\nu=1,...,D\right\}$ be two integral bases of $\Lambda$. Write $V=(\bs{v}^{(1)},...,\bs{v}^{(D)})$ and $W=(\bs{w}^{(1)},...,\bs{w}^{(D)})$. We can then find matrices $Z,Z'\in GL_D(\mb{A})$ such that $V=ZW$ and $W=Z'V$. Thus $V=ZZ'V$, which in turn yields $ZZ'=I_D$. It follows that $\det Z$ is a unit in $\mb{A}$, which has absolute value $1$ in $\mb{K}_{\infty}$. This shows that the determinant of a lattice as defined above is independent of the choice of integral basis. Meanwhile, it follows from the definition \eqref{eq:size of tuple} of the absolute value on $\mb{K}_{\infty}^D$ that every $D$-dimensional box, in the sense of Appendix A, is an additive subgroup of $\mb{K}_{\infty}^D$. 
\par We first need a lemma concerning the volume of any parallelepiped in the shape
\[
\mc{P}(\bs{x}^{(1)},...,\bs{x}^{(D)})=\left\{\alpha_1\bs{x}^{(1)}+...+\alpha_D\bs{x}^{(D)}:\alpha_1,...,\alpha_D\in\mb{T}\right\},\]
for any linearly independent vectors $\bs{x}^{(1)},...,\bs{x}^{(D)}\in\mb{K}_{\infty}^D$.
\begin{lem}\label{lem:parallelopiped}
Let $D$ be any positive integer. Let $\bs{x}^{(1)},...,\bs{x}^{(D)}$ be real $D$-tuples that are linearly independent over $\mb{K}_{\infty}$. Then
\[\mr{Vol}\;\mc{P}(\bs{x}^{(1)},...,\bs{x}^{(D)})=\size{\det(\bs{x}^{(1)},...,\bs{x}^{(D)})}.\]
\begin{proof}
The claim here is a straightforward consequence of Theorem \ref{thm:change of variables}. Write $M$ for the $D\times D$-matrix with columns $\bs{x}^{(1)},...,\bs{x}^{(D)}$. Since $\bs{x}^{(1)},...,\bs{x}^{(D)}$ are linearly independent, it follows that $M\in GL_D(\mb{K}_{\infty})$. So every element of $\mc{P}=\mc{P}(\bs{x}^{(1)},...,\bs{x}^{(D)})$ is in the form $\bs{\beta}=M\bs{\alpha}$, where $\bs{\alpha}\in\mb{T}^D$. Using Theorem \ref{thm:change of variables}, we can apply this substitution in the integral $\int_{\mb{T}^D}\mr{d}\bs{\alpha}$, which equals 1 by normalisation. Hence
\[1=\int_{\mb{T}^D}\mr{d}\bs{\alpha}=\size{\det M}^{-1}\int_{\mc{P}}\mr{d}\bs{\beta}=\size{\det M}^{-1}\vol\mc{P}.\]
The required equality thus follows by rearranging this equality, and recalling the definition of $M$.
\end{proof}
\end{lem}
In the remainder of this appendix, we fix a lattice $\Lambda$ in $\mb{K}_{\infty}^D$. Using this lemma, we derive the following analogue of Minkowski's first theorem over function fields.
\begin{lem}\label{lem:minkowski I}
Suppose $\Gamma$ is an $D$-dimensional box with
\begin{equation}\label{eq:assumption in minkowski I}
\vol{\Gamma}>d(\Lambda).\end{equation}
Then $\Gamma$ contains a non-zero point of $\Lambda$.
\begin{proof}
Let $\bs{x}^{(1)},...,\bs{x}^{(D)}$ denote an integral basis for $\Lambda$. Put $\mc{P}=\mc{P}(\bs{x}^{(1)},...,\bs{x}^{(D)})$. For each $\bs{x}\in\Lambda$, let
$\mc{R}(\bs{x})=\mc{P}\cap(\Gamma-\bs{x})$. By the translation-invariance of the measure on $\mb{K}_{\infty}^D$, we have
\[\vol\mc{R}(\bs{x})=\vol((\mc{P}+\bs{x})\cap\Gamma)\]
for each such $\bs{x}$. Summing over all such $\bs{x}$ gives
\[
\sum_{\bs{x}\in\Lambda}\vol\mc{R}(\bs{x})=\sum_{\bs{x}\in\Lambda}\int_{\mc{P}+\bs{x}}1_{\Gamma}(\bs{\alpha})\;\mr{d}\bs{\alpha}.\]
Here $1_{\Gamma}$ denotes the indicator function on $\Gamma$. The compactness of $\Gamma$ implies that, on both sides of the last equation, only finitely many non-zero terms occur, so there is no issue with the  convergence of the two infinite sums involved. Also, from the definition of $\mc{P}$, it transpires that the translates $\mc{P}+\bs{x}$ of $\mc{P}$ by lattice points $\bs{x}\in\Lambda$ are disjoint, and their union coincides with $\mb{K}_{\infty}^D$. It thus follows from the last displayed equality, together with \eqref{eq:assumption in minkowski I} and Lemma \ref{lem:parallelopiped}, that
\[\sum_{\bs{x}\in\Lambda}\vol\mc{R}(\bs{x})
=\int_{\mb{K}_{\infty}^s}1_{\Gamma}(\bs{\alpha})\;\mr{d}\bs{\alpha}
=\vol\Gamma>d(\Lambda)=\vol\mc{P}.\]
Since each $\mc{R}(\bs{x})$ is a subset of $\mc{P}$, it follows from this inequality that two of the $\mc{R}(\bs{x})$ must intersect. There thus exist distinct lattice points $\bs{x}_1,\bs{x}_2\in\Lambda$ and $\bs{\alpha}\in\mc{P}$ for which $\bs{\alpha}\in\mc{R}(\bs{x}_1)\cap\mc{R}(\bs{x}_2)$. Tthe definition of $\mc{R}(\bs{x})$ gives $\bs{\alpha}+\bs{x}_1\in\Gamma$ and $\bs{\alpha}+\bs{x}_2\in\Gamma$. Since $\Gamma$ is an additive group, we must have
\[\bs{x}_1-\bs{x}_2=(\bs{\alpha}+\bs{x}_1)-(\bs{\alpha}+\bs{x}_2)\in\Gamma.\]
This point must be non-zero since $\bs{x}_1\neq\bs{x}_2$. Since $\bs{x}_1,\bs{x}_2\in\Lambda$, we also have  $\bs{x}_1-\bs{x}_2\in\Lambda$, as required.
\end{proof}
\end{lem}
Before advancing further, we need to introduce some notation and terminology. When $\mu$ is a positive integer, and $\bs{y}^{(1)},...,\bs{y}^{(\mu)}\in\mb{K}_{\infty}^D$, let $\size{\bs{y}^{(1)},...,\bs{y}^{(\mu)}}_{\mb{K}}$ denote the set of all linear combinations of $\bs{y}^{(1)},...,\bs{y}^{(\mu)}$ over $\mb{K}$. Meanwhile, define $\wh{R}_1$ to be the infimum of $\size{\bs{x}}$ over all non-zero points $\bs{x}\in\Lambda$. Take $\bs{x}^{(1)}$ to be some point in $\Lambda$ with $\size{\bs{x}^{(1)}}=\wh{R}_1$. Inductively, for every $\nu=2,...,D$, define $\wh{R}_{\nu}$ to be the infimum of $\size{\bs{x}}$ over all points $\bs{x}\in\Lambda\backslash\size{\bs{x}^{(1)},...,\bs{x}^{(\nu-1)}}_{\mb{K}}$. Take $\bs{x}^{(\nu)}\in\Lambda$ with $\size{\bs{x}^{(\nu)}}=\wh{R}_{\nu}$. We call $\wh{R}_{\nu}$ the $\nu$-th \emph{successive minima} of $\Lambda$, and $\bs{x}^{(\nu)}$ a $\nu$-th minimal point of $\Lambda$. Evidently $\wh{R}_1\leq...\leq\wh{R}_D$.
\begin{lem}\label{lem:minkowski II}
Let $\wh{R}_1,...,\wh{R}_D$ be the successive minima of $\Lambda$. Then
\[\wh{R}_1...\wh{R}_D=d(\Lambda).\]
\begin{proof}
By applying a suitable unimodular linear transformation if necessary, we can assume without loss of generality that $\Lambda$ has minimal points in the form
\begin{align*}
\bs{x}^{(1)}=\;&(x_1^{(1)},0,...,0),\\
\bs{x}^{(2)}=\;&(x_1^{(2)},x_2^{(2)},...,0)\\
\vdots&\\
\bs{x}^{(D)}=\;&(x_1^{(D)},x_2^{(D)},...,x_D^{(D)}),
\end{align*}
for some $x_{\mu}^{(\nu)}\in\mb{K}_{\infty}$ ($1\leq\mu\leq\nu\leq D$). Choose an arbitrary integral basis $\bs{y}^{(1)},...,\bs{y}^{(D)}$ of $\Lambda$. Then we can find a matrix $Z\in GL_D(\mb{A})$ such that
\[(\bs{x}^{(1)},...,\bs{x}^{(D)})=Z(\bs{y}^{(1)},...,\bs{y}^{(D)}).\]
Taking absolute values of the determinants of both sides above, and recalling that $d(\Lambda)=\size{\det(\bs{y}^{(1)},...,\bs{y}^{(D)})}$, we get
\[\size{\det(\bs{x}^{(1)},...,\bs{x}^{(D)})}=\size{\det Z}d(\Lambda).\]
Since $\det Z$ is a non-zero element of $\mb{A}$, it follows that $\size{\det Z}\geq1$. This together with our initial assumption on the points $\bs{x}^{(1)},...,\bs{x}^{(D)}$ gives
\[d(\Lambda)\leq\size{x_1^{(1)}...x_D^{(D)}}.\]
Courtesy of \eqref{eq:size of tuple}, the above implies that
\[d(\Lambda)\leq\size{\bs{x}^{(1)}}...\size{\bs{x}^{(D)}}.\]
From the definition of minimal points, this gives
\[d(\Lambda)\leq\wh{R}_1...\wh{R}_D.\]
It thus remains to establish the reverse inequality.
\par To this end, let $\Gamma$ be the box given by the inequalities \eqref{eq:box}. We claim that $\Gamma$ contains no non-zero point of $\Lambda$. Indeed, suppose $\Gamma$ contains a non-zero point $\bs{x}$ of $\Lambda$. From the definition of minimal points, we can find a $\nu=1,...,D$ such that $\bs{x}$ is linearly dependent over $\mb{K}$ on $\bs{x}^{(1)},...,\bs{x}^{(\nu)}$, but not on $\bs{x}^{(1)},...,\bs{x}^{(\nu-1)}$. So on one hand, it follows from the definition of $\wh{R}_{\nu}$ that $\size{\bs{x}}\geq\wh{R}_{\nu}$. On the other hand, by the way in which we have chosen the minimal points $\bs{x}^{(1)},...,\bs{x}^{(\nu)}$, we deduce that $x_{\nu+1}=...=x_D=0$. Also, since $\bs{x}\in\Gamma$, we have $\size{x_u}<\wh{R}_u$ for $u=1,...,\nu$. On recalling that $\wh{R}_1\leq...\leq\wh{R}_{\nu}$, we see that $\size{\bs{x}}<\wh{R}_{\nu}$. This is a contradiction to what we mentioned above. It thus follows from \eqref{eq:box} that $\bs{x}\notin\Gamma$, which proves the claim. This claim, together with Lemma \ref{lem:minkowski I}, implies that $\vol{\Gamma}\leq d(\Lambda)$. But it follows easily from \eqref{eq:extended measure} that $\vol\Gamma=\wh{R}_1...\wh{R}_D$, so we have the desired inequality
\[d(\Lambda)\geq\wh{R}_1...\wh{R}_D.\]
\end{proof}
\end{lem}
Although the choice of minimal points for a lattice is not unique, the following variant of Lemma 12.3 of \cite{davenport} gives an essentially canonical choice that is convenient to work with.
\begin{lem}\label{lem:normalisation}
Let $\Lambda$ and $\wh{R}_1,...,\wh{R}_D$ be as in Lemma \ref{lem:minkowski II}. Then $\Lambda$ has an integral basis $\bs{X}^{(\nu)}$ ($1\leq\nu\leq D$) in the form
\begin{equation}\label{eq:normalisation}
\bs{X}^{(\nu)}=(X_1^{(\nu)},...,X_{\nu}^{(\nu)},0,...,0)\qquad(\nu=1,...,D)\end{equation}
such that for each such $\nu$, we have
\begin{equation}\label{eq:size of normalised min pts}
\size{\bs{X}^{(\nu)}}=\wh{R}_{\nu}\qquad\text{and}\qquad\size{X_{\nu}^{(\nu)}}=\wh{R}_{\nu}.\end{equation}
\begin{proof}
Let $\bs{x}^{(1)},...,\bs{x}^{(D)}$ be chosen as in the proof of Lemma \ref{lem:minkowski II}. Take $\bs{X}^{(1)}=\bs{x}^{(1)}$, and choose $\bs{X}^{(2)}\in\Lambda\cap\size{\bs{x}^{(1)},\bs{x}^{(2)}}_{\mb{K}}$, such that $\bs{X}^{(1)}$ and $\bs{X}^{(2)}$ generate integrally all the points in $\Lambda\cap\size{\bs{x}^{(1)},\bs{x}^{(2)}}_{\mb{K}}$. The choice of $\bs{X}^{(2)}$ is arbitrary to the extent of added integral multiples of $\bs{x}^{(1)}$. Since $\bs{x}^{(1)}$ and $\bs{x}^{(2)}$ generate rationally (not necessarily integrally) all the points of $\Lambda$, in particular we can find $m,u_1,u_2\in\mb{A}$ with $m\neq0$ such that
\begin{equation}\label{eq:mX(2)}
m\bs{X}^{(2)}=u_1\bs{x}^{(1)}+u_2\bs{x}^{(2)}.\end{equation}
\par Our aim now is to show that the polynomials $u_1$ and $u_2$ above can be assumed to have absolute values at most $\size{m}$. This then implies that
\begin{align*}
\size{\bs{X}^{(2)}}\leq\;&\ma{\size{u_1/m}\size{\bs{x}^{(1)}},\size{u_2/m}\size{\bs{x}^{(2)}}}\\
\leq\;&\ma{\wh{R}_1,\wh{R}_2}=\wh{R}_2.\end{align*}
To prove our claim, note that $\bs{x}^{(2)}$ is generated integrally by $\bs{X}^{(1)}=\bs{x}^{(1)}$ and $\bs{X}^{(2)}$, so we can find $v_1,v_2\in\mb{A}$ such that
\[\bs{x}^{(2)}=v_1\bs{x}^{(1)}+v_2\bs{X}^{(2)}.\]
Combining this with \eqref{eq:mX(2)} gives
\begin{equation}\label{eq:X(2) in terms of x(1)}
(m-u_2v_2)\bs{X}^{(2)}=(u_1+u_2v_1)\bs{x}^{(1)}.\end{equation}
Due to our choice of $\bs{x}^{(1)}$ in Lemma \ref{lem:minkowski II}, by comparing the second coordinates on both sides above, we get
\begin{equation}\label{eq:m=u2v2}
m=u_2v_2.\end{equation}
Meanwhile, on considering the first coordinates in \eqref{eq:X(2) in terms of x(1)}, we see that $u_1=-u_2v_1$. Putting this and \eqref{eq:m=u2v2} into \eqref{eq:mX(2)} yields
\begin{equation}\label{eq:u2v2X(2)}
u_2v_2\bs{X}^{(2)}=u_2(-v_1\bs{x}^{(1)}+\bs{x}^{(2)}).\end{equation}
If $u_2=0$, then by looking at the second coordinates on both sides of \eqref{eq:mX(2)}, we get $m=0$, a contradiction. Hence $u_2\neq0$, and \eqref{eq:u2v2X(2)} gives
\begin{equation}\label{eq:v2X(2)}
v_2\bs{X}^{(2)}=-v_1\bs{x}^{(1)}+\bs{x}^{(2)}.\end{equation}
Arguing as before, we have $v_2\neq0$. Write $-v_1=w_1v_2+r_1$, where $w_1,r_1\in\mb{A}$ with $\size{r_1}<\size{v_2}$. Using this in \eqref{eq:v2X(2)}, we have
\begin{equation}\label{eq:new v2X(2)}
v_2\bs{X}^{(2)}=r_1\bs{x}^{(1)}+\bs{x}^{(2)}.\end{equation}
Here we have redefined $\bs{X}^{(2)}$ by subtracting  $w_1\bs{x}^{(1)}$ from the old $\bs{X}^{(2)}$. On comparing \eqref{eq:new v2X(2)} with \eqref{eq:mX(2)}, our claim in this paragraph is thus verified.
\par In the same vein, for each $\nu=3,...,D$, we define $\bs{X}^{(\nu)}$ to be a point in $\Lambda\cap\size{\bs{x}^{(1)},...,\bs{x}^{(\nu)}}_{\mb{K}}$ such that $\bs{X}^{(1)},...,\bs{X}^{(\nu)}$ generate integrally all points in $\Lambda\cap\size{\bs{x}^{(1)},...,\bs{x}^{(\nu)}}_{\mb{K}}$, with $\bs{X}^{(\nu)}$ arbitrary with respect to added integral multiples of $\bs{x}^{(1)},...,\bs{x}^{(\nu-1)}$. A similar argument to the one used in the last paragraph shows that
\[\size{\bs{X}^{(\nu)}}\leq\wh{R}_{\nu}\qquad\text{for all }\nu=1,...,D,\]
and consequently
\[\size{X_{\nu}^{(\nu)}}\leq\size{\bs{X}^{(\nu)}}\leq\wh{R}_{\nu}\]
for all such $\nu$. Our choice of the $\bs{X}^{(\nu)}$ ($1\leq\nu\leq s$) indicates that these points are in the form given in \eqref{eq:size of normalised min pts}.
\par The lower bound $\size{X_{\nu}^{(\nu)}}\geq\wh{R}_{\nu}$ follows by comparing determinants. Our choice of $\bs{X}^{(1)},...,\bs{X}^{(D)}$ above gives rise to an integral basis of $\Lambda$, so
\[d(\Lambda)=\size{\det(\bs{X}^{(1)},...,\bs{X}^{(D)})}=\size{X_1^{(1)}...X_D^{(D)}}.\]
Recalling the definition of successive minima, we have
\[d(\Lambda)\leq\wh{R}_1...\wh{R}_{\nu-1}\size{X_{\nu}^{(\nu)}}\wh{R}_{\nu+1}...\wh{R}_D.\]
Applying Lemma \ref{lem:minkowski II}, we obtain
\[\size{X_{\nu}^{(\nu)}}\geq\wh{R}_{\nu}\qquad\text{for all }\nu=1,...,D.\]
Hence $\size{\bs{X}^{(\nu)}}\geq\size{X_{\nu}^{(\nu)}}\geq\wh{R}_{\nu}$, as required.
\end{proof}
\end{lem}
The following lemma counts exactly the number of lattice points lying within a given distance from the origin.
\begin{lem}\label{lem:lattice points in a cube}
Let $\Lambda$ and $\wh{R}_1,...,\wh{R}_D$ be as in Lemma \ref{lem:minkowski II}. Then
\[\#\left\{\bs{x}\in\Lambda:\size{\bs{x}}<\wh{U}\right\}=\begin{cases}
1,& \text{if }U<R_1,\\
\wh{U}^{\nu}(\wh{R}_1...\wh{R}_{\nu})^{-1},& \text{if }R_{\nu}\leq U<R_{\nu+1}.\end{cases}\]
\begin{proof}
When $U<R_1$,the only point in $\Lambda$ with absolute value less than $\wh{U}$ is $0$, so the lemma follows immediately in this case. Now suppose $R_{\nu}\leq U<R_{\nu+1}$. By Lemma \ref{lem:normalisation}, every $\bs{x}\in\Lambda$ with $\size{\bs{x}}<\wh{U}$ has the form
\begin{equation}\label{eq:nu-dimensional x}
\bs{x}=v_1\bs{X}^{(1)}+...+v_{\nu}\bs{X}^{(\nu)},\end{equation}
where $v_1,...,v_{\nu}\in\mb{A}$. The normalisation \eqref{eq:normalisation} in that lemma implies further that
\begin{align}
\size{x_1}=\;&\size{v_1X_1^{(1)}&+v_2X_1^{(2)}&+...&+v_{\nu}X_1^{(\nu)}}&<\wh{U},\nonumber\\
\size{x_2}=\;&&\size{v_2X_2^{(2)}&+...&+v_{\nu}X_2^{(\nu)}}&<\wh{U},\nonumber\\
&&&\vdots&&\nonumber\\
\size{x_{\nu}}=\;&&&&\size{v_{\nu}X_{\nu}^{(\nu)}}&<\wh{U}.\label{eq:choice for v_nu}\end{align}
The last inequality above together with \eqref{eq:size of normalised min pts} yields
\[\size{v_{\nu}}<\wh{U}\wh{R}_{\nu}^{-1}.\]
Inductively, for any $\sigma=2,...,\nu$, once we obtain
\[\size{v_{\mu}}<\wh{U}\wh{R}_{\mu}^{-1}\]
for all $\mu=\sigma,...,\nu$, we know from Lemma \ref{lem:normalisation} that
\[\size{v_{\mu}X_{\sigma-1}^{(\mu)}}<\wh{U}\wh{R}_{\mu}^{-1}\size{\bs{X}^{(\mu)}}=\wh{U}\]
for all such $\mu$. The $(\sigma-1)$-th inequality in \eqref{eq:choice for v_nu}, namely
\[\size{x_{\sigma-1}}=\size{v_{\sigma-1}X_{\sigma-1}^{(\sigma-1)}+v_{\sigma}X_{\sigma-1}^{(\sigma)}+...+v_{\nu}X_{\sigma-1}^{(\nu)}}<\wh{U},\]
then produces the upper bound
\[\size{v_{\sigma-1}X_{\sigma-1}^{(\sigma-1)}}<\wh{U}.\]
Applying Lemma \ref{lem:normalisation} again gives
\[\size{v_{\sigma-1}}<\wh{U}\wh{R}_{\sigma-1}^{-1}.\]
Hence $v_1,...,v_{\nu}$ are subject to the constraints
\begin{equation}\label{eq:constraints for v}
\size{v_{\mu}}<\wh{U}\wh{R}_{\mu}^{-1}\qquad\text{for }\mu=1,...,\nu.\end{equation}
Therefore
\[\#\left\{\bs{x}\in\Lambda:\size{\bs{x}}<\wh{U}\right\}\leq\prod_{\mu=1}^{\nu}(\wh{U}/\wh{R}_{\mu})=\wh{U}^{\nu}/(\wh{R}_1...\wh{R}_{\nu}).\]
The reverse inequality can be obtained by counting only those points $\bs{x}$ in the form \eqref{eq:nu-dimensional x} which satisfy \eqref{eq:constraints for v}.
\end{proof}
\end{lem}
Before we proceed further, we introduce the notion of adjoint lattices. We say that two lattices $\Lambda$ and $M$ in $\mb{K}_{\infty}^D$ are adjoint if their underlying matrices satisfy
\begin{equation}\label{eq:adjoint lattices}
\Lambda^{T}M=I_D,\end{equation}
where $\Lambda^T$ denotes the transpose of the matrix $\Lambda$. The points on $\Lambda$ and $M$ have the respective forms $\bs{x}=\Lambda\bs{u}$ and $\bs{y}=M\bs{v}$, where $\bs{u},\bs{v}\in\mb{A}^D$. Equation \eqref{eq:adjoint lattices} then implies that
\begin{equation}\label{eq:perpendicularity}
\bs{x}\cdot\bs{y}=\bs{u}\cdot\bs{v}\in\mb{A}.\end{equation}
The following lemma gives an inverse relation between the successive minima of two adjoint lattices.
\begin{lem}\label{lem:adjoint lattices}
Let $\Lambda$ and $M$ be adjoint lattices in $\mb{K}_{\infty}^D$, with successive minima $\wh{R}_1,...,\wh{R}_D$ and $\wh{S}_1,...,\wh{S}_D$. Then for each $\nu=1,...,D$, we have
\[\wh{R}_{\nu}=\wh{S}_{D-\nu+1}^{-1}.\]
\begin{proof}
Let $\bs{x}^{(1)},...,\bs{x}^{(D)}$ and $\bs{y}^{(1)},...,\bs{y}^{(D)}$ denote respectively minimal points of $\Lambda$ and $M$, with $\size{\bs{x}^{(\mu)}}=\wh{R}_{\mu}$ and $\size{\bs{y}^{(\mu)}}=\wh{S}_{\mu}$ for any $\mu=1,...,D$. Fix $\nu=1,...,D$, and let
\[\mc{U}=\left\{\bs{y}\in M:\bs{x}^{(\mu)}\cdot\bs{y}=0\text{ for all }\mu=1,...,\nu\right\}.\]
This is an $(D-\nu)$-dimensional subspace of $\mb{K}_{\infty}^D$. So not all of $\bs{y}^{(1)},...,\bs{y}^{(D-\nu+1)}$ lie in $\mc{U}$. There thus exist positive integers $m$ and $n$ with $m\leq\nu$ and $n\leq D-\nu+1$ such that $\bs{x}^{(m)}\cdot\bs{y}^{(n)}\neq0$. Owing to \eqref{eq:perpendicularity}, this implies that
\[\size{\bs{x}^{(m)}\cdot\bs{y}^{(n)}}\geq1.\]
But for each $\mu=1,...,D$ we have
\[\size{x_{\mu}^{(m)}}\leq\size{\bs{x}^{(m)}}=\wh{R}_m\leq\wh{R}_{\nu}\]
and
\[\size{y_{\mu}^{(n)}}\leq\size{\bs{y}^{(n)}}=\wh{S}_n\leq\wh{S}_{D-\nu+1}.\]
Therefore
\begin{equation}\label{eq:lower bound for adjoints}
\wh{R}_{\nu}\wh{S}_{D-\nu+1}\geq1\qquad\text{for all }\nu=1,...,D.\end{equation}
The reverse inequality can be proved as follows. From \eqref{eq:adjoint lattices}, we have
\[d(\Lambda)d(M)=1.\]
Applying Lemma \ref{lem:minkowski II} on the lattices $\Lambda$ and $M$ in turn gives
\[(\wh{R}_1...\wh{R}_D)(\wh{S}_1...\wh{S}_D)=1.\]
Rearranging this yields
\[\wh{R}_{\nu}\wh{S}_{D-\nu+1}=\prod_{\substack{\mu=1\\ \mu\neq\nu}}^{D}(\wh{R}_{\mu}\wh{S}_{D-\mu+1})^{-1}.\]
Applying \eqref{eq:lower bound for adjoints} leads to
\[\wh{R}_{\nu}\wh{S}_{D-\nu+1}\leq1\qquad\text{for all }\nu=1,...,D,\]
as required.
\end{proof}
\end{lem}

\end{document}